\documentclass[12pt,a4paper]{article}
\usepackage[utf8]{inputenc}
\usepackage{enumerate,booktabs,comment}
\usepackage{hyperref} 	
\usepackage{amssymb}
\usepackage{amsmath}
\usepackage{amsthm}
\usepackage{authblk}

\usepackage{mathrsfs} 	
\usepackage{graphicx} 	
\usepackage[margin=2.5cm]{geometry}
\usepackage{setspace}
\usepackage{adjustbox}
\usepackage{appendix}
\usepackage{mathabx}	
\usepackage{easybmat}	
\usepackage{booktabs}
\usepackage{listings}
\usepackage{color}
\definecolor{dkgreen}{rgb}{0,0.6,0}
\definecolor{gray}{rgb}{0.5,0.5,0.5}
\definecolor{mauve}{rgb}{0.58,0,0.82}
\renewcommand{\epsilon}{\varepsilon}
\renewcommand{\phi}{\varphi}
\makeatletter
\renewcommand\section{\@startsection
  {section}{2}{0mm}
  {-\baselineskip}
  {0.5\baselineskip}
  {\bfseries\large}}
\makeatother
\makeatletter
\renewcommand\subsection{\@startsection
  {subsection}{2}{0mm}
  {-\baselineskip}
  {0.5\baselineskip}
  {\bfseries\normalsize}}
\makeatother
\newtheoremstyle{plain}
  {0.7cm}			
  {0.7cm}			
  {\normalfont} 
  {0pt}       	
  {\bfseries} 	
  {\\}         	
  {0pt}
  {}          
\newtheorem{theorem}{Theorem}[section]
\newtheorem{lemma}[theorem]{Lemma}
\newtheorem{proposition}[theorem]{Proposition}
\newtheorem{corollary}[theorem]{Corollary}
\newtheorem{conjecture}[theorem]{Conjecture}
\newtheorem{remark}[theorem]{Remark}
\newtheorem{notation}[theorem]{Notation}

\newtheorem{definition}[theorem]{Definition}

\newtheorem{constr}[theorem]{Construction}

\DeclareMathOperator*{\code}{\Gamma}

\DeclareMathOperator*{\F}{\mathbb{F}}

\DeclareMathOperator*{\calq}{\mathcal{Q}}
\DeclareMathOperator*{\calc}{\mathcal{C}}

\DeclareMathOperator*{\aut}{Aut}
\DeclareMathOperator*{\Aut}{Aut}

\DeclareMathOperator*{\sym}{Sym}
\DeclareMathOperator*{\isom}{Isom}
\DeclareMathOperator*{\Tr}{Tr}
\DeclareMathOperator*{\tr}{Tr}

\DeclareMathOperator*{\sing}{sing}

\DeclareMathOperator*{\gl}{GL}
\DeclareMathOperator*{\gaml}{\mathrm{\Gamma L}}

\DeclareMathOperator*{\syp}{Sp}
\DeclareMathOperator*{\oo}{GO}

\DeclareMathOperator*{\gamsp}{\mathrm{\Gamma Sp}}
\DeclareMathOperator*{\gamoo}{\mathrm{\Gamma O}}

\newcommand{\psl}{\mathrm{PSL}}

\newcommand{\pgaml}{\mathrm{P\Gamma L}}
\newcommand{\psigmal}{\mathrm{P\Sigma L}}

\newcommand{\agl}{\mathrm{AGL}}

\newcommand{\agaml}{\mathrm{A\Gamma L}}

\newcommand{\blocky}{\Sigma}
\newcommand{\qbinom}[2]{\genfrac{[}{]}{0pt}{}{#1}{#2}}

\allowdisplaybreaks

\numberwithin{theorem}{section}
\numberwithin{equation}{section}
\numberwithin{figure}{section}
\numberwithin{table}{section}

\title{Codes and Designs in Johnson Graphs From Symplectic Actions on Quadratic Forms}
\author{John Bamberg, Alice Devillers, Mark Ioppolo, Cheryl E. Praeger}
\affil{Centre for the Mathematics of Symmetry and Computation,\\
Department of Mathematics and Statistics,\\
The University of Western Australia.}

\date{}

\begin{document}
\maketitle
\begin{abstract}
The Johnson graph $J(v, k)$ has as vertices the $k$-subsets of $\mathcal{V}=\{1,\ldots, v\}$,
and two vertices are joined by an edge if their intersection has size $k-1$.
An  \emph{$X$-strongly incidence-transitive code} in $J (v, k)$ is a proper vertex subset $\Gamma$ such that the subgroup $X$ of graph
automorphisms leaving $\Gamma$ invariant is transitive on the set $\Gamma$ of `codewords', and for each codeword $\Delta$, the setwise stabiliser $X_\Delta$ is transitive on $\Delta \times (\mathcal{V}\setminus \Delta)$. We classify the \emph{$X$-strongly incidence-transitive codes} in $J(v,k)$ for which $X$ is the symplectic group $\syp_{2n}(2)$ acting as a $2$-transitive permutation group of degree $2^{2n-1}\pm 2^{n-1}$, where the stabiliser $X_\Delta$ of a codeword $\Delta$ is contained in a \emph{geometric} maximal subgroup of $X$. In particular, we construct two new infinite families of strongly incidence-transitive codes associated with the reducible maximal subgroups of $\syp_{2n}(2)$.

\medskip\noindent{\it Keywords:}\quad codes in graphs; Johnson graphs; finite symplectic groups; Aschbacher classification; strongly incidence-transitive.

\medskip\noindent{\it MSN Classification:}\quad  05C25; 05E18; 20B25; 94B25
\end{abstract}

\section{Introduction}\label{intro}

Let $\mathcal{V}$ be a finite set of cardinality $v > 1$ and let $k$ be an
integer which satisfies $1 < k < v-1$. We denote by $\binom{\mathcal{V}}{k}$
the set of all $k$-sets ($k$-element subsets) of $\mathcal{V}$. The \emph{Johnson graph}
$J(\mathcal{V},k)$ is a graph with vertex set $\binom{\mathcal{V}}{k}$, where
a vertex pair $\Delta_1,\Delta_2 \in \binom{\mathcal{V}}{k}$ is joined by an
edge if and only if $|\Delta_1 \cap \Delta_2| = k-1$. If $\Delta$ is a $k$-set then we write $\overline{\Delta} = \mathcal{V}\setminus\Delta$, the complement of $\Delta$ in $\mathcal{V}$. The {\em complementing map} $s$ which swaps each $k$-set with its complement defines a graph isomorphism from $J(\mathcal{V},k)$ to $J(\mathcal{V}, v-k)$, and in particular $s$ is an automorphism of  $J(\mathcal{V},k)$ if $k=v/2$. In addition, the natural action of 
the symmetric group $\sym(\mathcal{V})$ on $\binom{\mathcal{V}}{k}$ induces a subgroup of automorphisms of $J(\mathcal{V},k)$. If $k\ne v/2$ then $\sym(\mathcal{V})$ is the full automorphism group, while if $k=v/2$ then the automorphism group of $J(\mathcal{V},k)$ is twice as large, namely  $\sym(\mathcal{V})\times\langle s\rangle$,  see \cite[Theorem 9.1.2]{drg}.

A \emph{code} in $J(\mathcal{V},k)$ is a non-empty vertex subset $\Gamma \subseteq \binom{\mathcal{V}}{k}$ and the elements of $\Gamma$ are called \emph{codewords}. We call $\Gamma$ \emph{complete} if $\Gamma = \binom{\mathcal{V}}{k}$. 
The \emph{automorphism group} $\aut(\Gamma)$ of a code $\Gamma$ in  $J(\mathcal{V},k)$ is the full subgroup of graph automorphisms which leave $\Gamma$ invariant. 
\textit{Strongly incidence-transitive codes} in Johnson graphs were introduced by Robert Liebler and the fourth author in connection with their investigations of \textit{neighbour-transitive codes} in Johnson graphs \cite{liebprae}, and are defined as follows. 

\begin{definition}\label{defsit}
A code $\code$ in $J(\mathcal{V},k)$ is called \emph{$X-$strongly incidence-transitive} if $X$ is a subgroup of $\aut(\Gamma)\cap \sym(\mathcal{V})$ which acts transitively on $\code$, and for each codeword $\Delta \in \code$, the setwise stabiliser $X_\Delta$ acts transitively on $\Delta \times \overline{\Delta}$.
\end{definition}
The complementing map $s$ maps each code $\Gamma$ in $J(\mathcal{V},k) $ to its {\em complementary code} $\overline{\Gamma}\subseteq J(\mathcal{V}, v-k)$, where  $\overline{\Gamma}=\{\overline{\Delta} \mid \Delta\in\Gamma \} =\{\Delta^s \mid \Delta\in\Gamma \}.$ 
From Definition~\ref{defsit} it is clear that $\Gamma \subseteq J(\mathcal{V},k)$ is $X-$strongly incidence-transitive if and only if its complementary code $\overline{\Gamma}$ is $X-$strongly incidence-transitive.
The research of Liebler and Praeger \cite{liebprae}, of Neunh{\"o}ffer and Praeger \cite{sporadicnt}, and of Durante~\cite{ND} (see also the exposition in \cite{BCC}), led to the discovery of several new infinite families of strongly incidence-transitive codes, but also left a major unresolved question concerning the existence of strongly incidence-transitive codes constructed from the $2$-transitive actions of the symplectic group $\syp_{2n}(2)$ of degrees $2^{n-1}(2^n+e)$, for $e \in \{1,-1\}$. We provide a brief overview of these actions for the purpose of stating the Main Theorem of this paper. Further discussion is available in Section \ref{jordansteiner}.

Let $X = \syp_{2n}(2)$ and let $B$ denote an $X$-invariant symplectic form defined on $V = \mathbb{F}_2^{2n}$.  We say $(V,B) = (\mathbb{F}_2^{2n},B)$ is a symplectic space. We denote by $\mathcal{Q}$ (or $\mathcal{Q}(V)$ if we need to specify the vector space) the set of all quadratic forms $\phi:V\to \mathbb{F}_2$ which polarise to $B$, that is, which satisfy the equation $B(x,y) = \phi(x+y) - \phi(x) - \phi(y)$, for all $x,y \in V$. For $\epsilon \in \{ +,- \}$ we denote by $\mathcal{Q}^\epsilon$ the set of all elements in $\mathcal{Q}$ of type $\epsilon$. If $\phi \in \mathcal{Q}$ and $U$ is a subspace of $V$ then we denote the restriction of $\phi$ to $U$ by $\phi_U: U \to \mathbb{F}_2$. The symplectic group $X$ admits a pair of distinct $2$-transitive actions on the sets $\mathcal{Q}^\epsilon$ of degrees $2^{n-1}(2^n+e)$ where $e=1$ if $\epsilon=+$ and $e=-1$ if $\epsilon=-$; these are known as the \emph{Jordan-Steiner actions}. From this point on, we use a slight abuse of notation and replace $e$ with $\epsilon$, so that for instance, we can write $|\mathcal{Q}^\epsilon| = 2^{n-1}(2^n+\epsilon)$. In other words, depending on the context, $\epsilon\in \{ +,- \}$ (when describing the type of a quadratic form) or $\epsilon\in \{ 1,-1 \}$ (when counting anything), and this is not ambiguous. For $\phi \in \mathcal{Q}$ we define $\sing(\phi) = \{ x \in V \mid \phi(x) = 0 \}$. Following \cite[Problem 1]{liebprae}, the present article contributes to the classification of strongly incidence-transitive codes associated with the Jordan-Steiner actions. We construct two new infinite families of $X$-strongly incidence-transitive codes, and prove that these are the only examples for which the stabiliser of a codeword is contained in a geometric maximal subgroup of $\syp_{2n}(2)$ (see Theorem \ref{mainthm} below). Note that $\syp_{2n}(2)$ is $2$-transitive on $\mathcal{Q}^\epsilon$, and therefore each code can be viewed as the block-set of a balanced incomplete block design with point-set $\mathcal{Q}^\epsilon$.
%

In the remainder of this paper, $(\mathbb{F}_2^{2n},B)$ is a symplectic space and  $\mathcal{Q}^\epsilon$ is the set of quadratic forms on $\mathbb{F}_2^{2n}$, of type $\varepsilon$, polarising to $B$. Next we give the two constructions mentioned above.  

\begin{constr}\samepage
\label{ndcode} 
Let $\epsilon, \epsilon' \in \{ +,- \}$ and $n,d \in \mathbb{Z}$ such that $n \geq 2$, $1 \leq d \leq n-1$ and $(n,d,\epsilon)\neq (2,1,+)$, and let $k= 2^{n-2}(2^d+\varepsilon')(2^{n-d}+\varepsilon\varepsilon')$. 
The code $\Gamma(n,d,\epsilon,\epsilon')$ in $J(\calq^\epsilon,k)$ consists of one codeword $\Delta(U)$ for each $2d$-dimensional nondegenerate subspace $U$, namely $\Delta(U)$ is the set of all the quadratic forms $\phi \in \calq^\epsilon$ such that $\varphi_U$ is of type $\varepsilon'$ and $\varphi_{U^\perp}$ is of type $\varepsilon\varepsilon'$.
\end{constr}

To clarify, type $\varepsilon\varepsilon'$ means $+$ if $\varepsilon=\varepsilon'$ and $-$ if $\varepsilon\neq \varepsilon'$ (we think of $\varepsilon\varepsilon'$ as a product, using the convention that $+$ corresponds to $1$ and $-$ corresponds to $-1$, as above).
The construction makes sense even if  $(n,d,\epsilon)=(2,1,+)$, but in this case it produces codes in the complete graph $K_{10}=J(10,k)$ with $k\in\{1, 9\}$. Thus we disregard this case.
We note that in Construction~\ref{ndcode}, the code $\Gamma(n,d,\epsilon,-\epsilon')$ is the complementary code of  $\Gamma(n,d,\epsilon,\epsilon')$. Moreover, by replacing $U$ by $U^\perp$ in Construction~\ref{ndcode}, we see that  
\[
\Gamma(n,d,\epsilon,\epsilon')=\Gamma(n,n-d,\epsilon,\epsilon\epsilon').
\] 
In particular we note that, for $n$ even, $\Gamma(n,\frac{n}2,-,-)=\Gamma(n,\frac{n}2,-,+)$, and that this code is self-complementary (that is to say, it is equal to its complementary code).
In Theorem~\ref{p:nd}, we determine the full automorphism groups of the codes $\Gamma(n,d,\epsilon,\epsilon')$, and prove that they are $\syp_{2n}(2)$-strongly incidence-transitive.  We have not been able to determine the minimum distance in general, but we make the following conjecture.

\begin{conjecture}
The minimum distance of the code $\Gamma(n,d,\epsilon,\epsilon')$ is $2^{2n-4}$ unless $\varepsilon=+$ and $d=(n\pm 1)/2$, in which case it is $2^{2n-4}-2^{n-3}$.
\end{conjecture}

\noindent
In Corollary~\ref{cor:distanceND} we prove that these quantities are upper bounds for the minimum distance, and computations using Magma~\cite{Magma} show that they are the correct values in all cases when $n\leq 5$.

\begin{constr}
\label{ticode}
Let $\epsilon \in \{ +,- \}$, $\delta \in \{ 0, 1\}$, and $n,d \in \mathbb{Z}$ such that $n \geq 2$, $1 \leq d \leq n$, and $(d,\varepsilon)\neq  (n,-)$, and let 
\begin{equation}\label{eqticode}
k=
\begin{cases}
2^{n-1}(2^{n-d}+\varepsilon) & \text{ if } \delta=0 \\
2^{2n-d-1}(2^{d}-1) & \text{ if } \delta=1.
\end{cases}    
\end{equation}
The code $\Gamma(n,d,\epsilon,\delta)$ in $J(\calq^\epsilon,k)$ consists of one codeword 
$\Delta^\delta(U)$ for each $d$-dimensional totally-isotropic subspace $U$, namely $\Delta^\delta(U) = \{ \phi \in \mathcal{Q}^\epsilon \mid \dim(\sing(\phi)\cap U) = d-\delta \}$. 
\end{constr}
We note that  in  Construction~\ref{ticode}, the code  $\Gamma(n,d,\epsilon,\delta)$ is the complementary code of  $\Gamma(n,d,\epsilon,1-\delta)$. Also, we note that the condition  $(d,\varepsilon)\neq  (n,-)$ is necessary since there are no $n$-dimensional totally isotropic subspaces if $\varepsilon=-$.
In Theorem~\ref{p:ti}, we determine the full automorphism groups of the codes $\Gamma(n,d,\epsilon,\delta)$, and prove that they are $\syp_{2n}(2)$-strongly incidence-transitive. Then in Theorem~\ref{p:distanceti} we determine the minimum distance $D$; the value of $D$ depends on the parameters, but in all cases $D\geq 2^{n-2}$.

By a \emph{geometric subgroup} of $\syp_{2n}(2)$ we mean a subgroup of a maximal geometric subgroup contained in one of the `geometric Aschbacher classes' $\mathcal{C}_1$-$\mathcal{C}_8$ which are defined in \cite{aschbacher84}, and are summarised in Table \ref{syp2n2max}. Detailed information about the maximal subgroups of $\syp_{2n}(2)$ is readily available in \cite{colva,kl}.

\begin{theorem}[Main Theorem]
	\label{mainthm}
	Let $\Gamma$  be a code in $J(\calq^\epsilon,k)$, for some $k$ satisfying 
	$1< k < |\mathcal{Q}^\varepsilon|-1$. 
Suppose that $X = \syp_{2n}(2) \leq \aut(\Gamma) \cap \sym(\mathcal{Q}^\varepsilon)$, and that $\Gamma$ is $X$-strongly incidence-transitive. Suppose also that the stabiliser of  a codeword  $\Delta$ is  a geometric subgroup of $\syp_{2n}(2)$. Then  $\Gamma$ arises from Construction \ref{ndcode} or \ref{ticode}. 
\end{theorem}


%



\begin{remark} We make a few comments.
\begin{enumerate}[(a)]
\item In Theorem~\ref{mainthm}, 	$2 \leq k \leq |\mathcal{Q}^\varepsilon|-2$, so $|\mathcal{Q}^\varepsilon|\geq 4$, and it follows that  $n\geq 2$ and $|\mathcal{Q}^\varepsilon|\geq 6$.

\item If $\Gamma$ is a nonempty subset of $\binom{\mathcal{Q}^\varepsilon}{k}$ with $k \in \left\{2,|\mathcal{Q}^\epsilon|-2\right\}$, then it follows from the $2$-transitivity of $X$ on $\mathcal{Q}^\epsilon$ that $\Gamma = \binom{\mathcal{Q}^\varepsilon}{k}$. Thus $\Aut(\Gamma)$ contains $\sym(\mathcal{Q}^\varepsilon)$, and in fact equality holds since 
$k\ne |\mathcal{Q}^\varepsilon|/2$ because of the restrictions noted in part (a). If $\Gamma$ were $X$-strongly incidence-transitive then,  by Definition~\ref{defsit}, the stabiliser in $X$ of a codeword $\Delta$ would be transitive on $\Delta\times (\mathcal{Q}^\varepsilon\setminus\Delta)$. Hence if $|\Delta|=2$ then the pointwise stabiliser of $\Delta$ would be transitive on $\mathcal{Q}^\varepsilon\setminus\Delta$, and similarly, if $|\mathcal{Q}^\varepsilon\setminus\Delta|=2$ then the pointwise stabiliser of $\mathcal{Q}^\varepsilon\setminus\Delta$ would be transitive on $\Delta$. In either case it would follow that $X$ is $3$-transitive on $\mathcal{Q}^\varepsilon$. This is the case if and only if $(n,\varepsilon)=(2,-)$ and the code $\Gamma$ is therefore 
$\Gamma(2,1,-,0)=\binom{\calq^-}{2} $ or its complement $\Gamma(2,1,-,1)=\binom{\calq^-}{4}$  as in Construction \ref{ticode}. In these cases, $k$ is $2$ or $|\mathcal{Q}^-|-2=4$ respectively, and $\aut(\Gamma)=\syp_4(2)=\sym(\mathcal{Q}^\varepsilon)\cong\sym(6)$.
These examples are complete codes, which are not very interesting.
\item  There is another complete code: $\Gamma(2,1,-,+)=\Gamma(2,1,-,-)=\binom{\calq^-}{3}$  as in Construction \ref{ndcode}. In this case, $k=3$, and $\aut(\Gamma)=\syp_4(2)=\sym(\mathcal{Q}^\varepsilon)\cong\sym(6)$. All codes with $(n,d,\epsilon)\neq (2,1,-)$ are non-complete, see Lemma \ref{le:noncomplete}.
\item There is an exceptional isomorphism between codes from Construction~\ref{ticode}, namely $\Gamma(2,1,+,1)$ and $\Gamma(2,2,+,0)$ are isomorphic, under an outer automorphism of $\syp_4(2)\cong\sym(6)\cong\psigmal_2(9)$.
These codes have $15$ codewords of size $4$. This isomorphism arises because such an outer automorphism is a so-called `graph automorphism' of $\syp_4(2)$ which interchanges totally isotropic $1$-subspaces and totally isotropic $2$-subspaces. 

\item By Theorem \ref{p:distanceti}, the minimum distances of codes arising from Construction \ref{ticode} are unbounded, and in particular are at least $2^{n-2}$ for any $d, \epsilon, \delta$. Hence any desired error-correction capability can be achieved, However as the minimum distance grows larger so does $n$, and hence the cost of constructing or storing the codewords also grows. 
For instance, $\Gamma(n,1,+,1)$ is a code in $J(\mathcal{Q}^+,2^{2n-2})$ with $2^{2n}-1$ codewords, minimum distance of $2^{2n-3}$ and $|\mathcal{Q}^+|=2^{n-1}(2^n+1)$.
\end{enumerate}
\end{remark}
%


%




%
Background on bilinear and quadratic forms is provided in Section \ref{background}. In Section \ref{jordansteiner}, we describe the Jordan-Steiner actions and outline some key properties. Section~\ref{sec:sitcodes} contains two useful general results about strongly incidence-transitive codes.
In Section \ref{sec:reducible} we study the subgroups of $\syp_{2n}(2)$ which act reducibly on $\mathbb{F}_2^{2n}$. In Section \ref{sec:irred} we study the geometric subgroups of $\syp_{2n}(2)$ which act irreducibly on $\mathbb{F}_2^{2n}$. In Section \ref{proofmain} we pull everything together to prove Theorem~ \ref{mainthm}. Detailed information about the families of codes is provided in Theorem \ref{p:nd} and Theorem \ref{p:ti}. 
\begin{table}
\small
\label{syp2n2max}
\begin{tabular}{|c|c|c|p{0.43\textwidth}|}
\toprule
Case & Class & Structure  & Description and conditions \\
\midrule
(a) & $\mathcal{C}_1$ & $\syp_{2d}(2) \times \syp_{2(n-d)}(2)$  & Stabiliser of a $2d$-dimensional nondegenerate subspace with $1 \leq d \leq (n-1)/2$ \\
(b) & $\mathcal{C}_1$ & $2^{d(d+1)/2}.2^{2d(n-d)} \rtimes  \gl_d(2) \times \syp_{2(n-d)}(2)$  & Stabiliser of a $d$-dimensional totally-isotropic subspace with $1 \leq d \leq n$ \\
(c) & $\mathcal{C}_2$ & $\syp_{2m}(2) \wr S_t$  & Stabiliser of a decomposition $\oplus_{i=1}^t V_i$ into $2m$-dimensional nondegenerate subspaces where $m = n/t<n$ \\
(d) & $\mathcal{C}_3$& $\syp_{2m}(2^b) \rtimes C_b$  & Stabiliser of a $b$-spread with $b$ prime and $n = mb$ \\
(e) & $\mathcal{C}_8$ & $\oo^\varepsilon_{2n}(2), \hspace{0.1cm} \epsilon = \{+,-\}$  & Isometry group of an $\epsilon$-type quadratic form which polarises to $B$ \\
\bottomrule
\end{tabular}
\caption{The geometric maximal subgroups of the symplectic group $\syp_{2n}(2)$. See \cite[Theorem 3.7]{wilson}. In particular classes $\mathcal{C}_4$-$\mathcal{C}_7$ are empty for $\syp_{2n}(2)$, see \cite[Table 3.5.C]{kl}. }
\end{table}
\section{Background on bilinear and quadratic forms}
\label{background}

Let $V$ be a vector space over the finite field $\mathbb{F}_q$. A \emph{bilinear form} on $V$ is a map $B:V\times V \to \mathbb{F}_q$ such that for all $u,v,w \in V$ and $\alpha \in \mathbb{F}_q$,
\begin{align*}
B(\alpha u+v,w) &= \alpha B(u,w) + B(v,w), \\
B(u, \alpha v + w) &= \alpha B(u,v) + B(u,w).
\end{align*}
We say $B$ is \emph{reflexive} if $B(u,v) = 0$ implies $B(v,u) = 0$ for all $u,v \in V$. If $U$ is a subset of $V$ then we define
\[U^\perp = \{ v \in V \mid B(u,v) = 0 \text{ for all } u \in U \},\] which is a subspace of $V$.
We call $U^\perp$ the \emph{orthogonal complement} of $U$, and in particular, $V^\perp$ is called the \emph{radical} of $B$. Note that if $U=\{c\}$ is a singleton, we simply write $c^\perp$ instead of $\{c\}^\perp$.
If $U$ and $W$ are subsets of $V$ then we say they are \emph{orthogonal} if $U \subseteq W^\perp$ (or equivalently, $W \subseteq U^\perp$ if $B$ is reflexive).    We call $B$ \emph{nondegenerate} if $V^\perp = \{ 0 \}$. If $U$ is a subspace of $V$ and $B$ is nondegenerate, then $\dim(V) = \dim(U)+\dim(U^\perp)$.
In this case, we say $U$ is \emph{nondegenerate} if $U \cap U^\perp = \{0\}$ and \emph{totally-isotropic} if $U \leq U^\perp$. 
If $U$ and $W$ are subspaces of $V$ then we write $U+W = \{ u+w \mid (u,w) \in U\times W \}$. If in addition $U \cap W = \{ 0 \}$ then we write $U + W = U \oplus W$, and if $B(u,w) = 0$ for all $(u,w) \in U\times W$ and $U \cap W = \{ 0 \}$ then we write $U \oplus W = U \perp W$.

We call $B$ \emph{symmetric} if $B(u,v) = B(v,u)$ for all $u,v\in V$ and we call $B$ \emph{skew-symmetric} if $B(u,v) = -B(v,u)$ for all $u,v\in V$. If $q$ is even then $B$ is symmetric if and only if $B$ is skew-symmetric. We call $B$ \emph{alternating} if $B(v,v) = 0$ for all $v\in V$. Every alternating form is skew-symmetric, but the converse holds if and only if $q$ is odd. We call $B$ \emph{symplectic} if it is both nondegenerate and alternating; in this case we call the pair $(V,B)$ a \emph{symplectic space}.

A map $\phi: V \to \mathbb{F}_q$ is a \emph{quadratic form} if
\begin{enumerate}[(i)]
\item $\phi(\lambda v) = \lambda^2 \phi(v)$ for all $v \in V$ and $\lambda \in \mathbb{F}_q$, and
\item the map $B_\mathcal{\phi}:V\times V \to \mathbb{F}_q$ defined by
\begin{align}
\label{polar}
B_\phi(u,v) = \phi(u+v) - \phi(u) - \phi(v)
\end{align}
is a bilinear form.
\end{enumerate}
Equation \eqref{polar} is called the \emph{polarisation equation} and $B_\phi$ is called the \emph{polar form} of $\phi$. Note that if $\phi$ is a quadratic form then  $B_\phi$ is symmetric and $B_\phi(u,u)=2\phi(u)$. We say $\phi$ is nondegenerate if $B_\phi$ is nondegenerate. Since the operation of subtraction is equivalent to addition in characteristic two, we often rewrite Equation \eqref{polar} as
\begin{align*}
B_\phi(u,v) = \phi(u+v) + \phi(u) + \phi(v)
\end{align*}
when $q$ is even. 

The \emph{general linear group} $\gl(V)$ is the group of all invertible linear transformations of $V$ where, for $g\in\gl(V)$ and $v\in V$, $vg$ denotes the image of $v$ under $g$. 
Let $B$ be a nondegenerate reflexive bilinear form and let $\phi$ denote a nondegenerate quadratic form. An element $g \in \gl(V)$ is called an \emph{isometry} of $B$ if $B(ug,vg) = B(u,v)$ for all $u,v \in V$. Similarly, $g$ is called an isometry of $\phi$ if $\phi(ug)=\phi(u)$ for all $u \in U$. The set of isometries of a bilinear or quadratic form on $V$ is a subgroup of $\gl(V)$, which we denote by $\isom(B)$ or $\isom(\phi)$. Note $\isom(\phi)\leq \isom(B_\phi)$. If $B$ is a symplectic form, then $\isom(B)$ is denoted by $\syp(V)$ or $\syp_{2n}(q)$ if $V$ is $2n$-dimensional over the field of order $q$.

 If $U$ is a subspace of $V$ then we use $\phi_U$ and $B_{U}$ to denote the restrictions of $\phi$ and $B$ to $U$ and $U\times U$, respectively. It is possible therefore that $B$ is degenerate but $B_U$ is nondegenerate. If $\phi_U=0$ (meaning $\phi(u)=0$ for all $u\in U$), then $U$ is said to be $\phi$-\emph{singular}.

If $B$ is a symplectic form on $V$ then we call a basis $\{e_1,f_1,e_2,f_2,\dots,e_n,f_n\}$ for $V$ a \emph{symplectic basis} if for all $1 \leq i,j \leq n$ we have $B(e_i,e_j) = B(f_i,f_j) = 0$, $B(e_i,f_i) = 1$ and $B(e_i,f_j) = 0$ if $i \neq j$. Every symplectic space has a symplectic basis \cite[Proposition 2.5.3]{kl}. 

\section{Jordan-Steiner actions}
\label{jordansteiner}
Let $\mathbb{F}$ be a finite field of characteristic two, and let $V$ be a vector space over $\mathbb{F}$ (of even dimension) equipped with a symplectic form $B$. We denote the full isometry group of $B$ by $X$. Let $\mathcal{Q}(V)$ denote the set of all quadratic forms on $V$ which polarise to $B$. We write $\mathcal{Q} = \mathcal{Q}(V)$ when there is no danger of ambiguity.

 Given $\varphi \in \mathcal{Q}$ and $g \in X$, we define a function $\varphi^g: V \rightarrow \mathbb{F}$ by
\begin{equation}
	\label{action}
	\varphi^g(x) = \varphi (xg^{-1}).
\end{equation}
It is routine to verify that Equation \eqref{action} defines a group action of $X$ on $\mathcal{Q}$. 
Note that the stabiliser $X_\phi$ of a quadratic form $\phi$ consists of the maps $g\in X$ such that $\varphi (xg^{-1})=\varphi (x)$ for all $x\in V$, or equivalently such that $\varphi (x)=\varphi (xg)$ for all $x\in V$.  We partition $\mathcal{Q}$ into two subsets which we denote by $\mathcal{Q}^\epsilon$, with $\epsilon\in \lbrace +,- \rbrace$. The elements of $\mathcal{Q}^+$ are the \emph{hyperbolic quadratic forms} contained in $\mathcal{Q}$, or equivalently, the quadratic forms in $\mathcal{Q}$ whose isometry group is an orthogonal group of \emph{plus-type}. Similarly, the elements of $\mathcal{Q}^-$ are the \emph{elliptic quadratic forms} contained in $\mathcal{Q}$, or equivalently, the quadratic forms in $\mathcal{Q}$ whose isometry group is an orthogonal group of \emph{minus-type}.

The sets $\mathcal{Q}^\epsilon$ are $X$-invariant and $X$ acts transitively on each. If $\mathbb{F}=\mathbb{F}_2$, then $X$ acts $2$-transitively on $\mathcal{Q}^\epsilon$ for each $\epsilon \in \{+,-\}$ (see \cite[Section 7.7]{dixonmortimer} for details). The actions of $X$ in the case $q=2$ are called the \emph{Jordan-Steiner actions}. See \cite{delsarte-alternating,kantorquadratic,sastrysincode} for some applications of the Jordan-Steiner actions to coding and design theory. The submodule structure of the associated permutation modules is studied in \cite{sastrysin}.

\subsection{Forms polarising to a given symplectic form}
\label{affinestructure}

Let $\mathbb{F}$ be a finite field of order $q = 2^b$. Let $\alpha \in \mathbb{F}$.

Let $V$ be a finite even dimensional vector space over $\mathbb{F}$ equipped with a symplectic form $B$. Our goal in Section \ref{affinestructure} is to show that if $\phi$ and $\psi$ are quadratic forms on $V$ which polarise to $B$, then there exists a unique $c \in V$ such that the equation
\begin{equation}
\label{affinerelationbeforeproof}
\phi(v) = \psi(v) + B(v,c)^2
\end{equation}
holds for all $v\in V$. Equation \eqref{affinerelationbeforeproof} is discussed and used in \cite{sastrysin}.

\begin{lemma}\samepage
\label{linearfunctional}
If $\phi$ and $\psi$ are quadratic forms on $V$ which polarise to $B$, then the function $f:V \to \mathbb{F}$ defined by 
\begin{equation}
\label{eqn_linearfunctional}
f(v) = \phi(v)^{q/2}+\psi(v)^{q/2}
\end{equation}
is linear.
\end{lemma}
\begin{proof}
If $\alpha \in \mathbb{F}$ then let $h(\alpha) = \alpha^{q/2}$. The automorphism group of $\mathbb{F}$ is cyclic of order $b$, and it is generated by the Frobenius automorphism $\alpha \mapsto \alpha^2$, so $h$ is an automorphism of $\mathbb{F}$. Since $q$ is a power of $2$ and $\alpha^q = \alpha$ (see \cite[Theorem 1.2]{simeon}), we have that $h(\alpha^2)=\alpha$. By squaring both sides of Equation \eqref{eqn_linearfunctional}, we find $f(v)^2 = \phi(v) + \psi(v)$ for all $v\in V$.
If $u,v \in V$ then
\begin{align*}
f(u+v)^2 &= \phi(u+v)+\psi(u+v) 	& \text{(by definition)} \\
	&= \phi(u) + \phi(v) + B(u,v) + \psi(u) + \psi(v) + B(u,v) & \text{(polarisation)} \\
	&= \phi(u) + \psi(u) + \phi(v) + \psi(v) & \text{(characteristic two)} \\
		&= f(u)^2 + f(v)^2 & \text{(by definition)}\\
	&= (f(u) + f(v))^2. & \text{(characteristic two)}
\end{align*}
Applying $h$ to both sides of this equation, we find $f(u+v) = f(u) + f(v)$.
Similarly, if $\alpha \in \mathbb{F}$ and $v \in V$ then
\begin{align*}
f(\alpha v)^2 &= \phi(\alpha v)+\psi(\alpha v) & \text{(by definition)} \\
	&= \alpha^2 \phi(v)+\alpha^2 \psi(v) & \text{(quadratic forms)} \\
	&= \left(\alpha f(v)\right)^2.
\end{align*}
Applying $h$ to both sides of this equation, we find $f(\alpha v) = \alpha f(v)$. 
Therefore $f$ is linear.
\end{proof}
If $W$ is an $m$-dimensional vector space over a field $\mathbb{K}$ then the \emph{dual space} of $W$ is the vector space of linear functions from $W$ to $\mathbb{K}$, where addition and scalar multiplication of functions are defined pointwise. The dual space of $W$ is denoted $W^\star$. Refer to \cite[Chapter 2]{roman} for a detailed discussion.
\begin{lemma}
\label{dualspaceiso}
Let $W$ be an $m$-dimensional vector space over a field $\mathbb{K}$ and let $B$ be a nondegenerate bilinear form on $W$. For each $w \in W$, define a function $\lambda(w): W \to \mathbb{K}$ by $\lambda(w)(u) = B(u,w)$ for all $u \in W$. Then $\lambda: W \to W^\star$, given by $\lambda:w\to \lambda(w)$, is a vector space isomorphism.
\end{lemma}
\begin{proof}
Since $B$ is linear in the first variable, we have $\lambda(w) \in W^\star$ for all $w \in W$, and since $B$ is linear in the second variable, the map $\lambda$ is a linear transformation. 
Additionally, if $\lambda(w_1) = \lambda(w_2)$ then $B(u,w_1) = B(u,w_2)$ for all $u \in W$, which implies $B(u,w_1-w_2) = 0$ for all $u \in U$. Since $B$ is nondegenerate, this implies $w_1 = w_2$. Therefore $\lambda$ is injective. Finally, \cite[Theorem 3.11]{roman} shows that $\dim(W) = \dim(W^\star)$ and therefore $\lambda$ is an isomorphism.
\end{proof}

\begin{lemma}
\label{param}
Let $V$ be a vector space over a finite field $\mathbb{F}$ of characteristic two, and let $B$ be a symplectic form on $V$. If $\phi$ and $\psi$ are quadratic forms on $V$ which polarise to $B$, then there exists a unique $c \in V$ such that 
\begin{equation}
\label{parameq}
\phi(v) + \psi(v) = B(v,c)^2 \text{ for all } v \in V.
\end{equation}
Conversely, if $\phi$ is a quadratic form on $V$ which polarises to $B$ then for any $c \in V$, the function defined by $\psi(v) = \phi(v) + B(v,c)^2$ is a quadratic form on $V$ which polarises to $B$.
\end{lemma}

\begin{proof}
First, suppose that $\phi$ and $\psi$ are quadratic forms on $V$ which polarise to $B$. By Lemma \ref{linearfunctional}, the map $f(v) = \phi(v)^{q/2}+\psi(v)^{q/2}$ is linear, that is to say, $f\in V^*$, and therefore Lemma \ref{dualspaceiso} implies that there exists a unique $c \in V$ such that $f(v) = B(v,c)$ for all $v \in V$. Squaring both sides, we find that $\phi(v) + \psi(v) = B(v,c)^2$ for all $v \in V$.

Conversely, let $c \in V$, let $\phi$ be a quadratic form on $V$ which polarises to $B$ and define $\psi(v) = \phi(v) + B(v,c)^2$ for all $v \in V$. Then for all $\lambda \in \F$ and all $u, v\in V$, we have
\begin{align*}
\psi(\lambda v) &= \phi(\lambda v) + B(\lambda v, c)^2 \\
	&= \lambda^2 \phi(v) + \lambda^2 B(v,c)^2 \\
	&= \lambda^2 (\phi(v) + B(v,c)^2) \\
	&= \lambda^2 \psi(v),
\end{align*} 
and
\begin{align*} 
\psi(u+v) &= \phi(u+v) + B(u+v,c)^2 \\
	&= \phi(u)+\phi(v)+B(u,v)+B(u,c)^2+B(v,c)^2\\
	&= \psi(u) + \psi(v) + B(u,v).
\end{align*}
Therefore $\psi$ is a quadratic form on $V$ which polarises to $B$.
\end{proof}

\begin{notation}
\label{notation}
If $\phi \in \mathcal{Q}$ and $c \in V$ then we use $\phi_c$ to denote the unique quadratic form in $\mathcal{Q}$ such that 
\begin{equation}
\label{parameq2}
\phi_c(v)  = \phi(v) + B(v,c)^2 \text{ for all } v \in V.
\end{equation}
Similarly, if $\phi,\psi\in \mathcal{Q}$ then we use $c_{\phi\psi}$ to denote the unique vector $c\in V$ such that Equation \eqref{parameq} holds. 
\end{notation}

\begin{remark}\label{rem3.5}\leavevmode\vspace{-\baselineskip} 
\begin{enumerate}[(a)]
\item If $\F=\F_2$ then Equation \eqref{parameq} simplifies to $\psi(v) = \phi(v) + B(v,c)$ for all $v \in V$.
\item If $\phi, \psi$ and $\rho$ are quadratic forms in $\mathcal{Q}$ then $c_{\phi\rho} = c_{\phi\psi} + c_{\psi\rho}$.
\item If $\phi$ and $\psi$ are related by Equation \eqref{parameq} then $\phi(v) = \psi(v)$ if and only if $v \in c^\perp$.
\item Note $\phi_0=\phi$.
\item Lemma \ref{param} implies that $\mathcal{Q}$ has the structure of an \emph{affine geometry} in the sense of \cite[Chapter 2]{bergergeometry1}. The points of the geometry are the elements of $\mathcal{Q}$ and if $\phi,\psi \in \mathcal{Q}$ then the translation vector $\overrightarrow{\phi\psi}$ is $c_{\phi\psi}$. Alternatively, the discussion in \cite[Section 2]{sastrysin} explains that $\mathcal{Q}$ corresponds to a coset in the $\mathbb{F}$-vector space $\bigwedge^2(V^\star)$ of alternating bilinear forms on $V$, and therefore $\mathcal{Q}$ inherits its affine structure from $\bigwedge^2(V^\star)$.
\end{enumerate}
\end{remark}

\subsection{Singular vectors}
\label{subsec-coordinates}

Throughout this section (even if not mentioned explicitly), $V$ is a finite $2n$-dimensional vector space over a field $\F=\F_q$ of characteristic two, $B$ is a symplectic form on $V$, and  $\mathcal{Q}$ is the set of quadratic forms which polarise to $B$.
Recall that if $\phi \in \mathcal{Q}$ then we define $\sing(\varphi) = \{ x \in V \mid \phi(x) = 0 \}$. The elements of $\sing(\varphi)$ are called \emph{$\varphi-$singular} vectors. Note that by definition, $\sing(\varphi)$ contains the zero vector for all $\varphi \in \mathcal{Q}$. If $\phi \in \mathcal{Q}^\epsilon$ then from \cite[Theorem 1.14]{ht} we have
\begin{equation}
\label{numsingular}
|\sing(\phi)| = (q^{n-1}+\epsilon)(q^n-\epsilon)+1 = q^{n-1}(q^n+\epsilon(q-1)).
\end{equation}
%

\begin{lemma}
\label{lem_singrelations}
Let $q=2$, and let $\phi$ and $\psi$ be quadratic forms on $V$ which polarise to $B$ and suppose that Equation \eqref{parameq} holds for some $c \in V$. Then
\begin{equation}
		\label{singrelations}
		\sing(\psi) = \left\{
		\begin{array}{lr}
			\sing(\varphi)+c               & \mbox{\ if\ }c \in \sing(\varphi)    \\
			(V \setminus \sing(\varphi))+c & 
			\mbox{\ if\ }c \notin \sing(\varphi).
		\end{array}
		\right.
	\end{equation}
	
\end{lemma}

\begin{proof}



Since $x^2=x$ and $2x=0$ for all $x\in \F$, 
\[\psi(v+c) =  \phi(v+c) + B(v+c,c)^2                           
		=          \phi(v) + \phi(c) + B(v,c) + B(v,c)^2 + B(c,c)^2 
		=   \phi(v) + \phi(c). 
\]
Suppose first that $c \in \sing(\phi)$. Then $\psi(v+c)=\phi(v)$, and so $\sing(\psi)=\sing(\phi)+c$. 
If $c \notin \sing(\phi)$. Then $\psi(v+c)=\phi(v)+1$, and so the vectors $w\in\sing(\psi)$ are precisely the vectors of the form $w=v+c$ where $v\in V\setminus\sing(\phi)$. Hence $\sing(\psi)=(V\setminus\sing(\phi))+c$. 
\end{proof}


	%
	%


\begin{corollary}
\label{formvectorbijection}
Let 
$\phi \in \mathcal{Q}$. Define a map $\lambda_\phi: V \to \mathcal{Q}$ by $\lambda_\phi(c) = \phi_c$. Then $\lambda_\phi$ is a bijection, and the type of $\phi_c$ can be identified by the size of $\sing(\phi_c)$. Moreover, if $q=2$ and $\phi \in \mathcal{Q}^\epsilon$ then  $\lambda_\phi(c)\in \mathcal{Q}^\epsilon$ if and only if  $c\in\sing(\phi)$.
\end{corollary}

\begin{proof}
It follows immediately from Lemma \ref{param} and Notation~\ref{notation} that $\mathcal{Q} =\{ \phi_c \mid c\in V\}$, and hence $\lambda_\phi$ is a bijection.
By \eqref{numsingular}, the number of elements of $\sing(\phi)$ depends only on $n, q$ and $\epsilon$, and for fixed $n, q$, these numbers are different  for $\epsilon=+$ and $\epsilon=-$. Thus if $q=2$ and $\phi \in \mathcal{Q}^\epsilon$ then, by Lemma~\ref{lem_singrelations} (especially \eqref{singrelations}), $|\sing(\phi_c)|=|\sing(\phi)|$ if and only if $c\in\sing(\phi)$. Hence $\lambda_\phi(c)\in\mathcal{Q}^\epsilon$ if and only if $c\in\sing(\phi)$.
\end{proof}

\begin{lemma}
	\label{formvectorpermiso}
	For each $\phi \in \mathcal{Q}$, the action of $X_\phi$ on $\mathcal{Q}$ defined by Equation \eqref{action} is equivalent to the natural action on $V$.
\end{lemma}

\begin{proof}
	Let $\phi \in \mathcal{Q}$. By Corollary \ref{formvectorbijection}, $\lambda_{\varphi}$ is a bijection. By Lemma \ref{param}, for each $\psi \in \mathcal{Q}$ there exists a unique $c \in V$ such that $\psi = \phi_c$. Then for all $\psi = \phi_c \in \mathcal{Q}$, $g \in X_{\phi}$, and $v \in V$ we have
	\begin{align*}
		\varphi_c^g(v) = & \varphi_c(v g^{-1})
		= \varphi(vg^{-1})+B(vg^{-1},c)^2 
		= \varphi^g(v)+B(v,cg)^2 \\
		=            & \varphi(v)+B(v,cg)^2           
		=             \varphi_{cg}(v).
	\end{align*}
	Thus $(\varphi_c)^g = \varphi_{cg}$, that is, $\lambda_{\varphi}(c)^{g} = (\varphi _c) ^g = \varphi_{cg} = \lambda_{\varphi}(cg)$, so the action of $X_\phi$ on $\mathcal{Q}$ is equivalent to the natural action of $X_\phi$ on $V$.
\end{proof}

\begin{corollary}
\label{formsetpermiso}
The action of $X$ on $\mathcal{Q}$ is equivalent to the elementwise action of $X$ on the collection $\{ \sing(\varphi) \subset V \mid \varphi \in \mathcal{Q} \}$.
\end{corollary}
\begin{proof}
For all $g \in X$ we have
\begin{align*}
\sing(\varphi^g) &= \{ v \in V \mid \varphi^g(v)=0 \} = \{ v \in V \mid \varphi(vg^{-1})=0 \} \\
	&= \{ wg \in V \mid \varphi(w)=0 \} = \sing(\varphi)g.\qedhere
\end{align*}
\end{proof}

\begin{corollary}
	\label{affineactioncor}
	Let $\varphi, \varphi_d \in \mathcal{Q}$ with $d \in V$. Let $c \in V$. For all $g \in X$, if $\varphi^g = \varphi_d$ then $(\varphi_c)^g = \varphi_{cg+d}$.
\end{corollary}
\begin{proof}
	For all $\varphi_c \in \mathcal{Q}$ and $g \in X$ satisfying $\varphi^g = \varphi_d$, using \eqref{parameq2}, the following holds for all $v \in V$.
	\begin{align*}
		(\varphi_c)^g(v) = & \varphi_c(v{g^{-1}})\\ =&\varphi(v{g^{-1}}) + B(v{g^{-1}},c)^2 \\
		=             & \varphi_d(v) + B(v,cg)^2                 \\
		=             & \varphi(v) + B(v,d)^2+ B(v,cg)^2                 \\
		=             & \varphi(v) + B(v,cg + d)^2.
	\end{align*}
	so $	(\varphi_c)^g = \varphi_{cg + d}$ as claimed.
\end{proof}
Note that if $g \in X_{\varphi}$ then Corollary \ref{affineactioncor} reduces to $(\varphi_c)^g = \varphi_{cg}$, as in the proof of  Lemma \ref{formvectorpermiso}.

\begin{lemma}\label{uniqueon1space}
    Let $\phi\in \mathcal{Q}$ and let $v$ be a nonsingular vector for $\phi$. Then for every $\lambda\in \mathbb{F}_q$, there exists a unique vector $w$ in the $1$-space spanned by $v$ with $\phi(w)=\lambda$.
\end{lemma}
\begin{proof}
Let $\lambda\in \mathbb{F}_q$ and assume there are two vectors $w_1,w_2$ in the $1$-space spanned by $v$ with $\phi(w_1)=\phi(w_2)=\lambda$. 
Since they are on the same $1$-space, we must have $w_2=\mu w_1$. 
Therefore $\lambda=\phi(w_2)=\mu^2\phi(w_1)=\mu^2\lambda$.
If $\lambda\neq 0$, then $\mu^2=1$ so $\mu=1$ since $q$ is even, and so $w_1=w_2$.
If $\lambda= 0$ and $w_i\neq 0$, then all multiples of $w_i$ are singular, contradicting $v$ being nonsingular. Thus $w_1=w_2=0$.

Hence there is at most one vector $w$ in the $1$-space spanned by $v$ with $\phi(w)=\lambda$.
Since the $1$-space contains $q$ vectors, the statement follows.
\end{proof}

%

We will sometimes be able to work with a specific form of $\epsilon $ type, as in the next lemma. In these cases, assuming $q=2$, we will choose the following forms, with respect to  the standard  symplectic basis $\{ e_i,f_i \mid 1 \leq i \leq n\}$ for $V$ : for all $x = \sum_{i=1}^n(x_i e_i + y_i f_i) \in V$ we set
	\begin{equation}\label{eq-forms}
		\left\{
		\begin{array}{ll}
		\phi_0^+(x) =	\sum_{i=1}^n x_i y_i   \\
		\phi_0^-(x) =	x_n + y_n + \sum_{i=1}^n x_i y_i 
		\end{array}
	\right.
	\end{equation}
Note that, with Notation \ref{notation},
\begin{equation}\label{linktwoforms} \text{ if } \phi=\phi_0^+, \text{ then } \phi_0^-=\phi_{e_n+f_n}.
\end{equation} 
Indeed $\phi_0^-(x)+\phi_0^+(x)=x_n+y_n=B(x,e_n+f_n)=B(x,e_n+f_n)^2$. Since $\phi_0^+(e_n+f_n)=1$, the two forms do indeed have different types, by Corollary \ref{formvectorbijection}, and it is well known that $\phi_0^+\in \mathcal{Q}^+$.

\begin{lemma}\label{le:noncomplete}
    Let $n\geq 2$, and $X=\syp_{2n}(2)$. Let $\Gamma$ be an $X$-strongly incidence-transitive code in $J(\mathcal{Q}^\varepsilon,k)$ for some $k$ satisfying $2 < k < |\mathcal{Q}^\varepsilon|-2$. Then $\Gamma\neq \binom{\mathcal{Q}^\epsilon}{k}$ unless $(n,\epsilon)=(2,-)$.
\end{lemma}
\begin{proof}
Assume that $\Gamma$ is a complete code, that is, $\Gamma= \binom{\mathcal{Q}^\epsilon}{k}$. Let $K=\min\{k,|\calq^\epsilon|-k\}$ so that $|\calq^\epsilon|\geq 2K$.
Since $X$ is transitive on codewords (by Definition~\ref{defsit}), it follows that $X$ is $K$-homogeneous in its action on $\calq^\epsilon$ (that is, $X$ is transitive on $K$-sets). Then by \cite[Theorem 9.4.B]{dixonmortimer}, $X$ is $(K-1)$-transitive. However (see \cite[Table 7.4]{cameronpermutation})  $X$ is $2$-transitive but not $3$-transitive, unless $(n,\epsilon)=(2,-)$ in which case  $X\cong S_6$, 
$|\calq^\epsilon|=6$, $k=3$, and $X$ is indeed transitive on $\Gamma=\binom{\mathcal{Q}^\epsilon}{k}$.  Thus from now on we assume that $(n,\epsilon)\ne (2,-)$, so $K\leq 3$ and $X$ is not $3$-transitive.
If $K=3$, then by \cite[Theorem 9.4.B]{dixonmortimer} or \cite{kantorhomogeneous}, since $X$ is not $3$-transitive, we must have
\[
X\in\{\agl_1(8), \agaml_1(8),\agaml_1(32)\}\]
or  $\psl_2(q) \leq X \leq \pgaml_2(q)$, where $ q = 3\pmod 4$ (all in their  natural action). 
However $X=\syp_{2n}(2)$ is not isomorphic to any of these groups (noting that $X=\syp_4(2)\cong \psl_2(9).2$ does not arise since $9\neq 3\pmod 4$). This contradiction complete the proof.
%
\end{proof}

\section{Strongly incidence-transitive codes}\label{sec:sitcodes}

In this section, we state and prove two simple lemmas that we will need on strongly incidence-transitive codes in Johnson graphs in $J(\mathcal{V},k)$. The first is
about $G$-strongly incidence-transitive codes where a codeword stabiliser is contained in a subgroup $M$ acting transitively and imprimitively on the underlying set $\mathcal{V}$.

\begin{lemma}
	\label{unionsofblocks}
	Let $\code$ be a $G$-strongly incidence-transitive code in $J(\mathcal{V},k)$ with $\Delta \in \code$. Let $M$ be a subgroup of $G\leq \aut(\code)$ which acts transitively on $\mathcal{V}$, while leaving invariant a nontrivial partition $\mathcal{I}$ of $\mathcal{V}$. If $G_\Delta < M \leq G$ then $\Delta$ is a union of parts of $\mathcal{I}$.
\end{lemma}

\begin{proof}
	Let $\Delta \in \code$. Since $\code$ is $G$-strongly incidence-transitive, $G_\Delta$ has two orbits in $\mathcal{V}$, namely $\Delta$ and $\overline{\Delta}$. We assume without loss of generality that $|\Delta| \leq |\overline{\Delta}|$. Suppose that $\Delta$ is not a union of parts of $\mathcal{I}$, so also $\overline{\Delta}$ is not a union of parts either. Therefore there exists $\blocky \in \mathcal{I}$ such that both $\blocky \cap \Delta$ and $\blocky \cap \overline{\Delta}$ are non-empty. Also, since $|\Delta| \leq \frac{1}{2}|\mathcal{V}|$, there exists $\blocky' \in \mathcal{I}$ such that $\blocky'\neq \blocky$ and $\blocky' \nsubseteq \Delta$. Choose  $\omega_0, \omega,\omega' \in \mathcal{V}$ as follows: $\omega_0 \in \blocky \cap \Delta, \omega \in \blocky \cap \overline{\Delta}$ and $\omega' \in \blocky' \cap \overline{\Delta}$. Since $\code$ is $G$-strongly incidence-transitive, there exists a permutation $h \in G_{\Delta,\omega_0}$ such that $\omega^h = \omega'$. Since $\omega_0$ is fixed by $h$ it follows that the part $\Sigma$ is fixed setwise by $h$. On the other hand $h$ moves $\omega \in \Sigma$ to $\omega' \in \Sigma'$, a contradiction. 
\end{proof}

We use the following general lemma to analyse potential $G$-strongly incidence-transitive codes, again concerning the action of a subgroup containing a codeword stabiliser.

\begin{lemma}
	\label{choosecarefully}
	Suppose $\Gamma$ is a $G$-strongly incidence-transitive code in a Johnson graph $J(\mathcal{V},k)$, and that $M$ is a subgroup of $G$ that contains the stabiliser $G_\Delta$ of a codeword $\Delta\in\Gamma$. Then for each $\omega\in\Delta$, the complement $\overline{\Delta}$ is contained in an $M_\omega$-orbit in $\mathcal{V}$, say $\Theta_M(\omega)$, and  $\overline{\Delta} \subseteq \bigcap_{\omega \in \Delta} \Theta_M(\omega)$.
\end{lemma}

\begin{proof}
By assumption $G_\Delta\leq M$, and hence $G_{\Delta,\omega}\leq M_\omega$. Since $\Gamma$ is  $G$-strongly incidence-transitive, $G_{\Delta,\omega}$ is transitive on 
$\overline{\Delta}$. Thus $\overline{\Delta}$ is contained in an $M_\omega$-orbit, say $\Theta_M(\omega)$, and the result follows.
\end{proof}

\section{Codes with reducible codeword stabilisers}
\label{sec:reducible}
Let $(V,B) = (\mathbb{F}_2^{2n},B)$ be a symplectic space with $n \geq 2$ and let $X \cong \syp_{2n}(2)$ be the isometry group of $B$. Let $U$ be a nontrivial proper subspace of $V$. Recall 
that $U$ is \emph{nondegenerate} if $U \cap U^\perp = \{0\}$ and \emph{totally-isotropic} if $U \leq U^\perp$. A $\calc_1$-subgroup of $X$ is the full setwise stabiliser of a nondegenerate or totally-isotropic subspace $U$ of $V$. If $U$ is totally-isotropic, or if $U$ is nondegenerate and $\dim(U) \neq n$, then the associated $\mathcal{C}_1$-subgroup is maximal in $X$ (see \cite[Table 3.5.C]{kl}).

\subsection{Nondegenerate subspaces}
\label{sec:nondegen}
A subspace $U\leq V$ is nondegenerate if and only if $V = U\oplus U^\perp$, and if $U$ is nondegenerate then $\dim(U)$ is even.  We let $\mathcal{Q}^\epsilon_U$ denote the set of all quadratic forms of type $\epsilon$ on $U$ which polarise to the restriction $B_U$ of $B$ to $U \times U$. If $U$ and $W$ are subspaces of $V$ such that $U \cap W = \{ 0 \}$, and if $B_U, B_W$ are symplectic forms on $U, W$, respectively,  then by $B_U \oplus B_W$ we mean the bilinear map $B'$ on $U\oplus W$ defined by $B'(u+w,u'+w') = B_U(u,u')+B_W(w,w')$ for all $u,u' \in U$ and $w,w'\in W$. It is straightforward to see that $B'$ is a symplectic form on $U\oplus W$ and that $U$ and $W$ are orthogonal relative to $B'$.  Similarly, we write $\phi_U \oplus \phi_W$ for the map $\phi':U\oplus W \to \mathbb{F}_2$ such that $\phi'(u+w) = \phi_U(u)+\phi_W(w)$ for all $u \in U$ and $w\in W$.
\begin{lemma}
\label{typeproduct}
Let $U$ be a nondegenerate subspace of $V$ and $W = U^\perp$. Let $(\phi_U,\phi_{W}) \in \mathcal{Q}^\epsilon_U\times\mathcal{Q}^{\epsilon'}_{W}$ and let $B_U$ and $B_{W}$ denote the polar forms of $\phi_U$ and $\phi_{W}$. Then $\phi_U \oplus \phi_W$ is a quadratic form of type $\epsilon\epsilon'$ on $V=U\oplus W$ which polarises to $B_U \oplus B_W$.
\end{lemma}
\begin{proof}
Let $\phi = \phi_U \oplus \phi_W$, let $B = B_U \oplus B_W$, and note that $\phi(0) = \phi_U(0)+\phi_W(0) = 0$. It is straightforward to check that $\phi$ is a quadratic form on $V=U\oplus W$. We show that $\phi$ polarises to $B$ as follows. For $i\in \{1,2\}$, let $x_i = u_i + w_i$, where $u_i \in U$ and $w_i \in W$. Then we have
	\begin{equation*}
	\begin{split}
		\varphi(x_1 + x_2) + \varphi(x_1) + \varphi(x_2) =& \varphi(u_1 + w_1 + u_2 + w_2) + \varphi(u_1+w_1)+\varphi(u_2 + w_2) \\
		=                                         & \varphi_U(u_1 + u_2) + \varphi_W(w_1 + w_2) + \phi_U(u_1) + \phi_W(w_1) \\
												 & + \phi_U(u_2) + \phi_W(w_2) \\
		=										 & B_U(u_1,u_2) + B_W(w_1,w_2) \\
		=										 & B(u_1+w_1,u_2+w_2) \\
		=                                         & B(x_1,x_2).
		\end{split}
	\end{equation*}
Therefore $\phi$ polarises to $B$. The fact that $\phi$ is of type $\epsilon \epsilon'$ follows from  \cite[Proposition 2.5.11]{kl}.
\end{proof}
\begin{lemma}\samepage
	\label{uniqueform}
	Let $V = \oplus_{i=1}^t V_i$ be a vector space equipped with a symplectic form  $B$, where each $V_i$ is nondegenerate and $V_i$ is orthogonal to $V_j$ for $i \neq j$. For each integer $i$ such that $1 \leq i \leq t$, let $\varphi_i$ be a quadratic form on $V_i$ which polarises to $B_{V_i}$. Then $\phi = \oplus_{i=1}^t \phi_i$ is the unique quadratic form on $V$ such that $\varphi_{V_i} = \varphi_i$ for $1 \leq i \leq t$, and $\phi$ polarises to $B$. Moreover $\phi$ is of type $\prod_{i=1}^t \epsilon_i$, where $\epsilon_i$ is the type $\phi_i$.
\end{lemma}

\begin{proof}
	We proceed by induction on $t$. If $t=1$ there is nothing to prove.  Let $t>1$ and suppose that the result is true for orthogonal direct sums of $t-1$ subspaces. Let $W = \oplus_{i=1}^{t-1} V_i$ and note that $V = W \oplus V_t$ and $W$ is orthogonal to $V_t$. By induction, $\varphi':= \oplus_{i=1}^{t-1}\phi_i$ is the unique quadratic form on $W$ such that $\varphi'_{V_i} = \varphi_i$ for all $i \leq t-1$, and $\phi'$ polarises to $B_W$, and $\phi'$ is of type $\prod_{i=1}^{t-1} \epsilon_i$. Let $\phi := \phi' \oplus \phi_t$. By Lemma \ref{typeproduct}, $\phi$ is a quadratic form  of type $\prod_{i=1}^t \epsilon_i$ on $V$, and $\phi$ polarises to $B = B_W \oplus B_{V_t}$, and the inductive hypothesis implies that $\phi_{V_i} = \phi_i$ for $1 \leq i \leq t$. 
	
	It remains to show that $\phi$ is unique. Let $\psi$ be a quadratic form on $V$ which polarises to $B$ such that $\psi_{V_i} = \phi_i$ for $1 \leq i \leq t$. Since $V = \oplus_{i=1}^t V_i$, each $x\in V$ has a unique expression $x = \sum_{i=1}^t x_i$ with $x_i \in V_i$. In particular, we have $w = \sum_{i=1}^{t-1} x_i\in W$ and $x=w+x_t$, and also $B(w,x_t)=0$ since, by assumption $W$ is orthogonal to $V_t$. Further, by the induction hypothesis (in particular the uniqueness of $\phi'$ mentioned in the previous paragraph), $\psi_W=\oplus_{i=1}^{t-1}\phi_i=\phi'$. Thus 
	\[
		\psi(x) =  \psi(w+x_t)	= \psi(w) + \psi(x_t) + B(w,x_t) = \varphi'(w) + \varphi_t(x_t) +0 = (\phi'\oplus \phi_t)(x) = \phi(x).
	\]
Therefore $\psi = \phi$, and $\phi$ is the unique quadratic form on $V$ which polarises to $B$ such that $\phi_i = \phi_{V_i}$ for $1 \leq i \leq t$.
\end{proof}

\begin{lemma}\label{lem:stabU}
Let $U$ be a nondegenerate subspace of $V$ and choose a basis for $V$ of the form $(b, b')$ where $b, b'$ are bases of $U, U^\perp$, respectively. Then, relative to this basis, the setwise stabiliser $X_U$ consists of all block diagonal matrices with diagonal $(A,B) \in \syp(U) \times \syp(U^\perp)$.
\end{lemma}

\begin{proof}
Let $b=(e_1,f_1,\dots,e_d,f_d)$ and $b'=(e_{d+1},f_{d+1},\dots,e_n,f_n)$  be ordered symplectic bases for $U$ and $U^\perp$, respectively, where $\dim(U)=2d$. The subgroup of $\gl_{2n}(2)$ which stabilises both $U$ and $U^\perp$ setwise is given by
\[
\left\{ \begin{pmatrix}R&O\\O&S\end{pmatrix} \bigg\rvert\, (R,S) \in \textstyle{\gl_{2d}(2)}\times\textstyle{\gl_{2(n-d)}(2)} \right\}.
\]
The Gram matrix $\mathcal{J}$ for $B$ is block diagonal with $n$ blocks of the form $\begin{pmatrix}
0&1\\
1&0
\end{pmatrix}$. 
Let $\mathcal{J}_U$ denote the Gram matrix for the restriction $B_U$, and similarly let $\mathcal{J}_{U^\perp}$ denote the Gram matrix for the restriction to $U^\perp$. Then $\mathcal{J}$ is the block diagonal matrix with blocks $(\mathcal{J}_U, \mathcal{J}_{U^\perp})$. Let $M \in X_U$, so $M$ leaves both $U$ and $U^\perp$ invariant, and hence is a block diagonal matrix with blocks $(R,S)$, say. Since $M \in \syp_{2n}(2)$ we also have $M\mathcal{J}M^T = \mathcal{J}$. Combining these conditions, we find
\begin{equation}
\label{c1ndcond}
	\begin{pmatrix}\mathcal{J}_U&O\\O&\mathcal{J}_{U^\perp}\end{pmatrix}
    =
	\begin{pmatrix}R&O\\O&S\end{pmatrix}
	\begin{pmatrix}\mathcal{J}_U&O\\O&\mathcal{J}_{U^\perp}\end{pmatrix}
	\begin{pmatrix}R^T&O\\O&S^T\end{pmatrix}
	=
	\begin{pmatrix}R\mathcal{J}_UR^T&O\\O&S\mathcal{J}_{U^\perp}S^T\end{pmatrix}.
	\end{equation}
Equation \eqref{c1ndcond} holds if and only if $(R,S) \in \syp_{2d}(2)\times \syp_{2(n-d)}(2)$.
\end{proof}

As a corollary to above three lemmas we compute the number of forms $\varphi\in\mathcal{Q}^\varepsilon$ which restrict on a given nondegenerate subspace  $U$ to a form  
$\varphi_U$ of given type.

\begin{corollary}\label{cor:kforcon}
    Let $U$ be a nondegenerate $2d$-dimensional subspace of $V$, where $1\leq d\leq n-1$. Then, for given $\varepsilon$, $X_U= \syp(U) \times \syp(U^\perp)$ has exactly two orbits in   $\mathcal{Q}^\varepsilon$, namely the subsets $\Delta(U,\varepsilon, \varepsilon')$, for $\varepsilon'\in\{+,-\}$,  of 
    forms $\varphi\in \mathcal{Q}^\varepsilon$ for which $\varphi_U$ is of type $\varepsilon'$ (and hence $\varphi_{U^\perp}$ is of type $\varepsilon\varepsilon'$). Moreover  
    \[
    k(n,d,\varepsilon,\varepsilon') := |\Delta(U,\varepsilon, \varepsilon')|= 2^{n-2}(2^d+\epsilon')(2^{n-d}+\epsilon\epsilon').
    \]
\end{corollary}

\begin{proof}
Let $\varphi\in \mathcal{Q}^\varepsilon$.
 Then, by Lemma \ref{uniqueform}, $\varphi$ is (uniquely) expressible as $\varphi =\varphi_U\oplus\varphi_{U^\perp}$, where $\varphi_U, \varphi_{U^\perp}$ are nondegenerate and polarise to $B_U, B_{U^\perp}$, respectively. Further, if
$\varphi_U$ is of type $\varepsilon'$, then $\varphi_{U^\perp}$ is of type $\varepsilon\varepsilon'$, by Lemma~\ref{typeproduct}. Now $\syp(U), \syp(U^\perp)$ 
act transitively in their respective Jordan-Steiner actions on $\mathcal{Q}^{\varepsilon'}_U$ and $\mathcal{Q}^{\varepsilon \varepsilon'}_{U^\perp}$,  and by Lemma~\ref{lem:stabU}, $X_U= \syp(U) \times \syp(U^\perp)$. It follows that $X_U$ is transitive on the set of pairs  
 $(\varphi_U,\varphi_{U^\perp}) \in \mathcal{Q}^{\epsilon'}_U\times \mathcal{Q}^{\epsilon\epsilon'}_{U^\perp}$, and hence $X_U$ is transitive on the subset $\Delta(U,\varepsilon, \varepsilon')$. Since each form in  $\mathcal{Q}^\varepsilon$ lies in $\Delta(U,\varepsilon, \varepsilon')$ for some $\epsilon'$, this proves that $X_U$ has two orbits in  $\mathcal{Q}^\varepsilon$. Moreover, for $\varphi\in \Delta(U,\varepsilon, \varepsilon')$,  the stabiliser 
 \[
 X_{U,\varphi}= \syp(U)_{\varphi_U} \times \syp(U^\perp)_{\varphi_{U^\perp}}
 \]
 and, applying the Orbit-Stabiliser Theorem, we have, 
\begin{equation*}
	\begin{split}
	k(n,d,\varepsilon,\varepsilon') &= |X_U:X_{U,\varphi}| 
	= \frac{|\syp_{2d}(2)|}{|\oo_{2d}^{\epsilon'}(2)|}\frac{|\syp_{2(n-d)}(2)|}{|\oo_{2(n-d)}^{\epsilon\epsilon'}(2)|} \\
	&= 2^{d-1}(2^d+\epsilon') \cdot 2^{n-d-1}(2^{n-d}+\epsilon\epsilon') \\
	&= 2^{n-2}(2^d+\epsilon')(2^{n-d}+\epsilon\epsilon').\qedhere
	\end{split}
\end{equation*}
\end{proof}

We use these results to show that all $X$-strongly incidence-transitive codes in $J(\mathcal{Q}^\varepsilon,k)$ for which the codeword stabilisers leave invariant a proper nondegenerate subspace of $V$ arise from Construction~\ref{ndcode}.

\begin{theorem}
\label{nondegencase}
Let $n\geq 2$, and $X=\syp_{2n}(2)$. Let $\Gamma$ be an $X$-strongly incidence-transitive code in $J(\mathcal{Q}^\varepsilon,k)$ for some $k$ satisfying $2 \leq k \leq |\mathcal{Q}^\varepsilon|-2$, and let $\Delta\in\Gamma$. Suppose that $X_\Delta\leq X_U$ for some nondegenerate  $2d$-dimensional subspace $U$ of $V$, where $1\leq d\leq  n-1$.  Then $\Gamma = \Gamma(n,d,\varepsilon,\varepsilon')$, for some $\varepsilon'\in\{ +, -\}$, $k= 2^{n-2}(2^d+\varepsilon')(2^{n-d}+\varepsilon\varepsilon')$, and $(n,d,\varepsilon)\neq (2,1,+)$,  as in Construction~\ref{ndcode}.
\end{theorem}


\begin{proof}
By Corollary~\ref{cor:kforcon}, $X_U$ has two orbits in $\mathcal{Q}^\varepsilon$, namely the sets $\Delta(U,\varepsilon, \varepsilon')$ defined there, for $\varepsilon'\in\{+,-\}$. 
Since $\Gamma$ is $X$-strongly incidence-transitive, the stabiliser $X_\Delta$ also has two orbits in $\mathcal{Q}^\varepsilon$, namely $\Delta$ and $\overline{\Delta}$. Hence, since $X_\Delta\leq X_U$, it follows that $\Delta =\Delta(U,\varepsilon, \varepsilon')$  for some $\varepsilon'$, and since $X_U$ leaves $\Delta(U,\varepsilon, \varepsilon')$ invariant we have equality $X_\Delta= X_U$. Since $\Delta= \Delta(U,\varepsilon, \varepsilon')$ it follows by Corollary~\ref{cor:kforcon} that $k=|\Delta|= 2^{n-2}(2^d+\varepsilon')(2^{n-d}+\varepsilon\varepsilon')$,
and since $2\leq k\leq |\mathcal{Q}^\varepsilon|-2$, we have $(n,d,\varepsilon)\neq (2,1,+)$. 
Thus $\Delta(U,\varepsilon, \varepsilon')$ is the codeword $\Delta(U)$ of the code $\Gamma(n,d,\varepsilon,\varepsilon')$ of Construction~\ref{ndcode}, and we conclude that $\Gamma=\Gamma(n,d,\varepsilon,\varepsilon')$, proving the theorem. 
%
\end{proof}

We  now  analyse the family of codes in Construction~\ref{ndcode}. 


\begin{theorem}\label{p:nd}
Let $n\geq 2$, $1\leq d\leq  n-1$, $X=\syp_{2n}(2)$, and let $\Gamma = \Gamma(n,d,\varepsilon,\varepsilon')$ in  $J(\mathcal{Q}^\varepsilon,k)$,  for   some $\varepsilon, \varepsilon'\in\{ +, -\}$ with  $(n,d,\varepsilon)\neq (2,1,+)$,  as in Construction~\ref{ndcode}, and let $A=\Aut(\Gamma)$.
 Then: 
\begin{enumerate}[(a)]
       \item For $\Delta(U)\in\Gamma$, with $U$ a  nondegenerate  $2d$-dimensional subspace of $V$, the stabiliser $X_U\leq X_{\Delta(U)}$, and equality holds unless $(d, \varepsilon) = (n/2, +)$;  in this exceptional case, $X_{\Delta(U)}= X_{\Delta(U^\perp)}=X_{\{U, U^\perp\}}\cong \syp_{n}(2)\wr S_2$, a maximal $\calc_2$-subgroup of $X$.
        \item The code $\Gamma$ is $X$-strongly incidence-transitive.
    \item The integer $k$ satisfies $k= 2^{n-2}(2^d+\varepsilon')(2^{n-d}+\varepsilon\varepsilon')$,
    and $3 \leq k \leq |\mathcal{Q}^\varepsilon|-3$. 
          \item  The code $\Gamma$ is self-complementary if and only if $(d,\epsilon) = (n/2,-)$.  

    \item The automorphism group $A=X$, unless $(d, \varepsilon) = (n/2, -)$, and in this exceptional case,  $A=X \times \langle s\rangle$, where $s$ is the complementing map, and for $\Delta(U)\in \Gamma$, $A_{\Delta(U)} = X_U\rtimes \langle xs\rangle\cong \syp_{n}(2)\wr S_2$ for some involution  $x\in X_{\{U, U^\perp\}}$ which interchanges $U$ and $U^\perp$.
\end{enumerate}
\end{theorem}

\begin{proof}
 Let $\Delta=\Delta(U)\in\Gamma$, so  $U$ is a  nondegenerate  $2d$-dimensional subspace of $V$. 

(a) Then  $V = U \oplus U^\perp$, and $X_U = X_{U^\perp} \cong \syp(U) \times \syp(U^\perp)$, by Lemma~\ref{lem:stabU}. It follows that  $X_U\leq X_\Delta$.  Using the information in \cite[Main Theorem and Table 3.5.C]{kl} if $n\geq 7$, and from \cite[Tables 8.28, 8.48, 8.64, and 8.80]{colva} for $2\leq n\leq 6$, we see that either $X_U$ is a maximal subgroup of $X$, and hence $X_U= X_\Delta$, or $d= n/2$ and the only proper subgroup of $X$ properly containing $X_U$ is the maximal $\calc_2$-subgroup $X_{\{U, U^\perp\}}= \syp_{2d}(2)\wr S_2$. In the latter case, where $d=n/2$, an element of $X_{\{U, U^\perp\}}$ which interchanges $U$ and $U^\perp$ maps $\Delta(U)$ to the codeword $\Delta(U^\perp)$ of $\Gamma(n,n/2,\varepsilon,\varepsilon\varepsilon')$. Moreover $\Delta(U)=\Delta(U^\perp)$ and $\Gamma=\Gamma(n,n/2,\varepsilon,\varepsilon\varepsilon')$ if and only if $\varepsilon=+$. Thus $X_\Delta=X_U$ if $\varepsilon=-$, and $X_\Delta=X_{\{U, U^\perp\}}$ if $\varepsilon=+$. This completes the proof of part (a). 

(b)  By part (a), $X_U\leq X_{\Delta(U)}$, and it follows from Corollary~\ref{cor:kforcon} that  $\Delta=\Delta(U,\varepsilon, \varepsilon')$ and the $X_U$-orbits in $\mathcal{Q}^\varepsilon$ are $\Delta$ and $\overline{\Delta}$.
By Lemma \ref{uniqueform}, for $\varphi\in\Delta$, we have $\varphi=\varphi_U\oplus\varphi_{U^\perp}$ with $\varphi_U\in\calq^{\varepsilon'}_U$ and $\varphi_{U^\perp}\in\calq^{\varepsilon\varepsilon'}_{U^\perp}$, while for $\psi\in \overline{\Delta}$, we have $\psi=\psi_U\oplus\psi_{U^\perp}$ with $\psi_U\in\calq^{-\varepsilon'}_U$ and $\psi_{U^\perp}\in\calq^{-\varepsilon\varepsilon'}_{U^\perp}$. Now, by \cite[Table 1, and the details in  (3.2.4e)]{maximalfactorisations}, 
the group $\syp(U)$ factorises as 
$\syp(U)={\rm GO}_{2d}^{\varepsilon'}(2){\rm GO}_{2d}^{-\varepsilon'}(2)$, so the stabiliser $\syp(U)_{\varphi_U}={\rm GO}_{2d}^{\varepsilon'}(2)$ is transitive on $\calq^{-\varepsilon'}_U$. Similarly the stabiliser in $\syp(U^\perp)$ of $\varphi_{U^\perp}$ is transitive on $\calq^{-\varepsilon\varepsilon'}_{U^\perp}$. It follows that $X_{U,\varphi}$ is transitive on 
$\calq^{-\varepsilon'}_U\times \calq^{-\varepsilon\varepsilon'}_{U^\perp}$, and hence $X_{U,\varphi}$ is transitive on $\overline{\Delta}$. Thus $X_U$, and hence also $X_\Delta$, is transitive on $\Delta\times \overline{\Delta}$. Since, by Witt's Lemma, $X$ is transitive on the set of nondegenerate $2d$-dimensional subspaces of $V$, it follows from Definition~\ref{defsit} that $\Gamma$ is $X$-strongly incidence-transitive, proving  part (b). 

(c)  By Corollary~\ref{cor:kforcon}, $k=k(n,d,\varepsilon, \varepsilon')= 2^{n-2}(2^d+\varepsilon')(2^{n-d}+\varepsilon\varepsilon')$.
Note that Construction~\ref{ndcode}  excludes the case  $(n,d,\varepsilon)= (2,1,+)$  (where $|\mathcal{Q}^\varepsilon|=10$, and   $(k, \varepsilon')=( 9,+)$ or $(1,-)$), and so in all cases 
it is easy to check that $k\geq 3$. We claim that  $k\leq |\mathcal{Q}^\varepsilon|-3=2^{n-1}(2^n+\epsilon)-3$.
Now $k$ satisfies
\[
k=k(n,d,\varepsilon, \varepsilon') \leq  2^{2n-2}+2^{n-2}(\varepsilon 2^d + 2^{n-d}) + \varepsilon 2^{n-2} =k(n,d,\varepsilon, +)
\]
and $\varepsilon 2^d+2^{n-d} \leq 2\varepsilon + 2^{n-1}$, since $1\leq d\leq n-1$. Hence  
\[
k\leq k(n,1,\varepsilon, +)=  2^{2n-2}+2^{n-2}(2\varepsilon + 2^{n-1}) + \varepsilon 2^{n-2}
= 3\cdot 2^{2n-3}+3\varepsilon\cdot 2^{n-2}.
\]
It is therefore sufficient to prove that $k(n,1,\varepsilon, +)\leq 2^{2n-1}+\epsilon2^{n-1}-3$, and this inequality is equivalent to $\varepsilon 2^{n-2}\leq 2^{2n-3}-3$. The latter clearly holds for all $n\geq 3$ and $\varepsilon\in\{+,-\}$. If $n=2$, then we must have $d=1$ and $\varepsilon=-$, and the latter inequality also holds. 
This establishes the claim, and part (c) is proved. 

(d) Suppose that $\Gamma$ is self-complementary. Then in particular  $k=|\calq^\epsilon|/2$, that is to say, $2^{n-2}(2^d+\epsilon')(2^{n-d}+\epsilon\epsilon')=2^{n-2}(2^n+\epsilon)$, which is equivalent to $2^d\epsilon+2^{n-d}=0$. The only possible parameters are $(d,\epsilon) = (n/2,-)$. Conversely suppose that $(d,\epsilon) = (n/2,-)$. Let $\Delta$ be a codeword, so  $\Delta=\Delta(U)$ for some $n$-dimensional nondegenerate subspace $U$. The elements of $\Delta$ are the forms $\phi \in \calq^\epsilon$ such that $\varphi_U$ has type $\varepsilon'$ and $\varphi_{U^\perp}$ has type $-\varepsilon'$.
Let $s$ be the complementing map.
Then $\Delta^s$ consists of all forms $\phi \in \calq^\epsilon$ such that $\varphi_U$ has type $-\varepsilon'$ and $\varphi_{U^\perp}$ has type $\varepsilon'$. Thus $\Delta^s=\Delta(U^\perp)$. Since $U^\perp$ is also an  $n$-dimensional nondegenerate subspace, it follows that $\Delta^s\in \Gamma$. Hence $\Gamma$ is self-complementary and $s$ is an automorphism of $\Gamma$. 

(e) From the definition of $\Gamma$ it is clear that $X\leq \aut(\Gamma)=A$. 
Assume first that $(d, \varepsilon) \neq (n/2, -)$. Then, as we showed in the proof of part (d), $k\neq |\mathcal{Q}^\epsilon|/2$, and so  $A$  is a subgroup of $\sym(\mathcal{Q}^\epsilon)$. 
Moreover, by Lemma \ref{le:noncomplete},  $\Gamma\neq \binom{\mathcal{Q}^\epsilon}{k}$, and so  $A$ is a proper subgroup of $\sym(\mathcal{Q}^\epsilon)$ and does not contain $\mathrm{Alt}(\mathcal{Q}^\epsilon)$.
It follows from the classification of the maximal subgroups of the finite symmetric and alternating groups that $A=X$, see \cite[Chapter 9]{maximalfactorisations}, especially Tables II--VI.

Now assume that $(d, \varepsilon) = (n/2, -)$, so that $k=|\mathcal{Q}^\epsilon|/2$. In this case the full automorphism group of $J(\mathcal{Q}^\varepsilon,k)$  is $\sym(\mathcal{Q}^\epsilon)\times \langle s\rangle$, and we showed in part (d) that $s\in A$. In the special case where $(n,d, \varepsilon) = (2,1, -)$, we have $|\mathcal{Q}^\epsilon|=6$, $k=3$, and $X=\syp_4(2)=\sym(\mathcal{Q}^\epsilon)$ is $3$-transitive on $\mathcal{Q}^\epsilon$. Hence 
$\Gamma= \binom{\mathcal{Q}^\epsilon}{3}$, and so  $A=\aut(J(\mathcal{Q}^\varepsilon,k))= \sym(\mathcal{Q}^\epsilon)\times \langle s\rangle=X\times \langle s\rangle$. 
Suppose now that $n\neq 2$. 
Then by Lemma \ref{le:noncomplete}, $\Gamma\neq \binom{\mathcal{Q}^\epsilon}{k}$, and so $A\cap \sym(\mathcal{Q}^\epsilon)$ is a proper subgroup of $\sym(\mathcal{Q}^\epsilon)$ and does not contain $\mathrm{Alt}(\mathcal{Q}^\epsilon)$.  As in the previous paragraph we conclude that $A\cap \sym(\mathcal{Q}^\epsilon)=X$. Since $s\in A$, it follows that $A=X \times \langle s\rangle$.

Finally, in the case $(d, \varepsilon) = (n/2, -)$, we consider $A_{\Delta(U)}$ for $\Delta(U)\in\Gamma$. Since $|A:X|=2$ it follows that $|A_{\Delta(U)}:X_{\Delta(U)}|\leq 2$, and also, by part (a), $X_{\Delta(U)}=X_U$. 
We showed in part (d) that $s$ interchanges $\Delta(U)$ and $\Delta(U^\perp)$. Consider the involution $x\in X_{\{U, U^\perp\}}$ defined by $(y, y')^x=(y',y)$, for all $(y,y')$ in the base group $Y=\syp_{n}(2)\times \syp_{n}(2)$ of  $X_{\{U, U^\perp\}}$. Then  $x$ interchanges $U$ and $U^\perp$, and also therefore $x$ interchanges  $\Delta(U)$ and $\Delta(U^\perp)$. Thus $xs$ fixes $\Delta(U)$ setwise, and since $s$ centralises $X$, it follows that $A_{\Delta(U)} = X_U\rtimes \langle xs\rangle\cong \syp_{n}\wr S_2$. 
\end{proof}

We complete this subsection by studying the other parameters of the codes in Construction~\ref{ndcode}.
The number of codewords is equal to $|X:X_{\Delta(U)}|$ since $X$ is transitive on codewords by Witt's Lemma, and by Theorem~\ref{p:nd}, $X_{\Delta(U)}=X_U$ if $(d, \varepsilon) \ne (n/2, +)$, and otherwise is a subgroup twice the size.  We can simply compute $|\Gamma|$ by calculating $\frac{|\syp_{2n}(2)|}{|\syp_{2d}(2)|\cdot |\syp_{2(n-d)}(2)|}$,
just as was done in \cite[Table 4.1.2]{BurnessGiudici}. Thus if $(d, \varepsilon) \ne (n/2, +)$, then
\begin{align*}
|\Gamma|&=4^{d(n-d)}\qbinom{n}{d}_4
\end{align*}
where $\qbinom{n}{d}_q$ denotes a Gaussian binomial coefficient, and if $(d, \varepsilon) = (n/2, +)$, then $|\Gamma|$ is half of this quantity.

Finding the minimum distance for these codes is more delicate.
For distinct codewords $\Delta(U_1), \Delta(U_2)\in\Gamma$, with $U_1, U_2$  nondegenerate  $2d$-dimensional subspaces of $V$, the distance between $\Delta(U_1), \Delta(U_2)$ is smallest when $\Delta(U_1)\cap\Delta(U_2)$ is as large as possible. So we want to maximise the number of forms $\varphi\in\calq^\varepsilon$ such that both $\varphi_{U_1}$ and $\varphi_{U_2}$ have type $\varepsilon'$. In general, the subspaces $U_1\cap U_2$ and $U_1^\perp\cap U_2^\perp$ could be far from nondegenerate, and we have not analysed the complete situation. 
Computational evidence for small values of $n$ and $d$ suggests that the pairs of codewords analysed in the following lemma may provide examples of codeword pairs at minimum distance.


\begin{lemma}\label{le:distanceND}
    Let $n\geq 3$, $1\leq d\leq  n-1$, and $\varepsilon, \varepsilon'\in\{ +, -\}$, and let $\Gamma = \Gamma(n,d,\varepsilon,\varepsilon')$ as in Construction~\ref{ndcode} in the graph $J(\mathcal{Q}^\varepsilon,k)$, with $k$ as in Construction~\ref{ndcode}. Then 
    \begin{enumerate}[(a)]
        \item there exist pairs of codewords at distance $2^{2n-4}$; and 
        \item provided $(d, \varepsilon)\ne (n/2, +)$ there are also pairs of codewords at distance 
        \[
        2^{n-3}(2^{n-|n-2d|}-1)(2^{|n-2d|}-\epsilon).
        \] 
    \end{enumerate} 
\end{lemma}
\begin{proof} 
Recall that  if $k$-subsets $\Delta_1,\Delta_2$ are at distance $D$ in $J(\mathcal{V},k)$, then the complements $\overline{\Delta_1},\overline{\Delta_2}$ are also at distance $D$ in $J(\mathcal{V},|\mathcal{V}|-k)$.
Since  the code $\Gamma(n,d,\epsilon,-\epsilon')$ is the complementary code of  $\Gamma(n,d,\epsilon,\epsilon')$, we may also assume that $\epsilon'=+$. Thus for a nondegenerate subspace $U$ of dimension $2d$ the codeword $\Delta(U)$ is the set of quadratic forms $\phi\in\mathcal{Q}^\varepsilon$ such that $\phi_U$ is of $+$ type (and hence, by Lemma~\ref{uniqueform}, $\varphi_{U^\perp}$ has type $\epsilon$).

For distinct $\Delta(U_1), \Delta(U_2)\in\Gamma$, with $U_1, U_2$  nondegenerate  $2d$-dimensional subspaces of $V$, the distance between $\Delta(U_1), \Delta(U_2)$ is equal to 
\[
k-|\Delta(U_1)\cap\Delta(U_2)|=|\Delta(U_1)|-|\Delta(U_1)\cap\Delta(U_2)|=|\Delta(U_1)\setminus\Delta(U_2)|. 
\]
This is equal to the number of forms $\varphi\in\calq^\varepsilon$ such that  $\varphi_{U_1}$ has $+$ type  and $\varphi_{U_2}$ has $-$ type. 

(a)   Assume first that $d\leq n/2$   Let $W=\langle e_1, f_1, e_2, f_2\rangle$, $W'=\langle e_3, f_3,\ldots, e_{d+1}, f_{d+1}\rangle$,  $W''=\langle e_{d+2}, f_{d+2},\ldots, e_n, f_n\rangle$, where the $e_i$ and $f_i$ form a symplectic basis of $V$, so that $V$ admits an orthogonal decomposition $V=W\oplus W'\oplus W''$ with  $W, W', W''$ nondegenerate subspaces of dimensions $4,2(d-1), 2(n-d-1)$, respectively.  Since $d\leq n/2$ and $n\geq3$, the subspace $W''$ is non-trivial for all $d$, but $W'$ is trivial if $d=1$. We will compute the distance between  $\Delta(U_1)$ and  $\Delta(U_2)$ for $U_1=\langle e_1, f_1\rangle\oplus W'$ and $U_2=\langle e_1, f_1+f_2\rangle\oplus W'$ (both nondegenerate  $2d$-dimensional subspaces). 
        
       
       Let $\phi$ be as in \eqref{eq-forms}.
 This implies that $\phi_{U_1}$, $\phi_{U_2}$, and $\phi_{W\oplus W'}$ are all of $+$ type, by Corollary~\ref{formvectorbijection}, $\mathcal{Q}^\varepsilon$ comprises all the forms $\varphi_c$ with $c\in\sing(\varphi)$, and hence 
 \[
 \Delta(U_1)=\{\phi_c\mid c\in\sing(\varphi), (\phi_c)_{U_1}\text{ is of } + \text{ type } \}. 
 \] 
Let $c\in\sing(\varphi)$ and, for $i=1, 2$, write $c=c_i'+c_i''$ where $c_i'\in U_i, c_i''\in U_i^\perp$. By \eqref{parameq2},  $(\phi_c)_{U_i}$ depends only on $c_i'$, and by Corollary~\ref{formvectorbijection} applied to $\mathcal{Q}^+_{U_i}$, $(\phi_c)_{U_i}$ has $+$ type if and only if $c_i'\in\sing(\varphi)\cap U_i$. Now write  $c=w+w'+w''$ where $w=c_1e_1+d_1f_1+c_2e_2+d_2f_2\in W, w'\in W', w''\in W''$. 
Then $c_1'= c_1e_1+d_1f_1+w'$ and $c_2'=(c_1+c_2)e_1+d_1(f_1+f_2)+w'$.
Thus, for $c\in\sing(\varphi)$,  $\varphi_c\in \Delta(U_1)\setminus \Delta(U_2)$ if and only if $c_1'\in\sing(\varphi)$ (so that $\varphi_c\in\Delta(U_1)$) and $c_2'\not\in\sing(\varphi)$ (so that $\varphi_c\not\in\Delta(U_2)$), that is to say,
\begin{align*}
\varphi(c)&= \phi(c_1e_1+d_1f_1+c_2e_2+d_2f_2+w'+w'')=0\\
\varphi(c_1')&= \phi(c_1e_1+d_1f_1+w')=0\\
\varphi(c_2')&= \phi((c_1+c_2)e_1+d_1(f_1+f_2)+w')=1.
\end{align*}
Applying the definition of $\varphi$ in \eqref{eq-forms}, these conditions become:
\[
c_1d_1+c_2d_2+\phi(w')+\phi(w'')=0, \ c_1d_1+\phi(w')=0, (c_1+c_2)d_1+\phi(w')=1,
\]
which in turn are equivalent to $\phi(w'')=c_2d_2$, $\phi(w')=c_1d_1$, and $c_2d_1=1$. So we must have $c_2=d_1=1$, and  for arbitrary $w'\in W', w''\in W''$ we must have $c_1=\phi(w')$ and $d_2=\phi(w'')$. Thus the number of choices for $c$ is therefore $|W'|\cdot|W''|=2^{2n-4}$, and hence the distance between $\Delta(U_1)$ and $\Delta(U_2)$ is $2^{2n-4}$ (which we observe is independent of $d$).
       

This proves part (a) for $d\leq n/2$. Since $ \Gamma(n,d,\epsilon,\epsilon')=\Gamma(n,n-d,\epsilon,\epsilon\epsilon')$, if $d > n/2$, the code  $ \Gamma(n,d,\epsilon,\epsilon')$ also has codewords at distance $2^{2n-4}$.

    (b)  Assume first that $d< n/2$, and let $U_1,U_2$ be  nondegenerate $2d$-dimensional subspaces such  that $U_2\subseteq U_1^\perp$. Then  
 $V$ admits an orthogonal decomposition $V=U_1\oplus U_2\oplus W$ with $W=(U_1\oplus U_2)^\perp\ne 0$, so $U_1, U_2, W$ have dimensions $2d, 2d, 2(n-2d)$, respectively. For $\varphi\in\mathcal{Q}^\varepsilon$, it follows from  Lemma~\ref{uniqueform} that $\varphi= \varphi_{U_1}\oplus \varphi_{U_2}\oplus \varphi_{W},$ and   
as discussed above $|\Delta(U_1)\setminus\Delta(U_2)|$ is the number of forms $\varphi\in\calq^\varepsilon$ such that  $\varphi_{U_1}$ has type $\varepsilon'=+$ and $\varphi_{U_2}$ has type $-$. For such forms, $\phi_W$ has type $-\epsilon$, by Lemma~\ref{uniqueform}.
Thus $|\Delta(U_1)\setminus\Delta(U_2)|$ is the number of forms $\varphi\in\calq^\varepsilon$ such that  $\varphi_{U_1}, \varphi_{U_2}, \varphi_{W}$ have types $+, -, -\varepsilon$, respectively. Hence 
\[
|\Delta(U_1)\setminus \Delta(U_2)|=2^{d-1}(2^d+1)\cdot 2^{d-1}(2^d-1)\cdot 2^{n-2d-1}(2^{n-2d}-\epsilon)=2^{n-3}(2^{2d}-1)(2^{n-2d}-\epsilon).
\]
This proves part (b) when $d<n/2$ since $|n-2d|= n-2d>0$.

If $d>n/2$, then using the fact that  $ \Gamma(n,d,\epsilon,\epsilon')=\Gamma(n,n-d,\epsilon,\epsilon\epsilon')$, we get pairs of codewords at distance $2^{n-3}(2^{2(n-d)}-1)(2^{2d-n}-\epsilon)$, which matches the statement for $n-2d<0$.
Finally we consider the case $n=2d$.  Then $W$ is trivial in the decomposition above and $U_1^\perp=U_2$. If $\epsilon=+$ then  $\Delta(U_1)= \Delta(U_2)$, the same codeword,  while if $\epsilon=-$ then  $\Delta(U_1)$ and $\Delta(U_2)$ are disjoint so are at distance $k=2^{n-2}(2^d+1)(2^{n-d}-1)=2^{n-2}(2^{n}-1)$, which matches the statement for $n-2d=0$.
\end{proof}

Unfortunately there are too many orbits of $\syp_{2n}(2)$ on pairs of nondegenerate $2d$-dimensional subspaces of dimension, if $d$ is large, to be able to do a complete analysis. Thus we were not able to compute the minimum distance for this family of codes, but Lemma~\ref{le:distanceND} allows us to derive an upper bound. In fact we discovered computationally, using Magma \cite{Magma}, that the upper bound given in Corollary~\ref{cor:distanceND} is the correct minimum distance for all $n, d, \varepsilon, \varepsilon'$ with  $n\leq 5$. We conjecture  that the minimum distance of $\Gamma(n,d,\epsilon,\epsilon')$  is equal to this upper bound in general.

\begin{corollary}\label{cor:distanceND}
 Let $n\geq 2$, $1\leq d\leq  n-1$, and $\varepsilon, \varepsilon'\in\{ +, -\}$ such that $(n,d,\varepsilon)\neq (2,1,+)$. 
 Then the minimum distance $D$ of $\Gamma(n,d,\varepsilon,\varepsilon')$ as in Construction~\ref{ndcode} satisfies: 
 \[
D \leq \left\{ \begin{array}{ll}
     2^{2n-4}-2^{n-3}   & \mbox{if $|n-2d|= 1$ and $\epsilon=+$}, \\
     2^{2n-4}           & \mbox{otherwise.}
\end{array}     \right.
\]
\end{corollary}

\begin{proof}
 Note that $\Gamma(2,1,-,+)=\Gamma(2,1,-,-)=\binom{\mathcal{Q}^\varepsilon}{3}$, so the minimum distance is $1=2^{2n-4}$, so we may assume that $n\geq3$.
  By Lemma \ref{le:distanceND}(a), there are pairs of codewords at distance $2^{2n-4}$, so $D\leq 2^{2n-4}$.
 
 Suppose first that $d<n/2$. Then by Lemma \ref{le:distanceND}(b), there are also pairs of codewords at distance  $D':=2^{n-3}(2^{2d}-1)(2^{n-2d}-\epsilon)= 2^{n-3}(2^{n}-2^{n-2d}-\varepsilon 2^{2d}+\varepsilon)$. 
Now we compare this distance to the upper bound $2^{2n-4}$: the inequality $D'<2^{2n-4}$ 
 can be rewritten as $2^{n-1}+\epsilon<2^{n-2d}+\epsilon 2^{2d}$.
 Since $d\geq 1$, this never holds if $\epsilon=-$.
 If $\epsilon=+$ the condition can be written as $2^{n-1}-2^{n-2d}=2^{n-2d}(2^{2d-1}-1)< 2^{2d}-1=2(2^{2d-1}-1)+1$ and we see that this holds exactly when $n-2d=1$. In this case $D'= 2^{2n-4}-2^{n-3}$ is an improved upper bound, as in the statement.
 
 For $d>n/2$, using the fact that $ \Gamma(n,d,\epsilon,\epsilon')=\Gamma(n,n-d,\epsilon,\epsilon\epsilon')$, we also get this improved upper bound when $n-2(n-d)=1$ that is, when $n-2d=-1$.
\end{proof}


 



\subsection{Totally-isotropic subspaces}
\label{sec:ti}
Let $(V,B) = (\mathbb{F}_2^{2n},B)$ be a symplectic space with $n \geq 2$ and let $X \cong \syp_{2n}(2)$ be the isometry group of $B$. 
Let $U$ be a nontrivial proper totally-isotropic subspace of $V$ and set $d = \dim(U)$, so $d\leq n$. Let $X_U$ denote the setwise stabiliser of $U$. Throughout Section \ref{sec:ti} we use the notation $\mathcal{Q}^\varepsilon _{d} := \lbrace \varphi \in \mathcal{Q}^\varepsilon \mid U \subseteq \sing(\varphi) \rbrace$  and $\mathcal{Q}^\varepsilon _{d-1}  := \lbrace \varphi \in \mathcal{Q}^\varepsilon \mid U \nsubseteq \sing(\varphi) \rbrace$. Clearly $\mathcal{Q}^\varepsilon$ is the union of these two sets. If $\phi \in \mathcal{Q}^-$ then the dimension of a $\phi$-singular subspace is at most $n-1$, and in this case, if $\dim(U)=n$, then $\mathcal{Q}^\varepsilon_n =\emptyset$ and $\mathcal{Q}^\epsilon = \mathcal{Q}^\epsilon_{n-1}$. If $(d,\epsilon) \neq (n,-)$ then both $\mathcal{Q}^\epsilon_d$ and $\mathcal{Q}^\epsilon_{d-1}$ are nonempty.
\begin{lemma}
	\label{intersectionthm}
Let $U$ be a nontrivial totally-isotropic $d$-dimensional subspace of $V$.
\begin{enumerate}[(a)]
\item If $\varphi \in \mathcal{Q}$ then $\sing(\varphi)\cap U$ is a subspace of $U$. If $U$ is $\varphi$-singular then the dimension of $\sing(\varphi)\cap U$ is $d$, otherwise the dimension is $d-1$. In particular, for each $\varepsilon\in\{+, -\}$ and each $r\in\{d-1,d\}$, $\mathcal{Q}_r^\varepsilon = \lbrace \varphi \in \mathcal{Q}^\varepsilon \mid \dim(U \cap \sing(\varphi))=r \rbrace$.

\item If $\varphi \in \mathcal{Q}_d^\varepsilon$ and $c \in V$ then $\varphi_c \in \mathcal{Q}^\varepsilon_d$ if and only if $c \in \sing(\varphi)\cap U^\perp$, and $\varphi_c \in \mathcal{Q}^\varepsilon_{d-1}$ if and only if $c \in \sing(\varphi)\setminus(\sing(\varphi)\cap U^\perp)$. 

\item $X_U$ acts transitively on $\mathcal{Q}^\varepsilon _{d}$.
\end{enumerate}
\end{lemma}

\begin{proof}\  
(a) Let $u,v \in U$. Since $U \leq U^\perp$ we have
	\begin{equation}
	\label{addingint}
	\varphi(u+v) = \varphi(u) + \varphi(v) + B(u,v) = \varphi(u) + \varphi(v).
	\end{equation}
	If $u,v \in \sing(\varphi)\cap U$ then Equation \eqref{addingint} implies $\varphi(u+v) = 0$ and therefore $\sing(\varphi)\cap U$ is a subspace. Moreover, Equation \eqref{addingint} implies that the restriction $\varphi_U:U \to \mathbb{F}_2$  is a linear map and therefore $\dim(\varphi(U)) = 0$ or $1$, depending on whether or not $U$ is $\varphi$-singular. Then by the Rank-Nullity Theorem, $\dim(\sing(\varphi)\cap U) = d$ or $d-1$.

(b) Since $\varphi\in \mathcal{Q}^\epsilon_d$ we have $U\subseteq \sing(\varphi)$ by part (a). Corollary \ref{formvectorbijection} implies that $\phi_c \in \mathcal{Q}^\epsilon$ if and only if $c \in \sing(\varphi)$.   From the definition of $\mathcal{Q}^\epsilon_d$ and Lemma \ref{lem_singrelations} we have
\begin{align*}
\varphi_c \in \mathcal{Q}^\varepsilon_d &\Leftrightarrow c \in \sing(\varphi)\ \mbox{and}\ U \subseteq \sing(\varphi_c)\\
    &\Leftrightarrow c \in \sing(\varphi)\ \mbox{and}\ U \subseteq \sing(\varphi)+c \\
    &\Leftrightarrow c \in \sing(\varphi)\ \mbox{and}\  U+c \subseteq \sing(\varphi).     
\end{align*}  
However $U+c \subseteq \sing(\varphi)$ if and only if 
\begin{equation}
\label{conditiony}
\varphi(u+c) = 0 \text{ for all } u \in U.
\end{equation}
Expanding Equation \eqref{conditiony} using the polarisation identity and using the fact that $\varphi(u) = \varphi(c) = 0$, we have $\varphi(u+c) = B(u,c)$ for all $u\in U$. Therefore $\varphi_c \in \mathcal{Q}^\varepsilon_d$ if and only if $c \in \sing(\varphi) \cap U^\perp$. The second assertion now follows from part (a). 

(c) Assume $(d,\epsilon) \neq (n,-)$ so that $\mathcal{Q}_d^\epsilon$ is nonempty.
It is shown in \cite[Section 7.7]{dixonmortimer} that for each $c \in V$, the symplectic transvection $\tau_c: V \rightarrow V$ defined by $v\tau_c = v + B(v,c)c$ is an element of $X$ and an involution. Let $\phi,\phi_c \in \mathcal{Q}^\epsilon_d$. We show that $\phi$ can be mapped to $\phi_c$ by an element of $X_U$. By part (b), $c \in U^\perp$, so for all $u \in U$ we have $u \tau_c = u + B(u,c)c = u$. Therefore $\tau_c$ fixes $U$ pointwise, and $\tau_c\in X_U$. Since $\phi(c)=0$ and $x^2=x$ for each $x\in \F_2$, we have, for all $v\in V$,
	\begin{align*}
		\varphi^{\tau_c}(v) &= \varphi(v+B(v,c)c) \\ \nonumber
		&= \varphi(v) + B(v,c)^2\varphi(c) + B(v,c)^2 \\ \nonumber
		&= \varphi(v) + B(v,c)^2 \\ \nonumber
		&= \varphi_c(v). \nonumber
	\end{align*}
	Therefore $\varphi^{\tau_c}=\varphi_c$, and as $\varphi_c$ was an arbitrary element of $\mathcal{Q}^\epsilon_d$, it follows that  $X_U$ is transitive on $\mathcal{Q}^\varepsilon _d$.
\end{proof}
Lemma \ref{intersectionthm}(a) shows $\mathcal{Q}^\varepsilon_{d}$ and $\mathcal{Q}^\varepsilon_{d-1}$ are nonempty unless $(d,\varepsilon) = (n,-)$; the latter case is dealt with in Section \ref{transc1}. After discussing the structure of $X_U$, we prove in Corollary \ref{cor:transQd-1} that $X_U$ also acts transitively on $\mathcal{Q}^\varepsilon_{d-1}$.
\subsubsection{The structure of \texorpdfstring{$X_U$}{XU}}
\begin{lemma}[{\cite[Lemma 4.1.12]{kl}}]
	\label{levilemma}
	Let $(V,B) = (\mathbb{F}_2^{2n},B)$ be a symplectic space with full isometry group $X$ and let $U$ be a totally-isotropic $d$-dimensional subspace of $V$. Then there exist subspaces $U'$ and $W$ of $V$ such that the following hold:
	\begin{enumerate}[(a)]
		\item $U'$ is totally-isotropic of dimension $d$, $U \oplus U'$ is nondegenerate and $W = (U \oplus U')^\perp$,
		\item $X_U = R \rtimes L$, where $L$ fixes setwise each of the subspaces $U$, $U'$ and $W$, and $R$ acts trivially on the spaces $U$, $U^\perp /U$ and $V/U^\perp$,
		\item $L \cong \gl(U) \times \syp (W)$.
	\end{enumerate}
\end{lemma}

The semidirect product $X_U = R \rtimes L$ in Lemma \ref{levilemma}(b) is called a \emph{Levi decomposition} with \emph{Levi component} $L$ and \emph{unipotent radical} $R$.

Throughout Section \ref{sec:ti} we make use of the standard symplectic basis ordered as follows
\begin{equation}
	\label{basis}
	\mathscr{B} = \lbrace e_1, \dots, e_d, e_{d+1},f_{d+1}, \dots,e_n,f_n,f_1,\dots,f_d \rbrace,
\end{equation}
(which, if $d=n$, is $\{e_1,\dots, e_n,f_n\dots,f_1\}$).
If $x,x' \in V$ are given by $x = \sum_{i=1}^n (x_ie_i + y_i f_i)$ and $x' = \sum_{i=1}^n (x_i'e_i + y_i' f_i)$ then
\begin{equation}
\label{Bdef}B(x,x') = \sum_{i = 1}^n x_{i}y_{i}' + y_{i}x_{i}'.\end{equation}
The Gram matrix for $B$, relative to the ordered basis $\mathscr{B}$,  is therefore
\begin{equation}
\label{tigram}
\mathcal{J} = \left( \begin{array}{ccc}
				0   & 0          & I_d \\
				0   & J & 0   \\
				I_d & 0          & 0
			\end{array} \right)
\end{equation}
	where $I_d$ is a $d \times d$ identity matrix and $J$ is the $2(n-d) \times 2(n-d)$ block-diagonal matrix with blocks $\begin{pmatrix}
			0 & 1 \\ 1 & 0
		\end{pmatrix}$, and $J$ is an empty matrix if $d=n$.
Since $X_U$ acts transitively on the totally-isotropic $d$-dimensional subspaces of $V$, we may set $U = \langle e_i \mid 1 \leq i \leq d \rangle$, $U' = \langle f_i \mid 1 \leq i \leq d \rangle$ and $W = \langle e_i, f_i \mid d+1 \leq i \leq n \rangle$ throughout Section \ref{sec:ti}. Note that $W=\{0\}$ if $d=n$.
%
\begin{lemma}
	\label{calculatestabiliser}
For each integer $d$ satisfying $1\leq d \leq n$, the setwise stabiliser $X_U$ of the totally-isotropic space $U = \langle e_1, \cdots, e_d \rangle$ is the group of all $2n \times 2n$ block matrices of the form
	$$g= \left(
		\begin{array}{ccc}
				P^{-T} & 0 & 0 \\
				Y & N & 0 \\
				Q & Z & P
			\end{array}
		\right) $$
	which satisfy the following conditions:
	\begin{enumerate}[(i)]
		\item $P \in \gl_d(2)$,
		\item $N \in \isom(B_W) \cong \syp_{2(n-d)}(2)$ (and is an empty matrix if $d=n$),
		\item $QP^T + P Q^T = Z J Z^T$ if $d<n$, while $QP^T + P Q^T = 0$ if $d=n$,    and
		\item $P Y^T = Z JN^T$ if $d<n$, while $Y, Z$ are empty matrices if $d=n$.
	\end{enumerate}
In particular, $Y$ is determined uniquely by $Z$, $N$ and $P$.
\end{lemma}

\begin{proof}
Let $g\in \syp_{2n}(2)$. If $g$ fixes $U$ then it also fixes $U^\perp = \langle e_1, \cdots, e_d,  e_{d+1}, f_{d+1}, \cdots, e_{n}, f_{n} \rangle$. Expressing $g$ with respect to the ordered basis defined in Equation \eqref{basis} we have
	$$ g = \left(
		\begin{array}{ccc}
				M & 0 & 0 \\
				Y & N & 0 \\
				Q & Z & P
			\end{array}
		\right)$$
	where $M$ and $P$ lie in $\gl_d(2)$, $N$ lies in $\gl_{2(n-d)}(2)$ (and is an empty matrix if $d=n$), and $Q$, $Y$ and $Z$ are respectively $d\times d$, $2(n-d)\times d$ and $d \times 2(n-d)$ matrices. By definition, $g \in \syp_{2n}(2)$ if and only if $g\mathcal{J}g^T = \mathcal{J}$, where $\mathcal{J}$ is as in Equation \eqref{tigram}. We have
	\begin{align*}
		g\mathcal{J}g^T = & \left(
		\begin{array}{ccc}
				M & 0 & 0 \\
				Y & N & 0 \\
				Q & Z & P
			\end{array}
		\right)
		\left(
		\begin{array}{ccc}
				0   & 0 & I_d \\
				0   & J & 0   \\
				I_d & 0 & 0
			\end{array}
		\right)
		\left(
		\begin{array}{ccc}
				M^T & Y^T & Q^T \\
				0   & N^T & Z^T \\
				0   & 0   & P^T
			\end{array}
		\right)          \\
		=       &
		\left(
		\begin{array}{ccc}
				0 & 0   & M \\
				0 & N J & Y \\
				P & Z J & Q
			\end{array}
		\right)
		\left(
		\begin{array}{ccc}
				M^T & Y^T & Q^T \\
				0   & N^T & Z^T \\
				0   & 0   & P^T
			\end{array}
		\right)          \\
		=       &
		\left(
		\begin{array}{ccc}
				0    & 0                     & M P^T                        \\
				0    & N J N^T      & N J Z^T+YP^T       \\
				PM^T & P Y^T+Z JN^T & P Q^T+QP^T+ Z JZ^T
			\end{array}
		\right).
	\end{align*}
	Invoking the condition $\mathcal{J} = g\mathcal{J}g^T$, we have
	\[
		\left(
		\begin{array}{ccc}
				0   & 0            & I_d \\
				0   & J & 0   \\
				I_d & 0            & 0
			\end{array}
		\right) =
		\left(
		\begin{array}{ccc}
				0    & 0                     & M P^T                        \\
				0    & N J N^T      & N J Z^T+YP^T       \\
				PM^T & P Y^T+Z JN^T & P Q^T+QP^T+ Z JZ^T
			\end{array}
		\right).
	\]
	Therefore $g \in \syp_{2n}(2)$ if and only if $M=P^{-T}$ and conditions (i)-(iv) hold.
\end{proof}
\begin{corollary}
	\label{levi}
Let $d$ be an integer satisfying $1\leq d \leq n$ and set $U = \langle e_1, \cdots, e_d \rangle$. Then $X_U=R\rtimes L$ is a Levi decomposition with Levi component $L$ and unipotent radical $R$ given by
	\begin{align*}
		L & = \left\lbrace  
		\begin{pmatrix}
				P^{-T} & 0 & 0      \\
				0 & N & 0      \\
				0 & 0 & P
		\end{pmatrix} : P \in \textstyle{\gl_d(2)}, N \in \textstyle{\syp_{2(n-d)}(2)} \right\rbrace \\
		\\
		R & = \left\lbrace  \begin{pmatrix}
				I_d & 0          & 0   \\
			J Z^T   & I_{2(n-d)} & 0   \\
				Q   & Z          & I_d
		\end{pmatrix} : Q + Q^T = Z J Z^T  \right\rbrace
	\end{align*}
	where if $d=n$, then $N,  Z$ are empty matrices and the defining condition on $Q$ is $Q+Q^T=0$ ($Q$ is symmetric).
\end{corollary}

\begin{proof}
	We have $U = \langle e_1, \dots, e_d \rangle$, $U' = \langle f_1, \cdots , f_d \rangle$ and $W = \langle e_{d+1}, f_{d+1} , \cdots , e_n, f_n \rangle$. It is clear that property (a) of Lemma \ref{levilemma} holds. 
	
	Note that $L$  fixes setwise each of the subspaces $U$, $U'$ and $W$, and $R$  acts trivially on the spaces $U$, $U^\perp /U$ and $V/U^\perp$.
	
	Suppose $g \in X_U$. Then $g$ has the form of the matrix in Lemma \ref{calculatestabiliser} with entries $ N, P, Q, Y, Z$ satisfying the conditions of Lemma~\ref{calculatestabiliser}(i)--(iv). Further, we can factorise $g$ as follows: 
	\begin{align*}
		g = \left(
		\begin{array}{ccc}
				P^{-T} & 0 & 0 \\
				Y & N & 0 \\
				Q & Z & P
			\end{array}
		\right) = \left(
		\begin{array}{ccc}
				I_d      & 0          & 0   \\
				Y P^{T} & I_{2(n-d)} & 0   \\
				Q P^{T} & Z N^{-1}   & I_d
			\end{array}
		\right)
		\left(
		\begin{array}{ccc}
				P^{-T} & 0 & 0      \\
				0 & N & 0      \\
				0 & 0 & P
			\end{array}
		\right) = hh'.
		\end{align*}
    The second factor $h'$ lies in $L$. The first factor $h\in R$ if and only if $YP^T=J(ZN^{-1})^T$ and $(QP^{T})+(QP^{T})^T=(ZN^{-1})J(ZN^{-1})^T$. We check these conditions  as follows, using the conditions from  Lemma~\ref{calculatestabiliser}. 
    If $d=n$ then $Y, Z, N$ are empty matrices, and the matrix $QP^{T}$ 
    satisfies $(QP^{T})+(QP^{T})^T =QP^T+PQ^T$, which is a zero matrix by Lemma~\ref{calculatestabiliser}(iii). Thus $h'\in R$ in this case. Suppose now that $d<n$ and  
    note that $J=J^T$. Also  $N$, and hence also $N^{-1}$, lie in $\isom(B_W)$, so we have $N^{-1}JN^{-T}=J$. Now 
        	\begin{align*}
    J(ZN^{-1})^T &=J^T(ZN^{-1})^T = (ZN^{-1}J)^T&\\
    &= (ZJN^T)^T \       &\text{(since $N^{-1}JN^{-T}=J$)}\\ 
    &= (PY^{T})^T &\text{(by Lemma~\ref{calculatestabiliser}(iv))}  \\
    &= YP^T &
    \end{align*}
    so the first condition holds, and 
    	\begin{align*}
    (QP^{T})+(QP^{T})^T &=QP^T+PQ^T &\\
    &= ZJZ^T\       &\text{(by Lemma~\ref{calculatestabiliser}(iii))}\\ 
    &= Z (N^{-1} J N^{-T}) Z^T &\\
    &=  (ZN^{-1})J(ZN^{-1})^T&
    \end{align*}
    proving the second condition, so the matrix $h$ lies in $R$. 
Thus $X_U \subseteq RL$. It is easy to check that each matrix in $R$ or $L$ satisfies all the conditions of Lemma~\ref{calculatestabiliser} and hence both $R$ and $L$ are contained in $X_U$. Therefore $X_U = RL$.  Finally, $R \cap L = \lbrace I_{2n} \rbrace$, and  $R \triangleleft X_U$, so $X_U \cong R \rtimes L$.
\end{proof}
Recall the definition of $\phi_0^\epsilon$ in Equation \eqref{eq-forms}.
Provided that $(d,\varepsilon) \neq (n,-)$, we have $\varphi_0^\varepsilon \in \mathcal{Q}^\varepsilon_d$ since $U\subset\sing(\phi_0^\epsilon)$.

Our next goal is to provide a condition which is used in Lemma \ref{cosettrans} to show that elements of $L$ fix $\phi_0^+$ and $\phi_0^-$. Note that, with respect to the ordered basis $\mathscr{B}$ defined in \eqref{basis}, the vector $x = \sum_{i=1}^n (x_i e_i + y_i f_i)$ is represented by the $2n$-tuple
\begin{equation}
	\label{vecx}
x=(x_1,\dots,x_d,x_{d+1},y_{d+1},\dots,x_n,y_n,y_1,\dots,y_d).
\end{equation}
For each $\epsilon \in \{+,-\}$ we define a $2(n-d) \times 2(n-d)$ matrix
\begin{align*}
	K^+=\left( \begin{array}{ccccc}
			0 & 1 &        &   &   \\
			0 & 0 &        &   &   \\
			  &   & \ddots &   &   \\
			  &   &        & 0 & 1 \\
			  &   &        & 0 & 0
		\end{array} \right) \hspace{0.5cm} \text { or } \hspace{0.5cm}
	K^-=\left( \begin{array}{ccccccc}
			0 & 1 &        &   &   &   &   \\
			0 & 0 &        &   &   &   &   \\
			  &   & \ddots &   &   &   &   \\
			  &   &        & 0 & 1 &   &   \\
			  &   &        & 0 & 0 &   &   \\
			  &   &        &   &   & 1 & 1 \\
			  &   &        &   &   & 0 & 1
		\end{array} \right).
\end{align*}
Additionally, we define a $2n \times 2n$ matrix $K$ by
\begin{equation}
	\label{quadraticgram}
	K = \left( \begin{array}{ccc}
			0 & 0          & I_d \\
			0 & K^\varepsilon & 0   \\
			0 & 0          & 0
		\end{array} \right).
\end{equation}
\begin{lemma}
	\label{stabcorollary}
	Let $\varphi_0^\varepsilon\in\mathcal{Q}^\varepsilon$, $x\in V$, and $K$ be as in Equations $\eqref{eq-forms}$, $\eqref{vecx}$ and  $\eqref{quadraticgram}$, and let 
	\[
g=	\left(
		\begin{array}{ccc}
				P^{-T} & 0 & 0      \\
				0 & N & 0      \\
				0 & 0 & P
			\end{array}
		\right)
	\]
be an element of the Levi factor $L$ of $X_U$ defined in Corollary~$\ref{levi}$. Then
	\begin{enumerate}
	    \item[(a)] $\varphi_0^\varepsilon(x) = x K x^T$; 
	    \item[(b)] $g K g^T = K$ if and only if $N K^\varepsilon N^T = K^\varepsilon$; and
	    \item[(c)] if $N K^\varepsilon N^T = K^\varepsilon$ then $g$ is an isometry of $\phi_0^\varepsilon$.
	\end{enumerate}	
\end{lemma}

\begin{proof}
    (a) Write the tuple for $x$ in Equation \eqref{vecx} as $x=(a,b,c)$ where 
    \[
    \text{$a=(x_1,\dots,x_d)$, $b=(x_{d+1},y_{d+1},\dots,x_n,y_n)$ and $c=(y_1,\dots,y_d)$.}
    \]
    A straightforward computation shows that $xKx^T=bK^\varepsilon b^T + a c^T$, and a similar computation yields 
    \[
    bK^\varepsilon b^T = \left(\sum_{i=d+1}^nx_iy_i\right) + \delta(x_n^2+y_n^2)=\left(\sum_{i=d+1}^nx_iy_i\right) + \delta(x_n+y_n), 
    \]
    where $\delta=0$ if $\varepsilon=+$ and $\delta=1$ if $\varepsilon=-$. Since $a c^T =\sum_{i=1}^dx_iy_i$, we conclude that $\varphi_0^\varepsilon(x) = x K x^T$.

    (b) 
	Now
	\begin{align*}
		gKg^T = & \left(
		\begin{array}{ccc}
				P^{-T}& 0 & 0      \\
				0 & N & 0      \\
				0 & 0 & P
			\end{array}
		\right)
		\left(
		\begin{array}{ccc}
				0 & 0          & I_d \\
				0 & K^\varepsilon & 0   \\
				0 & 0          & 0
			\end{array}
		\right)
		\left(
		\begin{array}{ccc}
				P^{-1} & 0   & 0      \\
				0   & N^T & 0      \\
				0   & 0   & P^{T}
			\end{array}
		\right)          \\
		=       &
		\left(
		\begin{array}{ccc}
				0 & 0            & P^{-T}\\
				0 & N K^\varepsilon & 0 \\
				0 & 0            & 0
			\end{array}
		\right)
		\left(
		\begin{array}{ccc}
				P^{-1} & 0   & 0      \\
				0   & N^T & 0      \\
				0   & 0   & P^{T}
			\end{array}
		\right)          \\
		=       &
		\left(
		\begin{array}{ccc}
				0 & 0                & I_d \\
				0 & N K^\varepsilon N^T & 0   \\
				0 & 0                & 0
			\end{array}
		\right)
	\end{align*}
	and hence  $g K g^T = K$ if and only if $N K^\varepsilon N^T = K^\epsilon$.
	
    (c) Suppose now that $N K^\varepsilon N^T = K^\varepsilon$, so $g K g^T = K$ by part (b). 	Then by part (a), 
    \[
    \varphi_0^\varepsilon(xg) = (xg)K(xg)^T = x(gKg^T)x^T = xKx^T = \varphi_0^\varepsilon(x).
    \]
    Since this holds for all $x\in V$ it follows that $g$ is an isometry of $\phi_0^\varepsilon$.
\end{proof} 
Note that the converse of part (c) of Lemma \ref{stabcorollary} need not hold. A matrix $K$ with the property in Lemma \ref{stabcorollary}(a) is called a Gram matrix for $\varphi_0^\varepsilon$, and while each quadratic form in $\mathcal{Q}$ has a unique \emph{lower-triangular} Gram matrix, Gram matrices are not unique for quadratic forms in characteristic $2$. In particular not every element $h$ which stabilises $\varphi_0^\varepsilon$  in the Jordan-Steiner action is guaranteed to satisfy $hKh^T=K$.

\subsubsection{The action of \texorpdfstring{$X_U$}{XU} on quadratic forms}
\label{ticodesaresit}

In this subsection we show that the codes produced by Construction \ref{ticode} are $X$-strongly incidence-transitive. Recall that $U = \langle e_i \mid 1 \leq i \leq d \rangle$ and $\mathcal{Q}^\epsilon_d = \{ \phi\in \mathcal{Q}^\epsilon \mid U \subseteq \sing(\phi) \}$. We saw in Lemma \ref{intersectionthm} that $X_U$ is transitive on $\mathcal{Q}^\epsilon_d$, so let us set
\[
\text{$\Delta= \mathcal{Q}^\epsilon_d$\quad so that\quad $\overline{\Delta}= \mathcal{Q}^\epsilon \setminus\Delta = \mathcal{Q}^\epsilon_{d-1}$,}
\]
by Lemma~\ref{intersectionthm}(a). To achieve our goal, it is sufficient to show that,  for some $\phi\in \Delta$, $X_{U,\phi}$ acts transitively on $\overline{\Delta}$. As noted above, $\phi_0^\epsilon $ defined in \eqref{eq-forms} is in $\Delta$, so we may pick $\phi=\phi_0^\epsilon$. Further, Lemma \ref{formvectorpermiso} implies that, to do this it is sufficient to show that $X_{U,\phi}$ acts transitively on $\lambda_\phi^{-1}(\overline{\Delta}) = \{ c \in V \mid \phi_c \in \overline{\Delta} \}$, and by Lemma \ref{intersectionthm}(b),  $\lambda_\phi^{-1}(\overline{\Delta}) = \sing(\phi)\setminus (\sing(\phi)\cap U^\perp)$.
We prove transitivity of $X_{U,\phi}$ on $\sing(\phi)\setminus (\sing(\phi)\cap U^\perp)$ in two steps. First, in Corollary~\ref{singertransitive}(b), we show that $X_{U,\phi}$ acts transitively on the nontrivial $U^\perp$-cosets, and then, in Lemma~\ref{cosettrans}, we prove that, for a nontrivial $U^\perp$-coset $[w]$, the stabiliser $X_{U,\phi,[w]}$ acts transitively on the subset $[w]\cap\sing(\varphi)$. 

%
For a prime power $q$, an element $\alpha \in \mathbb{F}_q$ is a \emph{primitive element} if the multiplicative subgroup generated by $\alpha$ is the full multiplicative group $\mathbb{F}_q^\times$ of nonzero elements in $\mathbb{F}_{q}$. A monic polynomial $f \in \mathbb{F}_q[x]$ of degree $d \geq 1$ is a \emph{primitive polynomial} if it is the minimal polynomial of a primitive element of $\mathbb{F}_{q^d}$. For a primitive element $\alpha\in \mathbb{F}_{q^d}$, the set
$\{ \alpha^i \mid 0\leq i\leq d-1\}$ is a basis for $\mathbb{F}_{q^d}$ viewed as a vector space over $\mathbb{F}_{q}$.
An element of $\gl_d(q)$ of order $q^d-1$ is called a \emph{Singer cycle}. We describe a well known construction of a Singer cycle in $\gl_d(2)$ below.
\begin{constr}
\label{ssconstruction}
Let $\alpha$ be a primitive element in $\mathbb{F}_{2^d}$ and let $h(x) = x^d+\sum_{i=0}^{d-1} a_ix^i$ be the minimal polynomial of $\alpha$ over $\mathbb{F}_2$. Define  $f: \mathbb{F}_{2^d} \rightarrow  \mathbb{F}_2^d$ by
%
\[
	f: \sum_{i=0}^{d-1} c_i \alpha^i \mapsto  (c_0,...,c_{d-1}),
\]
%
so $f$ is an invertible linear transformation with $\mathbb{F}_{2^d}$ viewed as an $\mathbb{F}_2$-vector space.
The multiplicative group $\mathbb{F}_{2^d}^\times$ acts on $\mathbb{F}_{2^d}$ by multiplication modulo $h(x)$. We obtain a faithful linear representation $\rho: \mathbb{F}_{2^d}^\times \rightarrow \gl_d(2)$ by setting
\begin{equation}
	\label{rhorep}
	\rho(\alpha) =
	\left(
	\begin{array}{ccc}
			f(\alpha^1) \\
			f(\alpha^2) \\
			\vdots \\
			f(\alpha^{d})
		\end{array}
	\right)
	=
	\left(
	\begin{array}{cccc}
			0      &     &         &         \\
			\vdots &     & I_{d-1} &         \\
			0      &     &         &         \\
			a_0    & a_1 & \dots   & a_{d-1}
		\end{array}
	\right)
\end{equation}
where the $a_i$ are the coefficients of $h$ as above.
The multiplicative subgroup $\langle \rho(\alpha) \rangle < \gl_d(2)$ is called a \emph{Singer subgroup} and its generator $\rho(\alpha)$ is called a \emph{Singer cycle}.
\end{constr}

\begin{lemma}
	\label{singersubgroup}
	Let $\alpha, h, \rho$ be as in Construction~\ref{ssconstruction}, 
	let $s$ be the $2n \times 2n$ matrix defined by
	$$
	s = \left(
		\begin{array}{ccc}
				\rho(\alpha)^{-T} & 0          & 0       \\
				0            & I_{2(n-d)} & 0       \\
				0            & 0          & \rho(\alpha)
			\end{array}
		\right),
		$$
	and set $S = \langle s \rangle$. Then $S$ is a subgroup of the Levi factor $L$ of $X_U$
	, $S$ fixes $\phi_0^+$ and $\phi_0^-$ defined by Equation \eqref{eq-forms}, and $S$ is cyclic of order $2^d-1$.
\end{lemma}
\begin{proof}
    Let $S = \langle s \rangle$. Since $\alpha$ is a primitive element of $\mathbb{F}_{2^d}^\times \cong C_{2^d-1}$, it follows from Construction \ref{ssconstruction} that $S$ is cyclic of order $2^d -1$. Corollary \ref{levi} implies that $S \leq L< X_U$, and Lemma \ref{stabcorollary} implies that $S \leq X_{\varphi}$ for both of the forms $\phi$  in Equation \eqref{eq-forms}.
\end{proof}
\begin{lemma}
	\label{singerpermiso}
	Let $S = \langle s \rangle$ denote the subgroup of $X_U$ defined in Lemma \ref{singersubgroup}. The setwise action of $S$ on the quotient space $V/U^\perp$ is permutationally isomorphic to the action of $\mathbb{F}_{2^d}^\times$ on $\mathbb{F}_{2^d}$ by multiplication.
\end{lemma}
\begin{proof}
Let $\alpha, h, \rho$ be as in Construction~\ref{ssconstruction}, and choose bases $\{\alpha^i \mid 0 \leq i \leq d-1 \}$ for $\mathbb{F}_{2^d}$ and $\{ f_j + U^\perp \mid 1 \leq j \leq d \}$ for $V/U^\perp$. Define a map $\widetilde{f}: \mathbb{F}_{2^d} \rightarrow V/U^\perp$ by $\widetilde{f}(\alpha^i) = f_{i+1}+U^\perp$ for $0 \leq i \leq d-1$ and extending linearly to $\mathbb{F}_{2^d}$. Then $\widetilde{f}$ is a vector space isomorphism. Similarly, define $\widetilde{\rho}: \mathbb{F}_{2^d}^\times \rightarrow S$ by $\widetilde{\rho}(\alpha^i) = s^i$ for $0\leq i\leq 2^d-2$, and note that $\widetilde{\rho}$ is a group isomorphism. Now consider the natural action of $\mathbb{F}_{2^d}^\times$ on $\mathbb{F}_{2^d}$ by multiplication. For all $\sum_{i=0}^{d-1} c_i\alpha^i \in \mathbb{F}_{2^d}$ we have,  towards the end using the definition of $s$ in Lemma~\ref{singersubgroup} and the fact that $s$ leaves $U^\perp$ invariant,
	\begin{align}
	\label{singerpermisoverificiation}
	\widetilde{f} \left( \left( \sum_{i=0}^{d-1} c_i \alpha^i \right) \alpha \right) &= \widetilde{f} \left( \sum_{i=0}^{d-1} c_i \alpha^{i+1} \right) \nonumber \\
		& = \widetilde{f}\left( \sum_{i=0}^{d-2} c_i \alpha^{i+1} + c_{d-1} \sum_{i=0}^{d-1} a_i \alpha^i \right) \nonumber \\
		& = \sum_{i=0}^{d-2} c_i \widetilde{f}(\alpha^{i+1}) + c_{d-1} \sum_{i=0}^{d-1} a_i \widetilde{f}(\alpha^i)\nonumber \\
		& = \left( \sum_{i=0}^{d-2} c_i f_{i+2} + c_{d-1} \sum_{i=0}^{d-1} a_i f_{i+1} \right) + U^\perp  \nonumber \\
		& = \left( \left(\sum_{i=0}^{d-1} c_i f_{i+1}\right) \rho(\alpha) \right) + U^\perp  \nonumber \\
		& = \left( \sum_{i=0}^{d-1} c_i f_{i+1} \right)s + U^\perp   =  \left( \sum_{i=0}^{d-1} c_i f_{i+1}  + U^\perp \right)s                                     \nonumber \\
		& = \left(\widetilde{f}\left(\sum_{i=0}^{d-1} c_i \alpha^i \right) \right) s. 
	\end{align}

Equation \eqref{singerpermisoverificiation} implies that for all $\beta \in \mathbb{F}_{2^d}$ we have
\begin{equation}
\label{partialsingerpermiso}
\widetilde{f} \left( \beta \alpha \right) = \left(\widetilde{f} \left( \beta \right)\right) s = \left(\widetilde{f} \left( \beta \right)\right) \widetilde{\rho}(\alpha). 
\end{equation}
Now, if $\gamma \in \mathbb{F}_{2^d}^\times$ then there exists $j$ such that $\gamma=\alpha^{j}$ with $0 \leq j \leq 2^d-2$. Successively applying Equation \eqref{partialsingerpermiso} we find that, for all $\beta \in \mathbb{F}_{2^d}$, 
\[
\widetilde{f}(\beta \gamma) = \widetilde{f}(\beta \alpha^j) = \widetilde{f}(\beta \alpha^{j-1})\widetilde{\rho}(\alpha) = \dots = \widetilde{f}(\beta)\widetilde{\rho}(\alpha)^j= \widetilde{f}(\beta)\widetilde{\rho}(\alpha^j) = \widetilde{f}(\beta)\widetilde{\rho}(\gamma).
\]
Therefore the action of $S$ on $V/U^\perp$ is permutationally isomorphic to the action of $\mathbb{F}_{2^d}^\times$ on $\mathbb{F}_{2^d}$.
\end{proof}

For $w \in V$, we use $[w]:= w+U^\perp$ to denote the image of $w$ under the natural projection map $\pi: V \rightarrow V/U^\perp$; the coset $[w]$ is called nontrivial if $[w]\ne U^\perp$.

\begin{corollary}\leavevmode\vspace{-\baselineskip}	\label{singertransitive}
	\begin{enumerate}
	    \item[(a)] The Singer subgroup $S = \langle s \rangle$ acts faithfully and regularly on the nontrivial cosets of $U^\perp$ in $V$.
	    \item[(b)] For the form $\varphi=\varphi_0^\epsilon$ defined in Equation~\eqref{eq-forms}, the stabiliser $X_{U,\varphi}$ is transitive on the nontrivial cosets of $U^\perp$ in $V$.
	\end{enumerate}
\end{corollary}

\begin{proof}
	The group $\mathbb{F}_{2^d}^\times$ acts regularly on itself by multiplication so Lemma \ref{singerpermiso} implies that $S$ acts faithfully and regularly on the nontrivial cosets of $V/U^\perp$. By Lemma~\ref{singersubgroup}, $S\leq X_{U,\varphi}$, and hence $X_{U,\varphi}$ is transitive on the nontrivial cosets of $U^\perp$ in $V$.
\end{proof}

Now we complete the proof that $X_{U,\varphi}$ is transitive on $\lambda_\phi^{-1}(\overline{\Delta}) = \sing(\phi)\setminus (\sing(\phi)\cap U^\perp)$.
\begin{lemma}
	\label{cosettrans}
	Let $\varphi=\varphi_0^\epsilon$ be as in Equation~\eqref{eq-forms}. Then  $\sing(\phi)\not\subseteq U^\perp$, and the stabiliser $X_{U,\varphi}$ is transitive on $\sing(\phi)\setminus (\sing(\phi)\cap U^\perp)$.
\end{lemma}

\begin{proof}
    First we note that  $\sing(\varphi)\not\subseteq U^\perp$ since, by Equation~\eqref{eq-forms}, $\varphi(f_1)=0$ and $f_1\notin U^\perp$. By  Corollary~\ref{singertransitive}(b), $X_{U,\phi}$ acts transitively on the nontrivial $U^\perp$-cosets in $V$, so to prove that $X_{U,\varphi}$ is transitive on $\sing(\phi)\setminus (\sing(\phi)\cap U^\perp)$ it is sufficient to prove that $X_{U,\varphi,[w]}$ acts transitively on $[w]\cap \sing(\varphi)$, for some nontrivial coset 
    $[w]=w+U^\perp$. We may assume without loss of generality that $[w]=f_1+U^\perp$ with $w=f_1\in [w]\cap\sing(\varphi)$. Recalling that $U=\langle e_1,\dots, e_d\rangle$, it follows from the definition of $B$ in \eqref{Bdef} that each $u\in U^\perp$ has the form $u = \sum_{i=1}^{n} u_i e_i + \sum_{i=d+1}^n v_if_i$, and so each element of $[w]$ is $u+f_1$ for some $u$ of this form. Suppose that such an element $u+f_1\in [w]\cap\sing(\varphi)$, and consider the matrix 
	\[
	g=\left(
		\begin{array}{ccc}
				I_d & 0 & 0 \\
				JZ^T & I_{2(n-d)} & 0 \\
				Q & Z & I_d
			\end{array}
		\right)
	\]
	where $Q,Z$ are defined by

	\begin{center}
		$ Q = \left(
			\begin{array}{cccc}
					u_1    & u_2    & \dots  & u_d    \\
					u_2    & 0      & \dots  & 0      \\
					\vdots & \vdots & \ddots & \vdots \\
					u_d    & 0      & \dots  & 0
				\end{array}
			\right)$,
		$Z = \left(
			\begin{array}{ccccc}
					u_{d+1} & v_{d+1} & \cdots & u_n    & v_n    \\
					0       & 0       & \cdots & 0      & 0      \\
					\vdots  & \vdots  & \cdots & \vdots & \vdots \\
					0       & 0       & \cdots & 0      & 0
				\end{array}
			\right)$,
	\end{center}
and $J$ is the $2(n-d) \times 2(n-d)$ block-diagonal matrix with blocks $\begin{pmatrix}
			0 & 1 \\ 1 & 0
		\end{pmatrix}$.

	By Corollary \ref{levi}, $g \in R$ if and only if $Q+Q^T = Z J Z^T$. Since the matrix $Q$ is symmetric, $Q + Q^T$ is the $d \times d$ zero matrix $\textbf{0}_{d \times d}$ (recall that we are working over $\mathbb{F}_2$). Computing $Z J Z^T$, we have
	\begin{align*}
		Z J Z^T = & \left(
		\begin{array}{ccccc}
				u_{d+1} & v_{d+1} & \cdots & u_n    & v_n    \\
				0       & 0       & \cdots & 0      & 0      \\
				\vdots  & \vdots  & \cdots & \vdots & \vdots \\
				0       & 0       & \cdots & 0      & 0
			\end{array}
		\right)
		\left( \begin{array}{ccccc}
				0 & 1 &        &   &   \\
				1 & 0 &        &   &   \\
				  &   & \ddots &   &   \\
				  &   &        & 0 & 1 \\
				  &   &        & 1 & 0
			\end{array} \right)
		\left(
		\begin{array}{cccc}
				u_{d+1} & 0      & \cdots & 0      \\
				v_{d+1} & 0      & \cdots & 0      \\
				\vdots  & \vdots & \cdots & \vdots \\
				u_n     & 0      & \cdots & 0      \\
				v_n     & 0      & \cdots & 0
			\end{array}
		\right)                                        \\
		=                    &
		\left(
		\begin{array}{ccccc}
				v_{d+1} & u_{d+1} & \cdots & v_n    & u_n    \\
				0       & 0       & \cdots & 0      & 0      \\
				\vdots  & \vdots  & \cdots & \vdots & \vdots \\
				0       & 0       & \cdots & 0      & 0
			\end{array}
		\right)
		\left(
		\begin{array}{cccc}
				u_{d+1} & 0      & \cdots & 0      \\
				v_{d+1} & 0      & \cdots & 0      \\
				\vdots  & \vdots & \cdots & \vdots \\
				u_n     & 0      & \cdots & 0      \\
				v_n     & 0      & \cdots & 0
			\end{array}
		\right) =  \textbf{0}_{d \times d}.
	\end{align*}

Therefore $g \in R$, and so by Corollary~\ref{levi}, $g\in X_U$ and hence also $g$ stabilises $U^\perp$. Moreover,  by the definition of the ordered basis $\mathscr{B}$ in \eqref{basis} and of the matrix $g$, we see that $g$ maps $w=f_1$ to $u + w\in [w]$. Thus $g \in X_{U,[w]}$. It remains to confirm that $g$ fixes $\varphi=\phi_0^\epsilon$ in each case $\varepsilon = \pm$.
	With respect to the basis $\mathscr{B}$ in Equation \eqref{basis}, an arbitrary element of $V$ has the form  $x = \sum_{i=1}^n \left( x_i e_i + y_i f_i  \right)$, which as a row vector is $(a,b,c)$, where
	\[
	a=(x_1, \dots, x_d), \quad b=(x_{d+1}, y_{d+1},\dots, x_n,y_n),\quad\text{and}\quad c=(y_1,\dots,y_d). 
	\]
	The row vector representing the image $xg$ is therefore $(a+bJZ^T+cQ, b+cZ, c)$, and we note the following, using the definitions of $Q, J, Z$:
	\begin{align*}
	    bJZ^T& = (y_{d+1}, x_{d+1},\dots, y_n, x_n)Z^T = \left(\sum_{i=d+1}^n(y_iu_i+x_iv_i), 0, \dots, 0\right),\\
	    cQ&=\left(\sum_{i=1}^d y_iu_i, y_1u_2,\dots, y_1u_d\right)\quad \text{and}\quad 
	    cZ=(y_1u_{d+1}, y_1v_{d+1},\dots, y_1u_n, y_1v_n).
	\end{align*}
Thus
	\begin{align*}
	xg = \left( x_1 + \sum_{i=1}^n u_i y_i + \sum_{i=d+1}^n v_i x_i \right) e_1 
		      + \sum_{i=1}^d y_{i}f_i + \sum_{i=d+1}^{n} (y_{i}+y_1 v_{i})f_{i}  
		      + \sum_{i=2}^n (x_i + y_1 u_i)e_i.
	\end{align*}
	 Recall that $y^2 = y$ for each $y\in \mathbb{F}_2$. Using \eqref{eq-forms}, if $\epsilon = +$ then
	\begin{align*}
		\varphi_0^+(xg) = & x_1 y_1 + y_1 \left(\sum_{i=1}^n u_i y_i + \sum_{i=d+1}^n v_i x_i \right) + \sum_{i=2}^d \left(x_i y_i+ y_1 u_i y_i \right) \\
		                & + \left(\sum_{i=d+1}^n x_i y_i \right)+ y_1 \sum_{i=d+1}^n \left( v_i x_i + u_i y_i + y_1 u_i v_i \right)                                \\
		=               & \varphi_0^+(x) + y_1 \left( u_1 + \sum_{i=d+1}^n u_i v_i \right)                                                               \\
		=               & \varphi_0^+(x) + y_1 \varphi_0^+(u+f_1).
	\end{align*}
Now using \eqref{linktwoforms}, for $\epsilon = -$, we get
	\begin{align*}
		\varphi_0^-(xg) = & \varphi_0^+(xg)+B(xg,e_n+f_n)\\ 
		=&\varphi_0^+(x) + y_1 \varphi_0^+(u+f_1)+(x_n+y_1u_n)+(y_n+y_1v_n)\\
			=&\varphi_0^-(x) + y_1 \varphi_0^+(u+f_1)+y_1(u_n+v_n)\\
			=&\varphi_0^-(x) + y_1 \varphi_0^+(u+f_1)+y_1B(u+f_1,e_n+f_n)\\
				=&\varphi_0^-(x) + y_1 \varphi_0^-(u+f_1).
	\end{align*}
Thus $\varphi(xg)=\varphi(x) + y_1 \varphi(u+f_1)$, for both values of $\varepsilon$. Now recall our assumptions that $u\in U^\perp$ and $u+f_1\in [w]\cap\sing(\varphi)$. Thus $\varphi(u+f_1) = 0$, and hence $\varphi(xg) = \varphi(x)$. Since this holds for all $x \in V$, it follows that $\varphi^{g^{-1}} = \varphi$, so $g\in X_{U,\varphi, [w]}$. This means that, for each $u+f_1\in [w]\cap\sing(\varphi)$ we have constructed an element of $X_{U,\varphi, [w]}$ that maps $f_1$ to $u+f_1$. Therefore the stabiliser $X_{U, \varphi,[w]}$ acts transitively on $[w]\cap \sing(\varphi)$.
\end{proof}
%
\begin{corollary}\label{cor:transQd-1}
    For $\phi=\phi_0^\epsilon$, the group $X_{U,\phi}$ acts transitively on $\mathcal{Q}_{d-1}^\epsilon$.
\end{corollary}
\begin{proof}
 By Lemma~\ref{formvectorpermiso}, the statement is equivalent to proving that $X_{U,\varphi}$ is transitive on $\lambda_\phi^{-1}(\mathcal{Q}_{d-1}^\epsilon) := \{ c \in V \mid \phi_c \in \mathcal{Q}_{d-1}^\epsilon \}$,  and by Lemma \ref{intersectionthm}(b),  $\lambda_\phi^{-1}(\mathcal{Q}_{d-1}^\epsilon) = \sing(\phi)\setminus (\sing(\phi)\cap U^\perp)$. Finally, by Lemma~\ref{cosettrans}, the action of $X_{U,\varphi}$ on this set is transitive.
\end{proof}
We now show that the only $X$-strongly incidence-transitive codes for which a codeword stabiliser leaves invariant a totally isotropic   $d$-dimensional subspace  of $V$, where  $(d,\varepsilon)\ne (n,-)$, are those in Construction \ref{ticode}.

\begin{theorem}\samepage
	\label{ticodeissit}
Let $n\geq 2$, and $X=\syp_{2n}(2)$. Let $\Gamma$ be an $X$-strongly incidence-transitive code in $J(\mathcal{Q}^\varepsilon,k)$  for some $k$ satisfying $2 \leq k \leq |\mathcal{Q}^\varepsilon|-2$, and let $\Delta\in\Gamma$. Suppose that $X_\Delta\leq X_U$ for some totally isotropic   $d$-dimensional subspace $U$ of $V$, where $(d,\varepsilon)\ne (n,-)$. Then $\Gamma = \Gamma(n,d,\varepsilon,\delta)$, for some $\delta=0, 1$, and \[k=
\begin{cases}
2^{n-1}(2^{n-d}+\varepsilon) &\text{ if } \delta=0 \\
2^{2n-d-1}(2^{d}-1) & \text{ if } \delta=1,
\end{cases} \]    as in Construction~\ref{ticode}. 
\end{theorem}

%
\begin{proof}
By Lemma~\ref{intersectionthm}(a), $\mathcal{Q}^\varepsilon = \mathcal{Q}^\varepsilon_d\cup \mathcal{Q}^\varepsilon_{d-1}$, and since $(d,\varepsilon)\ne (n,-)$, each of the sets $\mathcal{Q}^\varepsilon_d$ and 
$\mathcal{Q}^\varepsilon_{d-1}$ is non-empty and $X_U$-invariant. Since $X_\Delta\leq X_U$, each of these sets is also $X_\Delta$-invariant, and since $\Gamma$ is $X$-strongly incidence-transitive, $X_\Delta$ has just two orbits in $\mathcal{Q}^\varepsilon$, namely $\Delta$ and $\overline{\Delta}$. It follows that  
\[
\Delta= \Delta^\delta(U) = \mathcal{Q}^\varepsilon_{d-\delta} = \{ \phi \in \mathcal{Q}^\epsilon \mid \dim(\sing(\phi)\cap U) = d-\delta \},
\]
where $\delta=0$ or $1$. Thus $\Gamma$ is the code $\Gamma(n,d,\varepsilon,\delta)$ in Construction~\ref{ticode}.

Finally we determine $k=|\Delta|$, where $\Delta = \Delta^\delta(U)$. Let $\phi\in\Delta$ and suppose first that $\delta=0$. By Lemma \ref{intersectionthm} and Corollary \ref{formvectorbijection}, $|\Delta|=|\sing(\phi)\cap U^\perp|$. Without loss of generality, we may take $U=\langle e_1,\dots, e_d\rangle$ and $\phi=\phi_0^\epsilon$ as in \eqref{eq-forms} (since $(d,\varepsilon)\ne (n,-)$).
Then $U^\perp=\langle e_1,\dots, e_d,e_{d+1},f_{d+1},\ldots,e_n,f_n \rangle$, so each $u\in U^\perp$ has the form $u = \sum_{i=1}^{n} u_i e_i + \sum_{i=d+1}^n v_if_i$. Then $u\in \sing(\phi)\cap U^\perp$ if and only if $\phi(u)=0$ and by \eqref{eq-forms}, 
\[
	\phi(u)=	\left\{
		\begin{array}{ll}
		\phi_0^+(u) =	\sum_{i=d+1}^n u_i v_i   \\
		\phi_0^-(u) =	u_n + v_n + \sum_{i=d+1}^n u_i v_i.
		\end{array}
			\right.
\]
For a solution to $\phi(u)=0$, the coefficients  $u_1,\dots,u_d$ may be chosen arbitrarily and $u-\sum_{i=1}^du_ie_i$ must be a singular vector for the nondegenerate quadratic form of type $\epsilon$  on $\F_2^{2(n-d)}$ as in \eqref{eq-forms}. 
Thus 
\[
k=|\sing(\phi)\cap U^\perp|=2^d\cdot 2^{n-d-1}(2^{n-d}+\epsilon)=2^{n-1}(2^{n-d}+\epsilon).
\]
Finally, if $\delta=1$ then, since $\Gamma(n,d,\epsilon,0)$ and $\Gamma(n,d,\epsilon,1)$ are complementary codes,
\begin{align*}
    k &= 2^{n-1}(2^n+\epsilon) - 2^{n-1}(2^{n-d}+\epsilon) \\
    &= 2^{n-1}(2^n-2^{n-d}) \\
    &= 2^{2n-d-1}(2^d-1).
\end{align*}
This completes the proof.
\end{proof}

We  now  analyse the family of codes in Construction~\ref{ticode}.

\begin{theorem}\label{p:ti}
Let $n\geq 2$, $1\leq d\leq  n$, $(d,\epsilon)\neq (n,-)$, $X=\syp_{2n}(2)$, and let $\Gamma = \Gamma(n,d,\varepsilon,\delta)$ in  $J(\mathcal{Q}^\varepsilon,k)$ as in Construction~\ref{ndcode},  for   some $\delta\in\{0,1\}$, with $k$ as in \eqref{eqticode}. Let $A=\aut(\Gamma)$.
 Then: 
\begin{enumerate}[(a)]
  
    \item For $\Delta^\delta(U)\in\Gamma$, with $U$ a  totally-isotropic  $d$-dimensional subspace of $V$, the stabiliser $X_U$ is equal to $X_{\Delta^\delta(U)}$.
    \item The code $\Gamma$  is $X$-strongly incidence-transitive and is not self-complementary. 
    \item The integer $k$ satisfies $3 \leq k \leq |\mathcal{Q}^\varepsilon|-3$ unless $(n,d,\epsilon)= (2,1,-)$; and 
    \item $A=X$.
\end{enumerate}
\end{theorem}

\begin{proof} Let $\Delta$ be a codeword of $\Gamma$, so  
$\Delta=\Delta^\delta(U)$ for some  totally-isotropic  $d$-dimensional subspace $U$ of $V$ and $\delta\in\{0,1\}$.

(a) By Lemma~\ref{intersectionthm}(a) and Construction~\ref{ticode},  $\calq^\epsilon=\calq_d^\epsilon\cup\calq_{d-1}^\epsilon=\Delta^0(U)\cup \Delta^1(U)$. Thus $X_U$ leaves invariant both $\Delta^0(U)$ and $\Delta^1(U)$, and so $X_U\leq X_\Delta$. 
 Using the information in \cite[Main Theorem and Table 3.5.C]{kl} if $n\geq 7$, and from \cite[Tables 8.28, 8.48, 8.64, and 8.80]{colva} for $2\leq n\leq 6$, we see that $X_U$ is a maximal subgroup of $X$, and hence $X_U= X_\Delta$,

(b) It is straightforward to check that $k\neq |\calq^\epsilon|/2$, so $\Gamma$ is not self-complementary. 
Since $\Gamma(n,d,\epsilon,0)$ and $\Gamma(n,d,\epsilon,1)$ are complementary codes (that is, the codewords of one code are the complements in $\mathcal{Q}^\varepsilon$ of the codewords of the other code),  
 $\Gamma(n,d,\epsilon,0)$ is $X$-strongly incidence-transitive if and only if $\Gamma(n,d,\epsilon,1)$ is $X$-strongly incidence-transitive. Therefore it is sufficient to prove 
that $\Gamma(n,d,\epsilon,0)$ is $X$-strongly incidence-transitive, and from now on we assume that $\Gamma = \Gamma(n,d,\epsilon,0)$. 

By Corollary \ref{formsetpermiso}, for each $g \in X$ we have
\begin{align*}
\Delta^0(U)g &= \{ \phi^g \in \mathcal{Q}^\epsilon \mid \dim(\sing(\phi)\cap U) = d \} =  \{\phi^g \in \mathcal{Q}^\epsilon \mid U \subseteq \sing(\phi)\} \\
&= \{\phi \in \mathcal{Q}^\epsilon \mid U \subseteq \sing(\phi^{g^{-1}})\} 
= \{\phi \in \mathcal{Q}^\epsilon \mid U \subseteq \sing(\phi)g^{-1}\} \\
&= \{\phi \in \mathcal{Q}^\epsilon \mid Ug \subseteq \sing(\phi)\} 
= \Delta(Ug).
\end{align*}
Therefore the action of $X$ on the codewords of $\Gamma$ is equivalent to the action of $X$ on the $d$-dimensional totally-isotropic subspaces of $V$. Witt's Lemma implies that the latter of these actions is transitive, and therefore the former is also transitive.
By Definition~\ref{defsit} therefore, our task is to prove that, for each  codeword $\Delta=\Delta^0(U)$,  the stabiliser $X_{\Delta}=X_U$ is transitive on $\Delta\times \overline{\Delta}$, where $\overline{\Delta}=\mathcal{Q}^\varepsilon\setminus {\Delta}=\mathcal{Q}_{d-1}^\varepsilon$.
By the transitivity of $X$ on $\Gamma$ it is sufficient to do this for a single codeword, and so we shall do this in the case where 
$U = \langle e_i \mid 1 \leq i \leq d \rangle$. By Lemma \ref{intersectionthm}(c),  $X_\Delta=X_U$ is transitive on $\Delta= \mathcal{Q}^\epsilon_d$, and the form $\phi=\varphi_0^\epsilon$ defined in Equation \eqref{eq-forms} lies in $\Delta$. 
By Corollary \ref{cor:transQd-1}, we see that $X_{U,\phi}$ is transitive on $\overline{\Delta}$,
 and we conclude that $\Gamma$ is  $X$-strongly incidence-transitive.
 
 (d) It is straightforward to check that  $3\leq k\leq |\mathcal{Q}^\varepsilon|-3$ unless $(n,d,\varepsilon)=(2,1,-)$. In the exceptional case we have $|\mathcal{Q}^-|=6$ and $k\in\{2,4\}$.
 
(e) From the definition of $\Gamma$ it is clear that $X\leq A=\Aut(\Gamma)$, and since in all cases $k\ne |\mathcal{Q}^\epsilon|/2$ by part (b), the automorphism group $A$ is a subgroup of $\sym(\mathcal{Q}^\epsilon)$.

Assume first that $(n,\epsilon)\neq (2,-)$. 
Since $\Gamma\neq \binom{\mathcal{Q}^\epsilon}{k}$ by Lemma \ref{le:noncomplete},  $A$ is a proper subgroup of $\sym(\mathcal{Q}^\epsilon)$ and does not contain $\mathrm{Alt}(\mathcal{Q}^\epsilon)$. It follows from the classification of the maximal subgroups of finite symmetric and alternating groups that $A=X$, see \cite[Chapter 9]{maximalfactorisations}, especially Tables II--VI.

Now assume $(n,\epsilon)= (2,-)$.
Then $|\mathcal{Q}^-|=6$, and $(d,\epsilon)\neq (n,-)$ implies that $d=1$. Thus $k=2$ if $\Gamma=\Gamma(2,1,-,0) $ and  $k=4=|\mathcal{Q}^-n|-2$ if $\Gamma=\Gamma(2,1,-,1)$. 
Now $X=\syp_4(2)=\sym(\mathcal{Q}^-)$ is $3$-transitive on $\mathcal{Q}^-$. It follows that $\Gamma= \binom{\mathcal{Q}^-}{2}$, and so  $A=\aut(J(\mathcal{Q}^-,2))= \sym(\mathcal{Q}^-)=X$.
\end{proof}

By Theorem~\ref{p:ti}, $X$ acts transitively on the code and the stabiliser of a codeword $\Delta^\delta(U)$ is equal to $X_U$.  Thus the number of codewords is equal to $|X:X_{U}|$. 
Since the size of the unipotent radical $R$ of $X_U$ is $2^{d/2-3d^2/2+2dn}$ (see \cite[Proposition 4.1.19(II)]{kl}), we have:
\begin{align*}
|\Gamma|&=\frac{|X|}{|X_U|}=\frac{|X|}{|R|\cdot |L|}=\frac{|\syp_{2n}(2)|}{ |R|\cdot|\gl_{d}(2)|\cdot |\syp_{2(n-d)}(2)|}\\
&=\frac{2^{n^2}\prod_{i=1}^n (2^{2i}-1)}{2^{d/2-3d^2/2+2dn}\cdot 2^{d(d-1)/2}\prod_{i=1}^d (2^{i}-1)\cdot 2^{(n-d)^2}\prod_{i=1}^{n-d} (2^{2i}-1)}\\
&=\frac{\prod_{i=n-d+1}^n (2^{2i}-1)}{\prod_{i=1}^d (2^{i}-1)}
=\frac{\prod_{j=1}^d (2^{2(n-d+j)}-1)}{\prod_{i=1}^d (2^{i}-1)}\\
&=\prod_{i=1}^d \frac{2^{2(n-d+i)}-1}{2^{i}-1}.
\end{align*}

\begin{theorem}\label{p:distanceti}
The minimum distance $D$ of $\Gamma(n,d,\varepsilon,\delta)$ satisfies
\[
D = \left\{ \begin{array}{ll}
     2^{2n-d-2}     &  \mbox{if  $\epsilon=+$ and $d<n$}\\
     2^{n-1}     &  \mbox{if  $\epsilon=+$ and $d=n$}\\
     2^{n-2}(2^{n-d}-1)     &  \mbox{if   $\epsilon=-$.}\\
\end{array}   \right. 
\]
\end{theorem}

\begin{proof}
As $\Gamma(n,d,\varepsilon,0)$ and $\Gamma(n,d,\varepsilon,1)$ are complementary, they have the same minimum distance. We will compute the minimum distance of  $\Gamma=\Gamma(n,d,\varepsilon,0)$. 

Let $U_1,U_2$ be two totally-isotropic subspaces of dimension $d$. 
Recall that in the Johnson graph $J(\mathcal{Q}^\varepsilon,k)$, the distance between two codewords $\Delta^0(U_1)$ and  $\Delta^0(U_2)$ is $k-|\Delta^0(U_1)\cap \Delta^0(U_2)|$, so we are looking for the largest possible $|\Delta^0(U_1)\cap \Delta^0(U_2)|$.
Let $\ell=\dim(U_1\cap U_2)$, so $\ell\geq0$, and since $U_1, U_2$ are distinct, also $\ell<d$. Let $\ell+m=\dim(U_1^\perp\cap U_2)$. Since $U_1\cap U_2\leq U_1^\perp\cap U_2\leq U_2$, the integer $m\geq 0$ and $\ell+m\leq d$, and hence $0\leq \ell\leq d-m$. Now $(U_1^\perp\cap U_2)/(U_1\cap U_2)\cong ((U_1^\perp\cap U_2) +U_1)/U_1$ is an $m$-dimensional totally isotropic subspace of the symplectic space $U_1^\perp/U_1$ of dimension $2(n-d)$, and the largest totally-isotropic subspace has dimension $n-d$, so $m\leq n-d$.
By Witt's lemma, for a fixed choice of $m$ and $\ell$, all the pairs $(U_1,U_2)$ are in the same orbit of $X$,
so we may assume that 
\begin{align*}
 U_1&=\langle e_1,e_2,\ldots, e_d\rangle\\ 
 U_2&=\langle e_1,e_2,\ldots, e_\ell,e_{d+1},\ldots, e_{d+m},f_{\ell+1},\ldots,f_{d-m}\rangle
\end{align*}
Assume first that $d+m<n$ or $\epsilon=+$. It follows that $U_1\cup U_2\subseteq \sing(\phi)$ for $\phi$ as in \eqref{eq-forms}.
Then $\Delta^0(U_i)=\mathcal{Q}^\epsilon_d=\{\phi_c|c\in \sing(\phi)\cap U_i^\perp\}$ by Lemma \ref{intersectionthm}, and so  $\Delta^0(U_1)\cap \Delta^0(U_2)=\{\phi_c|c\in \sing(\phi)\cap U_1^\perp\cap U_2^\perp\}$.
We have 
\begin{align*}
 U_1^\perp\cap U_2^\perp&=(U_1+U_2)^\perp\\
 &=\langle e_1,e_2,\ldots,  e_{d+m},f_{\ell+1},\ldots,f_{d-m}\rangle^\perp\\
 &=\langle e_1,\ldots, e_\ell,e_{d-m+1},\ldots, e_{n},f_{d+m+1},\ldots,f_{n}\rangle
\end{align*}
For $c=\sum_{i=1}^\ell u_i e_i+ \sum_{i=d-m+1}^n u_i e_i+\sum_{i=d+m+1}^n v_i f_i\in U_1^\perp\cap U_2^\perp$, we have 
\[
\phi(c)=\begin{cases}
\sum_{i=d+m+1}^n u_iv_i&\text{ if }\epsilon=+\\ 
\sum_{i=d+m+1}^n u_iv_i+u_n+v_n&\text{ if }\epsilon=- \\
\end{cases}
\]

Hence there is no condition on $\sum_{i=1}^\ell u_i e_i$, or on $\sum_{i=d-m+1}^{d+m} u_i e_i$, while  $\sum_{i=d+m+1}^n (u_ie_i+v_i f_i)$ must be singular for a form of type $\epsilon$ on $\mathbb{F}_2^{2(n-d-m)}$.
Therefore \[|\sing(\phi)\cap U_1^\perp\cap U_2^\perp|=
2^\ell\cdot 2^{2m}\cdot 2^{n-d-m-1}(2^{n-d-m}+\epsilon)=2^{\ell+n-d+m-1}(2^{n-d-m}+\epsilon).\]

Hence 
\[
|\Delta^0(U_1)\cap\Delta^0(U_2)|=
2^{\ell+n-d+m-1}(2^{n-d-m}+\epsilon)
\] 
in this case.

Now assume that $d+m=n$ and $\epsilon=-$ (which implies that $0<m\leq d$  and $0<n-d\leq d$) and take again $\phi$ as in \eqref{eq-forms}. 
Then $W=U_2\cap \sin(\phi)$ is a subspace of dimension $d-1$, and let $t\in U_2\setminus W$ so that $U_2=W+\langle t\rangle$ and $\phi(t)=1$.
Note that $t\in W^\perp$ since $U_2\leq U_2^\perp$. 
\begin{align*}
    \Delta^0(U_2)&=\{\phi_c|U_2\subseteq \sing(\phi_c) \text{ and }c\in \sing(\phi)\}\\
   &=\{\phi_c|U_2\subseteq \sing(\phi)+c \text{ and }c\in \sing(\phi)\}\\
   &=\{\phi_c|c\in W^\perp\setminus t^\perp  \text{ and }c\in \sing(\phi)\}.
\end{align*}

Here we have $W=\langle e_1,e_2,\ldots, e_\ell,e_{d+1},\ldots, e_{n-1},f_{\ell+1},\ldots,f_{d-m}\rangle$ and $t=e_n$.
Thus  $\Delta^0(U_1)\cap \Delta^0(U_2)=\{\phi_c|c\in \sing(\phi)\cap (U_1^\perp\cap W^\perp)\setminus t^\perp\}$. Now 
$U_1^\perp\cap W^\perp = (U_1+ W)^\perp$, and 
\[
(U_1+ W)^\perp = \langle e_1,e_2,\ldots,  e_{n-1},f_{\ell+1},\ldots,f_{d-m}\rangle^\perp
=\langle e_1,\ldots, e_\ell,e_{d-m+1},\ldots, e_{n},f_{n}\rangle,
\]
and hence 
\[
(U_1^\perp\cap W^\perp)\setminus t^\perp = (U_1+W)^\perp\setminus e_n^\perp = f_n+\langle e_1,\ldots, e_\ell,e_{d-m+1},\ldots, e_{n}\rangle.
\]
%
 For $c=f_n+\sum_{i=1}^\ell u_i e_i+ \sum_{i=d-m+1}^n u_i e_i\in (U_1^\perp\cap W^\perp)\setminus t^\perp$ we have, by \eqref{eq-forms}, $\phi(c)=1+u_n+u_n=1$, so $c\notin \sing(\phi)$. Thus  $\Delta^0(U_1)\cap \Delta^0(U_2)$ is empty.
 
We now want to maximize the value of $|\Delta^0(U_1)\cap\Delta^0(U_2)|$ over all possible $m\leq \min(d,n-d)$ and $0\leq \ell\leq d-m$ (or  $0\leq \ell< d$ if $m=0$). 

For $m=0$, clearly the largest value is when $\ell$ is  as large as  possible, that is, for $\ell=d-1$. 
This gives us 
\begin{equation}\label{e:m0}
    |\Delta^0(U_1)\cap\Delta^0(U_2)|=2^{n-2}(2^{n-d}+\epsilon).
\end{equation}

For fixed $m>0$ such that $m<n-d$ if $\epsilon=-$, clearly the largest value is when $\ell$ is as large as possible, that is, for $\ell=d-m$. 
This gives us $|\Delta^0(U_1)\cap\Delta^0(U_2)|=
2^{n-1}(2^{n-d-m}+\epsilon)$. This in turn is maximised when $m$ is as small as possible, that is for $m=1$, where we have 
\begin{equation}\label{e:m1}
    |\Delta^0(U_1)\cap\Delta^0(U_2)|=2^{n-1}(2^{n-d-1}+\epsilon).
\end{equation}
Note this case does not happen if $n-d=1$ and $\epsilon=-$.

If $n-d\leq d$ and $\epsilon=-$, we can also have $m=n-d$. However in this case $|\Delta^0(U_1)\cap\Delta^0(U_2)|=0$ for any $\ell$, which is clearly never going to be the largest value for $|\Delta^0(U_1)\cap\Delta^0(U_2)|$.  Thus we need to compare the quantities in \eqref{e:m0} and \eqref{e:m1} to determine which is larger.

Case 1:  $\epsilon=+$ and $n>d$. Here the larger is $2^{n-1}(2^{n-d-1}+1)$ (with $m=1$ and $\ell=d-1$), since $2^{n-2}(2^{n-d}+1)< 2^{n-1}(2^{n-d-1}+1)$.
The minimum distance is \[k-2^{n-1}(2^{n-d-1}+1)=2^{n-1}(2^{n-d}+1)-2^{n-1}(2^{n-d-1}+1)=2^{2n-d-2}.\]

Case 2: $\epsilon=+$ and $n=d$. Here there is only one case:  $m=0$, so the largest $|\Delta^0(U_1)\cap\Delta^0(U_2)|$ is $2^{n-2}(2^{n-d}+\epsilon)=2^{n-1}$, with $m=0$ and $\ell=d-1$.
The minimum distance is \[k-2^{n-1}=2^{n-1}(2^{n-d}+1)-2^{n-1}=2^{n-1}.\]

Case 3: $\epsilon=-$ and  $n-d>1$. Here the larger is $2^{n-2}(2^{n-d}-1)$ (with $m=0$ and $\ell=d-1$), since $ 2^{n-1}(2^{n-d-1}-1)< 2^{n-2}(2^{n-d}-1)$.
The minimum distance is again $2^{n-2}(2^{n-d}-1).$

Case 4: $\epsilon=-$ and  $n-d=1$. 
Here, as observed above, the case \eqref{e:m1} does not occur so the largest value is $2^{n-2}(2^{n-d}+\epsilon)$ (with $\epsilon=-$, $m=0$, and $\ell=d-1$).
The minimum distance is \[k-2^{n-2}(2^{n-d}-1)=2^{n-1}(2^{n-d}-1)-2^{n-2}(2^{n-d}-1)=2^{n-2}(2^{n-d}-1).\qedhere\]
\end{proof}

\subsection{Transitive \texorpdfstring{$\mathcal{C}_1$}{C1}-type codeword stabilisers}
\label{transc1}
In this subsection we prove that \emph{there are no $X$-strongly incidence-transitive codes
$\Gamma \subset \binom{\mathcal{Q}^\epsilon}{k}$ for which $\epsilon=-$ and a codeword stabiliser leaves invariant  an $n$-dimensional totally-isotropic subspace of $V$}. This is the case $(d,\varepsilon)=(n,-)$ not covered by Construction~\ref{ticode} or Theorem~\ref{ticodeissit}. This allows us to complete the proof of Theorem \ref{mainthm} in the case where a codeword stabiliser is reducible on $V$. 

Let $U$ be an $n$-dimensional totally-isotropic subspace of $V$ and let $X_U$ be the stabiliser of $U$ in $X = \syp_{2n}(2)$. Note that since $\dim(U) = n$ we have $U = U^\perp$. For $\varphi \in \mathcal{Q}^-$ the maximum dimension of a  $\varphi$-singular subspace is $n-1$, and therefore Lemma \ref{intersectionthm} implies that $\sing(\varphi) \cap U$ is an $(n-1)$-dimensional subspace of $V$ for all $\varphi \in \mathcal{Q}^-$. 
Recall from Corollary \ref{levi} that $X_U \cong R \rtimes L$, where $R$ is the unipotent radical and $L \cong \gl_n(2)$ is a Levi factor. In particular, setting $n = d$ in Corollary \ref{levi} we have
\[
R = \left\lbrace  \left(
		\begin{array}{ccc}
				I_d     & 0   \\
				X       & I_d
			\end{array}
		\right) \mid X + X^T = 0 \right\rbrace
\]
so $R \cong \mathbb{F}_2^{n(n+1)/2}$ and $|R| = 2^{n(n+1)/2}$.

\begin{lemma}
\label{blocksyscontr}
Let $U$ be a totally-isotropic $n$-dimensional subspace of $V$, with $R, L$ such that $X_U=R\rtimes L$ as in Corollary~\ref{levi},  and let $\mathcal{H}$ denote the set of all $(n-1)$-dimensional subspaces of $U$. For each $H \in \mathcal{H}$  define 
$
P_H := \{ \varphi \in \mathcal{Q}^- \mid \sing(\varphi)\cap U = H \}.
$
\begin{enumerate}[(a)]
\item The subgroup $X_U$ is transitive on $\mathcal{Q}^-$.
\item Let $\varphi \in \mathcal{Q}^-$, $H=\sing(\phi)\cap U$,  and $c\in V$. Then $\varphi_c \in P_H$ if and only if $c\in H$. In particular $|P_H|=2^{n-1}$.

\item $\mathcal{P}=\{P_H \mid H \in \mathcal{H}\}$ is a system of imprimitivity for the action of $X_U$ on $\mathcal{Q}^-$. Moreover the $X_U$-actions on $\mathcal{H}$ and $\mathcal{P}$ are permutationally isomorphic.

\item The transitive actions of $L$ on $\mathcal{P}$, and of $\gl(U)$ on $\mathcal{H}$, are permutationally isomorphic. In particular, $|P_H|=2^{n-1}$ for any $H\in\mathcal{H}$.

\item The kernel of the $X_U$-action on $\mathcal{P}$ is $R$, the group $X_U^\mathcal{P}$ induced by $X_U$ on $\mathcal{P}$ is $L^\mathcal{P}\cong L$.
\end{enumerate}
\end{lemma}

\begin{proof}\leavevmode
\begin{enumerate}[(a)]
\item It follows from \cite[Table 1 and (3.2.4a)]{maximalfactorisations} that $X=X_UX_{\varphi}$, and hence $X_U$ is transitive on $\mathcal{Q}^-$.

\item First note that  since $\varphi\in\mathcal{Q}^-$, the $n$-dimensional subspace $U$ cannot be contained in $\sing(\varphi)$, and it follows from Lemma~\ref{intersectionthm}(a) that $H$ has dimension $n-1$, and $P_H$ is well-defined and non-empty.

Assume first that $c\in H$. Since $H$ is $\phi$-singular, the quadratic form $\phi_c$ is of minus type by Corollary \ref{formvectorbijection}. Thus 
$H_c:=\sing(\varphi_c)\cap U$ has dimension $n-1$, and hence  $H_c\in\mathcal{H}$.
We claim $H_c=H$. Since $H$ and $H_c$ have the same dimension, it is enough to show that $H\subseteq H_c$, that is that $H$ is $\phi_c$-singular and contained in $U$. The latter statement follows from the definition of $H$. For the former statement we compute $\phi_c(h)$ for $h\in H$:
\[\phi_c(h)=\phi(h)+B(h,c)^2=B(h,c),
\]
since $H\subseteq \sing(\phi)$. Since $c\in U=U^\perp$ and $h\in U$, we have that $B(h,c)=0$. This proves the claim.

Now assume $\varphi_c \in P_H$. In particular $\varphi_c\in\mathcal{Q}^-$ and so $c\in\sing(\phi)$ by Corollary \ref{formvectorbijection}. Moreover, $\sing(\phi_c)\cap U=H$. For any $h\in H$, we have $0=\phi_c(h)=B(h,c)$ and so $c\in H^\perp$.
Since $X_U$ is transitive on $\mathcal{Q}^-$ by part (a), we may pick without loss of generality $\phi=\phi_0^-$  from \eqref{eq-forms}. Then $H=\langle e_1,e_2,\ldots,e_{n-1}\rangle$ and $H^\perp=H\perp \langle e_n, f_n\rangle$. Then $c=h'+x_ne_n+y_nf_n$ for some $h'\in H$. Now $0=\phi(c)=\phi(h')+\phi(x_ne_n+y_nf_n)=x_ny_n+x_n+y_n=(x_n+1)(y_n+1)+1$. It follows that $x_n=y_n=0$ and so $c\in H$.

Therefore $|P_H|=|H|=2^{n-1}$ by Corollary \ref{formvectorbijection}.

\item  Since $\phi\in \mathcal{Q}^-$ uniquely determines the element $\sing(\varphi)\cap U$ of $\mathcal{H}$, each $\varphi\in\mathcal{Q}^-$ lies in a unique $P_{H}\in\mathcal{P}$, and $\mathcal{P}$ is a partition of $\mathcal{Q}^-$. 

Let $g\in X_U$, $\varphi\in\mathcal{Q}^-$, and $H:=\sing(\varphi)\cap U\in\mathcal{H}$. Then $x\in H^g$ if and only if $x^{g^{-1}}\in H$, or equivalently $x\in U^g=U$ and $0=\varphi(x^{g^{-1}}) = \varphi^g(x)$, that is, $x\in \sing(\varphi^g)\cap U$.   
Thus $(\sing(\varphi)\cap U)^g= \sing(\varphi^g)\cap U$.

Now suppose that $\varphi', \varphi''\in P_{H'}$, for some $H'\in\mathcal{H}$. Then $H'=\sing(\varphi')\cap U=\sing(\varphi'')\cap U$, and we have just proved that $(H')^g 
=\sing(\varphi'^g)\cap U=\sing((\varphi'')^g)\cap U\in\mathcal{H}$. It follows that $\phi'^g,\phi''^g\in P_{(H')^g}$, and hence that $(P_{H'})^g=P_{(H')^g}$. Thus $\mathcal{P}$ is $X_U$-invariant, and also the  $X_U$-actions on $\mathcal{H}$ and on $\mathcal{P}$ are permutationally isomorphic.

\item We construct an explicit permutational isomorphism. The basis $\mathscr{B}$ of \eqref{basis} is 
\[\{e_1,\dots,e_n,f_1,\dots,f_n\}\] 
(since $d=n$ here), and $U = \langle e_i \mid 1 \leq i \leq n \rangle$. We identify each $v\in V$ uniquely with an ordered pair $(u,u')$, where $u \in U$ and $u' \in U' = \langle f_i \mid 1 \leq i \leq n \rangle$. By Corollary \ref{levi},
\[
L = \left\{  \ell(a) =  \begin{pmatrix}
a & 0 \\
0 & a^{-T}
\end{pmatrix} \mid  a\in {\rm GL}_n(2) \right\}\cong\mathrm{GL}_n(2),
\]
and we note that $(u,u')^{\ell(a)} = (ua,u'a^{-T})$. Therefore $L$ acts on $U$ in the natural action of $\gl_n(2)$ on $U$. In particular, $L$ acts on $\mathcal{H}$ in the natural action of $\gl_n(2)$ on the hyperplanes of $U$.
Hence $(P_H)^{\ell(a)}=P_{H^a}$. Thus $\ell:a\to \ell(a)$ and $f:H\to P_H$ define a group isomorphism $\ell:\gl_n(2)\to L$, and bijection $f:\mathcal{H}\to \mathcal{P}$, such that $(\ell,f)$ is a permutational isomorphism from the $\gl_n(2)$-action  on $\mathcal{H}$ to the $L$-action on $\mathcal{P}$. Since $\gl_n(2)$ is transitive  on $\mathcal{H}$, both actions are transitive. Then all sets $P_H$ have the same cardinality which is $2^{n-1}$ by part (b).

\item By Corollary \ref{levi}, 
\[
R = \left\{  r(Q) =  \begin{pmatrix}
I_n & 0 \\
Q & I_n
\end{pmatrix} \mid  Q=Q^T \right\},
\] which acts trivially on $U$. Thus for each symmetric $Q$ and each $H\in\mathcal{H}$,
$H^{r(Q)} =  H$,
and hence $(P_H)^{r(Q)}=P_{H}$ by part (c). Thus $R$ is contained in the kernel $K$ of the $X_U$-action on $\mathcal{P}$, and since $X_U=R\rtimes L$, we have $X_U=KL$ and so the group $(X_U)^\mathcal{P}=X_U/K \cong L/(L\cap K)$. It follows from part (d) that $K=R$ and $(X_U)^\mathcal{P}=L^\mathcal{P}\cong L$. \qedhere
\end{enumerate}
\end{proof}

Lemma \ref{blocksyscontr} enables us to prove that there are no examples in $J(\mathcal{Q}^-, k)$ with codeword stabilisers leaving invariant a maximal totally isotropic subspace.

\begin{theorem}
\label{almostsit}
Let $X=\syp_{2n}(2)$, let $\Gamma$ be a code in $J(\mathcal{Q}^-,k)$ such that, for $\Delta\in\Gamma$, $X_\Delta\leq X_U$ for some totally isotropic   $n$-dimensional subspace $U$ of $V$. Then $\Gamma$ is not $X$-strongly incidence-transitive. 
\end{theorem}
\begin{proof}
Suppose to the contrary that $\Gamma$ is $X$-strongly incidence-transitive.
Since $X_U$ is transitive on $\mathcal{Q}^-$ by Lemma \ref{blocksyscontr} but $X_\Delta$ is not,  we have a proper containment $X_\Delta < X_U$. It then follows, from Lemmas~\ref{unionsofblocks} and~\ref{blocksyscontr}, that $\Delta$ is a union of some of the parts of the $X_U$-invariant partition $\mathcal{P}$ defined in Lemma~\ref{blocksyscontr}, say $\Delta = \bigcup_{H \in \mathcal{S}}P_H$, where $\mathcal{S} \subset \mathcal{H}$. 
Then $k=2^{n-1}\widehat{k}$, where $\widehat{k}=|\mathcal{S}|$ satisfies $1\leq \widehat{k}\leq |\mathcal{P}|-1$. Recall that $X_U=R\rtimes L$ and that the group $X_U^\mathcal{P}=L^\mathcal{P}\cong L$, by Lemma~\ref{blocksyscontr}. 
Let $\Gamma_0=\{\Delta^\ell\mid \ell\in L\}$ (a possibly  proper subcode of $\Gamma$).
Then for each $\Delta_0\in\Gamma_0$ we have $X_{\Delta_0}<X_U$. Let $\widehat{\Delta}_0$ denote the $\widehat{k}$-subset of parts of $\mathcal{P}$ contained in $\Delta_0$, and let $\widehat{\Gamma}_0 = \{\widehat{\Delta}_0\mid \Delta_0\in\Gamma_0\}$, so $\widehat{\Gamma}_0 \subset \binom{\mathcal{P}}{\widehat{k}}$. Since $R$ acts trivially on $\mathcal{P}$, it follows that $\Gamma_0$ is $X_U$-invariant, and $\widehat{\Gamma}_0$ is $L$-invariant. Moreover, $L$ is transitive on $\widehat{\Gamma}_0$, and the stabiliser in $L$ of $\widehat{\Delta}$ is $L_\Delta$.

We claim that $\widehat{\Gamma}_0$ is $L$-strongly incidence-transitive. By Definition~\ref{defsit} and since $L$ is transitive on $\widehat{\Gamma}_0$, it is sufficient to prove that $L_\Delta$ is transitive on $\widehat{\Delta}\times (\mathcal{P}\setminus \widehat{\Delta})$. Let $P_{H_1}, P_{H_2} \in \widehat{\Delta}$ and  $P_{H_1'}, P_{H_2'} \in \mathcal{P}\setminus \widehat{\Delta}$ (not necessarily distinct), and for $i=1, 2$, let $\varphi_i\in P_{H_i}$ and $\varphi_i'\in P_{H_i'}$. Since $\Gamma$ is $X$-strongly incidence-transitive, $X_\Delta$ acts transitively on $\Delta\times (\mathcal{Q}^-\setminus \Delta)$, and hence there exists $g\in X_\Delta$ such that $(\varphi_1, \varphi_1')^g=(\varphi_2, \varphi_2')$. Since $\mathcal{P}$ is an $X_U$-invariant  partition of $\mathcal{Q}^-$, this implies that  
$(P_{H_1}, P_{H_1'})^g=(P_{H_2}, P_{H_2'})$. Now $g=r\ell$ for some $r\in R$ and $\ell\in L$, and since $r$ acts trivially on $\mathcal{P}$ (by Lemma~\ref{blocksyscontr}) it follows that $(P_{H_1}, P_{H_1'})^\ell=(P_{H_2}, P_{H_2'})$. This proves the claim.

By Lemma \ref{blocksyscontr}, the action of $L$ on $\mathcal{P}$ is permutationally isomorphic to the action of $\gl_n(2)$ on the $(n-1)$-dimensional subspaces of $U$, and this, in turn, is permutationally isomorphic to the action of $\gl_n(2)$ on the $1$-dimensional subspaces of $U$, that is on the non-zero vectors of $U$. Therefore  we may identify $\mathcal{P}$ with the set of $1$-dimensional subspaces of the vector space $U=\mathbb{F}_2^n$, or in other words with the projective space $\mathrm{PG}(n-1,2)$.

We now claim that   $\widehat{\Delta}$ or $\mathcal{P}\setminus \widehat{\Delta}$ is a projective subspace of $\mathcal{P}$. 
If $\widehat{k}=1$ or  $|\mathcal{P}|-1$, then $\widehat{\Delta}$ or its complement is a projective $1$-space so the claim holds. 
Now assume $2\leq \widehat{k} \leq|\mathcal{P}|-2$. In particular $n\geq 3$ since, for $n=2$, $|\mathcal{P}|=3$. If $\widehat{k}=2$ or  $|\mathcal{P}|-2$,  then $\widehat{\Gamma_0}$ cannot be $L$-strongly incidence-transitive, as $L_{{\Delta}}$ will fix a point of $\mathcal{P}$, namely the third point of the projective line containing $\widehat{\Delta}$ or its complement. This is a contradiction, so $3\leq \widehat{k} \leq|\mathcal{P}|-3.$
Then $\widehat{\Gamma}_0$ is equivalent to a `projective type code' analysed in \cite[Section 7]{liebprae}, and \cite[Proposition 7.4]{liebprae} applies. Since $q = 2$, it follows, interchanging  $\widehat{\Delta}$ and ${P}\setminus \widehat{\Delta}$ if necessary,  that either $\widehat{\Delta}$ is a projective subspace of $\mathcal{P}$, or $\widehat{\Delta}$ is a $[0,2,3]_1$-subset, that is to say, each 
projective line meets $\widehat{\Delta}$ in $0, 2$, or $3$ projective points. Since the projective lines have size $3$, the second possibility corresponds to the complement $\mathcal{P}\setminus \widehat{\Delta}$ being a  $[0,1,3]_1$-subset and hence a projective subspace. 

Thus the claim is proved in all cases, and hence the stabiliser $L_\Delta$ of a codeword is equal to the stabiliser $L_W$ of a $d$-dimensional subspace $W$ of $U$, where $1\leq d<n$. By Lemma \ref{blocksyscontr}, $R$ fixes $\mathcal{P}$ pointwise and therefore $R$ leaves $\Delta$ and $W$ invariant, and we have $X_\Delta=R\rtimes L_\Delta = X_{U,W}$. Since $d<n$, this implies that $X_\Delta$ is properly contained in $X_W$. However, by Theorem~\ref{ticodeissit}, if the stabiliser $X_\Delta$ of a codeword $\Delta$ of an $X$-strongly incidence-transitive code leaves invariant a $d$-dimensional totally isotropic subspace $W$ with $d<n$, then $X_\Delta$ is equal to $X_W$. However $X_W$ does not leave invariant an $n$-dimensional subspace. This contradiction completes the proof.   
\end{proof}

\section{Codes with irreducible codeword stabilisers}

Here we examine $X$-strongly incidence-transitive codes in the case where a codeword stabiliser is contained in a maximal irreducible geometric subgroup as in Table~\ref{syp2n2max}. First we obtain a lower bound for the size $k=|\Delta|$ of a codeword, and then we consider in separate subsections the cases where the maximal geometric subgroup is of type $\mathcal{C}_2$, $\mathcal{C}_3$ and $\mathcal{C}_8$.

\label{sec:irred}
\begin{lemma}
	\label{smallinequality}
	Let $X = \syp_{2n}(2)$ and let $\Gamma \subset \binom{\mathcal{Q}^\varepsilon}{k}$ be an $X$-strongly incidence-transitive code with $2 \leq k \leq |\mathcal{Q}^\epsilon|-2$. 
	If $X_\Delta$ acts irreducibly on $V = \mathbb{F}_2^{2n}$ for some $\Delta \in \Gamma$ then $n \geq 3$ and $2n+1\leq k \leq |\mathcal{Q}^\varepsilon|-(2n+1)$.
\end{lemma}

\begin{proof} Since $2 \leq k \leq |\mathcal{Q}^\epsilon|-2$, we must have $n\geq 2$. 
%
First assume that $2 \leq k \leq  \frac{1}{2}|\mathcal{Q}^\epsilon|=2^{n-2}(2^n+\epsilon)$.
Suppose $\varphi \in \Delta$ and let $C = \lbrace c \in V \mid \varphi_c \in \Delta\rbrace$, where we note that $\phi=\phi_0$ (see Remark~\ref{rem3.5}) so $0\in C$. Let $U$ denote the intersection of all subspaces in $V$ which contain $C$, that is the subspace of $V$ spanned by $C$.
We claim that  $X_\Delta$ preserves $U$. Let $g\in X_\Delta$. Then $\phi^g\in \Delta$, so $\phi^g = \phi_0^g=\phi_d$ for some $d\in C$. Then by Corollary \ref{affineactioncor}, for all $c\in C$ we have $\phi_c^g=\phi_{cg+d}$, and hence $cg+d\in  C$. Now $C\subseteq U$ by the definition of $U$, and hence, $Cg+d\subseteq U$, and since $d\in  C\subseteq U$, also $Cg\subseteq U$. Since $C$ spans $U$ it follows that $Ug=U$, proving the claim. 
%
Now $X_\Delta$ acts irreducibly on $V$ and $|C| \geq 2$, so $U = V$ and $\langle C \rangle = V$. Therefore $C\setminus\{0\}$ contains a spanning set for $V$,    which implies that $k=|\Delta|=|C| \geq 2n+1$. Note this is not possible if $(n,\varepsilon) = (2,-)$, since we assumed  $k=|C| \leq  \frac{1}{2}|\mathcal{Q}^-|=3$.

Now assume that $ \frac{1}{2}|\mathcal{Q}^\epsilon| \leq k\leq |\mathcal{Q}^\epsilon|-2 $. The complementary code $$\overline{\Gamma}\subseteq J(\mathcal{Q}^\epsilon,|\mathcal{Q}^\epsilon|-k)=J(\mathcal{Q}^\epsilon,\overline{k})$$ is also $X-$strongly incidence-transitive and satisfies $2 \leq \overline{k} \leq  \frac{1}{2}|\mathcal{Q}^\epsilon|$. Therefore, by the previous argument $(n,\varepsilon) \neq (2,-)$ and $\overline{k}\geq 2n+1$. Thus $k\leq |\mathcal{Q}^\epsilon|-(2n+1)$. 

In either case we conclude that $2n+1\leq k \leq |\mathcal{Q}^\varepsilon|-(2n+1)$ and $(n,\epsilon)\ne (2,-)$.
Finally suppose that  $(n,\varepsilon) = (2,+)$. Then $|\mathcal{Q}^+|=10$ and hence $k=5$, so $|\Delta| =|\overline{\Delta}|= 5$. However this means that
%
$|\Delta \times \overline{\Delta}| = 25$, which does not divide $|\syp_4(2)|$, and hence 25 does not divide $|X_\Delta|$. This contradicts the assumption that the code is  $X$-strongly incidence-transitive. Therefore $n \geq 3$.
\end{proof}

\subsection{\texorpdfstring{$\mathcal{C}_2$}{C2}-type codeword stabilisers}
\label{imprimitivesection}

Here we compute the orbits of the maximal $\calc_2$-subgroups of $\syp_{2n}(2)$ in $\mathcal{Q}^\varepsilon$ and apply the results to the classification of strongly incidence-transitive codes. A maximal $\calc_2$-subgroup of $\syp_{2n}(2)$ is the stabiliser of a direct sum decomposition $V = \oplus_{i=1}^t V_i$ into nondegenerate subspaces $V_i$ of (even) dimension $\dim(V_i) = 2n/t$, with $t\geq 2$ and $t\mid n$, and $V_i$ is   orthogonal to $V_j$ for all $i \neq j$. We denote such a decomposition by $\mathcal{D}$, and the subgroup of $\syp_{2n}(2)$ which preserves $\mathcal{D}$ by $X_\mathcal{D}$, so $X_\mathcal{D}\cong \syp_{2n/t}(2)\wr S_t$.

\begin{definition}
	Let $\mathcal{D}$ denote a decomposition $V = \oplus_{i=1}^t V_i$ and let $\varphi: V \rightarrow \F_2$ be a nondegenerate    quadratic form of type $\varepsilon$ on $V$. Denote by $\varepsilon_i = sgn(V_i)$ the type of $\varphi |_{V_i}$, and write $\varepsilon_\mathcal{D} = (\varepsilon_1, \dots, \varepsilon_t)$. Let $\mathcal{E}(\varphi)$ denote the number of $i$ such that $\varepsilon_i = -$.
\end{definition}

Note that, by Lemma \ref{uniqueform}, we have $\epsilon=\prod_{i=1}^t \epsilon_i$, so if $\epsilon=+$ then $\mathcal{E}(\varphi)$ is even, and if $\epsilon=-$ then $\mathcal{E}(\varphi)$ is odd.

\begin{lemma}
	\label{dstaborbs} 
	The $X_\mathcal{D}-$orbits in $\mathcal{Q}^\varepsilon$ are the subsets $\mathcal{O}_m = \lbrace \varphi \in \mathcal{Q}^\varepsilon \mid \mathcal{E}(\varphi) = m \rbrace$ for all integers $m$ such that $0 \leq m \leq t$ and $(-1)^m = \varepsilon$. In particular if $X_\mathcal{D}$ has at most two orbits in $\mathcal{Q}^\varepsilon$, then one of the following holds for $t$ and the $X_\mathcal{D}$-orbits in $\mathcal{Q}^\varepsilon$:
	\begin{enumerate}[(i)]
		\item $t=2$ and $\varepsilon = +$, with orbits $\mathcal{O}_0$ and $\mathcal{O}_2$,
		\item $t=3$ and $\varepsilon = \pm$, with orbits $\mathcal{O}_0$ and $\mathcal{O}_2$ if $\varepsilon = +$, or $\mathcal{O}_1$ and $\mathcal{O}_3$ if $\varepsilon = -$,
		\item $t=4$ and $\varepsilon = -$, with orbits $\mathcal{O}_1$ and $\mathcal{O}_3$,
		\item $t=2$ and $\varepsilon = -$, and in this case $X_\mathcal{D}$ acts transitively and $\mathcal{O}_1 = \mathcal{Q}^-$; moreover $X_{\mathcal{D}, V_1}$ has two orbits in $\mathcal{Q}^-$, namely the sets $\mathcal{O}_{+,-}$ and $\mathcal{O}_{-,+}$ of forms with type $(+,-)_\mathcal{D}$, and  $(-,+)_\mathcal{D}$, respectively.
	\end{enumerate}	
	\end{lemma}

\begin{proof}
	For every $\varphi \in \mathcal{Q}^\varepsilon$ and every $i$ with $1 \leq i \leq t$, the group $\syp(V_i) \cong \syp_{2n/t}(2)$ preserves the type $\varepsilon_i$ of $\varphi_i = \varphi|_{V_i}$. Moreover, the top group $S_t$ of $X_\mathcal{D}$ permutes the components $V_i$ of the decomposition $\mathcal{D}$ and the corresponding direct factors of the base group $\syp_{2n/t}(2)^t$ of $X_{\mathcal{D}}$, and therefore permutes the entries of the vector $(\varepsilon_1, \dots, \varepsilon_t)$ while preserving $\mathcal{E}(\varphi)$. This implies the subsets $\mathcal{O}_m$ are $X_\mathcal{D}-$invariant. It remains to show that $X_\mathcal{D}$ acts transitively on each non-empty subset $\mathcal{O}_m$. Let $\varphi, \varphi' \in \mathcal{O}_m$. By the comment above   $\varepsilon = \Pi_{i=1}^t \varepsilon_i = \Pi_{i=1}^t \varepsilon_i'$, where $\varepsilon_i $ the type of $\varphi |_{V_i}$ and  $\varepsilon'_i $ the type of $\varphi' |_{V_i}$. Since $\mathcal{E}(\varphi) = \mathcal{E}(\varphi')$ and $S_t$ acts transitively on the $m$-subsets of $\lbrace 1, \dots , t \rbrace$, there exists $\sigma \in S_t$ such that $\varepsilon_{i^{\sigma^{-1}}} = \varepsilon_i '$ for all $i \leq t$, and by Lemma \ref{uniqueform}, we can represent $\varphi^\sigma$ uniquely as $\varphi_{1^{\sigma^{-1}}}\oplus \dots \oplus \varphi_{t^{\sigma^{-1}}}$. Further, for each $i$,  $\syp(V_i)$ acts transitively on the quadratic forms on $V_i$ of type $\varepsilon_i' = \varepsilon_{i^{\sigma^{-1}}}$, so there exists $g_i \in \syp(V_i)$ such that $\varphi_{i^{\sigma^{-1}}}^{g_i} = \varphi'_i$. Therefore $\sigma(g_1, \dots, g_t) \in X_\mathcal{D}$ maps $\varphi$ to $\varphi'$ and the non-empty $\mathcal{O}_m$ are the $X_\mathcal{D}$-orbits in $\mathcal{Q}^\varepsilon$. Finally $\mathcal{O}_m \neq \emptyset$ if and only if $(-1)^m = \varepsilon$ by Lemma \ref{uniqueform}.
	
	Checking the number of integers $m$ satisfying $(-1)^m= \varepsilon$ with $0 \leq m \leq t$, we see that if $\varepsilon = +$ and $t \geq 4$, then $\mathcal{O}_m$ is non-empty for $m \in \lbrace 0,2,4 \rbrace$, while if $\varepsilon = -$ and $t \geq 5$, then $\mathcal{O}_m$ is non-empty for $m \in \lbrace 1,3,5 \rbrace$. Further, if $X_\mathcal{D}$ has at most two orbits in $\mathcal{Q}^\varepsilon$, then one of (i)--(iv) holds, and in case (iv), the $X_{\mathcal{D}, V_1}$-orbits are the sets $\mathcal{O}_{+,-}$ and $\mathcal{O}_{-,+}$ described there. 
\end{proof}

\begin{theorem}	\label{c2smallcases}
Let $X=\syp_{2n}(2)$ for $n\geq 2$, let $\Gamma$ be an $X$-strongly incidence-transitive code in $J(\mathcal{Q}^\varepsilon,k)$ 
and $\Delta\in\Gamma$, such that $1<|\Delta|<|\mathcal{Q}^\epsilon -1$ and $X_\Delta$ is irreducible on $V$. Suppose that $X_\Delta\leq X_\mathcal{D}$ for some direct sum decomposition $\mathcal{D}$ of $V = \oplus_{i=1}^t V_i$, as above. 
Then, $t=2$, and   $\Gamma=\Gamma(n,n/2,+,+)$ or $\Gamma(n,n/2,+,-)$ as in Construction \ref{ndcode}, with $\Delta=\mathcal{O}_0$ or $\mathcal{O}_2$, respectively, and $X_\Delta=X_{\mathcal{D}}$, as in case (i) of Lemma~\ref{dstaborbs}. 
\end{theorem}

\begin{proof}
Let $X, \Gamma$, $\Delta$, and $\mathcal{D}$ be as in the statement, with $X_\Delta\leq X_{\mathcal{D}}$. Then $X_\mathcal{D}$ has at most two orbits in $\mathcal{Q}^\varepsilon$, and hence one of the cases (i)--(iv) of Lemma~\ref{dstaborbs} holds. We use the definition of the $\mathcal{O}_i$ in Lemma \ref{dstaborbs}. We note first that, in cases (i)--(iii), $X_{\mathcal{D}}$ has two orbits in $\mathcal{Q}^\epsilon$, namely $\mathcal{O}_i, \mathcal{O}_j$, for certain $i, j$. Since $X_\Delta$ also has two orbits in $\mathcal{Q}^\epsilon$, it follows that $\Delta=\mathcal{O}_i$ or $\mathcal{O}_j$ and $X_\Delta=X_{\mathcal{D}}$. 
We consider the four cases separately.

	In case (i), $n$ is even, $\epsilon=+$, and $V = V_1 \oplus V_2$ with $V_1$ and $V_2$ both nondegenerate of dimension $n$ and $V_2 = V_1^\perp$. Also $X_\Delta=X_{\mathcal{D}}$, and $\Delta$ is $\mathcal{O}_0$ or $\mathcal{O}_2$. 
	Suppose first that $\Delta = \mathcal{O}_0$ and $\overline{\Delta} = \mathcal{O}_2$, so each  $\varphi \in \Delta$ has the form $\varphi=\varphi_1\oplus\varphi_2$ with each $\varphi_i=\varphi|_{V_i}$ of type $+$, in other words $\phi$ has type $(+,+)_\mathcal{D}$. Now if $\phi_{V_1}$ has type $+$ then $\phi_{V_2}$ must have type $+$ too by Lemma \ref{uniqueform}. Thus $\Delta=\Delta(V_1)$ for $\epsilon=\epsilon'=+$ as defined in Construction \ref{ndcode}, and $\Gamma=\Gamma(n,n/2,+,+)$. 
	Similarly, if $\Delta = \mathcal{O}_2$ so $\overline{\Delta} = \mathcal{O}_0$,  then $\Gamma$ is the complementary code $\Gamma(n,n/2,+,-)$ by the comment after Construction \ref{ndcode},  and $\Gamma=\Gamma(n,n/2,+,-)$.
	
	In case (ii), $V = V_1 \oplus V_2 \oplus V_3$ and $\varepsilon \in \{+,-\}$. Replacing $\Gamma$ by its complementary code if necessary (see the remarks on \cite[page 6]{liebprae}), we may assume that $\Delta =\mathcal{O}_0$ if $\varepsilon=+$ and $\Delta =\mathcal{O}_3$ if $\varepsilon=-$. Therefore, in either case, $\Delta$ consists of all quadratic forms in $\mathcal{Q}^\varepsilon$ satisfying $\varepsilon_1 = \varepsilon_2 = \varepsilon_3 = \varepsilon$, while $\overline{\Delta}$ consists of all quadratic forms satisfying $\lbrace \varepsilon_1, \varepsilon_2, \varepsilon_3 \rbrace = \lbrace \varepsilon,-\varepsilon,-\varepsilon \rbrace$. 
	Note that, if $\varepsilon = -$ and $\dim(V_i) =2$, then the number of choices for $\varphi|_{V_i}$, where $\varphi\in \Delta$, is $1$ (the size of $\mathcal{Q}^-$ when the dimension is 2) and hence $|{\Delta}| =1$, contradicting our assumption that $|\Delta|>1$. Thus if $\varepsilon = -$ then $\dim(V_i) \geq 4$. 
	It follows that we may choose forms $\varphi \in \Delta$ and $\mu \in \overline{\Delta}$ such that $\mu$ has type $(\varepsilon, -\varepsilon, -\varepsilon)$ and $\varphi|_{V_1} \neq \mu|_{V_1}$. Since $X_\Delta = X_\mathcal{D}$ acts transitively on $\Delta \times \overline{\Delta}$, Lemma \ref{uniqueform} allows us to define uniquely a third form $\nu \in \overline{\Delta}$ such that $\nu|_{V_1} = \varphi|_{V_1}$ and $\nu|_{V_2} = \mu|_{V_2}$, $\nu|_{V_3} = \mu|_{V_3}$. Since $X_\mathcal{D}$ acts transitively on $\Delta \times \overline{\Delta}$, there exists $g \in X_{\mathcal{D},\varphi}$ such that $\mu^g = \nu$,  in particular, $\mu^g|_{V_1} = \varphi|_{V_1}$ has type $\varepsilon$. Since for both $\mu$ and $\nu$, the only restriction $\mu^g|_{V_i}, \varphi|_{V_i}$ of type $\varepsilon$ is the restriction to $V_1$, it follows that $g$ fixes $V_1$. Hence
	$g$ maps $\mu|_{V_1}$ to $\nu|_{V_1} = \varphi|_{V_1}$. However, this is a contradiction since  $\mu|_{V_1} \ne \varphi|_{V_1}$, and $g$ fixes $\varphi$ and $V_1$. So case (ii) gives no examples.

	Case (iii)  of Lemma \ref{dstaborbs} is dealt with in a similar manner to case (ii), but the details need care. Here $X_\Delta = X_\mathcal{D}$ and $t=4$,
	$\varepsilon=-$, and we suppose first that $\Delta = \mathcal{O}_1$, $\overline{\Delta} = \mathcal{O}_3$. Select $\varphi \in \Delta$ and $\mu \in \overline{\Delta}$ with respective types $(+,+,+,-)_\mathcal{D}$ and $(+,-,-,-)_\mathcal{D}$ such that $\varphi|_{V_1} \neq \mu|_{V_1}$ (this is possible since even if $\dim(V_1)=2$ there are three choices for this restriction). Using Lemma \ref{uniqueform} we define a unique quadratic form $\nu \in \mathcal{Q}^-$ by $\nu = \varphi_1 \oplus \mu_2 \oplus \mu_3 \oplus \mu_4$. By definition, $\nu \in \overline{\Delta}$, and hence by assumption there exists $g\in X_{\mathcal{D},\varphi}$ such that $\mu^g = \nu$; in particular $g$ fixes $V_1$ since $V_1$ is the only component of $\mathcal{D}$ where both $\mu$ and $\nu$ have $+$ type restriction, and since $g$ also fixes $\varphi$, it follows that g fixes $\varphi|_{V_1}$.  However $g$ maps $\mu|_{V_1}$ to $\nu|_{V_1}=\varphi|_{V_1}$, and this is a contradiction.
		If instead $\Delta = \mathcal{O}_3$ and $\overline{\Delta} = \mathcal{O}_1$,  then we use the argument just given on the complementary code. Thus there are no examples in case (iii).

In case (iv), $V = V_1 \oplus V_2$, $\varepsilon = -$, and $X_\mathcal{D} \cong \syp_n(2) \wr C_2$ is transitive on $\mathcal{Q}^-$, while $X_\Delta<X_\mathcal{D}$. Also $\mathcal{O}_{+,-}$ and $\mathcal{O}_{-,+}$ are the two orbits  in $\mathcal{Q}^-$ of the index $2$ subgroup  $X_{V_1}$ of $X_\mathcal{D}$, and these two sets form a system of imprimitivity for the action of $X_\mathcal{D}$ on $\mathcal{Q}^-$.
	By Lemma \ref{unionsofblocks}, $\Delta$ is a union of blocks of this imprimitivity system, and hence  $\Delta=\mathcal{O}_{-,+}$ or $\mathcal{O}_{+,-}$. In either case $X_{V_1}$ leaves $\Delta$ invariant, so $X_{V_1}\leq X_\Delta <X_\mathcal{D}$. It follows that $X_\Delta=X_{V_1}$. This however contradicts our assumption that $X_\Delta$ is irreducible on $V$, and so there are no examples in case (iv).  
\end{proof}

\subsection{\texorpdfstring{$\mathcal{C}_3$}{C3}-type codeword stabilisers}
\label{spreadsection}
In Section \ref{spreadsection} we show that there do not exist any strongly incidence-transitive codes of the type described in Theorem \ref{mainthm} whose codeword stabilisers are irreducible on $V$ and contained in a $\mathcal{C}_3$-subgroup of $\syp_{2n}(2)$. We open with some notes on finite fields and their automorphisms. Further details are available in \cite{handbookff}.
\begin{remark}
\label{finitefields}
Let $b \geq 1$ and $q=2^b$. There are $b$ distinct field automorphisms of $\mathbb{F}_q$, namely $\sigma_j$, for $0 \leq j \leq b-1$, where $\sigma_j(\alpha) = \alpha^{2^j}$ for all $\alpha \in \mathbb{F}_q$. The automorphism group of $\mathbb{F}_q$, denoted $\aut(\mathbb{F}_q)$, is cyclic  of order $b$ and is generated by $\sigma_1$. For $1\leq j\leq b-1$, the set of fixed points of $\sigma_j$ in $\mathbb{F}_{q}$ is the subfield $\mathbb{F}_{2^{(j,b)}}$ where $(j,b)$ denotes the greatest common divisor of $j$ and $b$. In particular $\aut(\mathbb{F}_q)$ fixes $\mathbb{F}_2$ pointwise. 

Since $q$ is even, every element of $\mathbb{F}_{q}$ is a square. In particular, for $\lambda \in \mathbb{F}_q^\times$, $\alpha = \lambda^{q/2}$ is the unique element in $\mathbb{F}_q^\times$ satisfying $\alpha^2 = \lambda$. The elements of $\mathbb{F}_q$ lying in the orbit $\lambda^{\aut(\mathbb{F}_q)}$ are called the \emph{conjugates} of $\lambda$ in $\mathbb{F}_q$. By \cite[Lemma 2.1.75]{handbookff}, $|\lambda^{\aut(\mathbb{F}_q)}| = b$ if and only if $\lambda$ is not contained in a proper subfield of $\mathbb{F}_q$. In particular, if $b$ is prime then $|\lambda^{\aut(\mathbb{F}_q)}| = b$ for all $\lambda \in \mathbb{F}_q \setminus \mathbb{F}_2$.

The \emph{absolute trace} $\Tr: \mathbb{F}_{2^b} \rightarrow \mathbb{F}_2$ is the $\mathbb{F}_2$-linear map defined by \[\Tr(\alpha)=  \sum_{i=0}^{b-1}\alpha^{2^i}=\alpha+\alpha^2+\alpha^4+\ldots+\alpha^{q/2}.\]
 The absolute trace map $\Tr$ is additive and invariant under field automorphisms so in particular $\Tr(\alpha^2) = \Tr(\alpha)^2$  for all $\alpha \in \F_{2^b}$. 
 Let $K = \ker(\Tr)$ and $K^\# = K \setminus \lbrace 0 \rbrace$. Note that $K$ is a codimension $1$ subspace of $\mathbb{F}_q$ (considered as a vector space over $\mathbb{F}_2$) and therefore $|K| = 2^{b-1}$. Additionally, note that $\Tr(1) = 0$ if and only if $b$ is even, and if $b = 2$ then $K = \mathbb{F}_2$ is fixed pointwise by $\aut(\mathbb{F}_4)$. Finally, if $b$ is an odd prime then by Fermat's Little Theorem we have $2^{b-1} \equiv 1 \pmod{b}$, and $\aut(\mathbb{F}_q)$ has $(2^{b-1}-1)/b$ orbits in $K^\#$, each of which has length $b$. 
\end{remark}

\subsubsection{Forms under field reduction}
\label{formsfieldred}
Let $W$ be a $2m$-dimensional vector space over $\mathbb{F}_q=\mathbb{F}_{2^b}$. Let $V$ denote the set of vectors in $W$ equipped with the same addition operation as in $W$, but with scalar multiplication restricted so that elements of $V$ may only be scaled by elements of $\mathbb{F}_2$. Then $V$ is a $2n$-dimensional vector space over $\mathbb{F}_2$, where $n = mb$. We let $W$ and $V$ denote the vector spaces described above throughout Subsections \ref{formsfieldred} and \ref{c3nocodes}. Moreover, we assume that $W$ is equipped with a symplectic form $\widetilde{B}: W \times W \to \mathbb{F}_{q}$, and we define $B: V \times V \to \mathbb{F}_2$ by $B = \tr\circ \widetilde{B}$, where $\Tr:\mathbb{F}_{q}\to \mathbb{F}_2$ is the absolute trace map defined in Remark~\ref{finitefields}. For $\epsilon \in \{ +,- \}$, we denote by $\mathcal{Q}^\epsilon(W)$ the set of all quadratic forms $\Phi: W \to \mathbb{F}_{q}$ of type $\epsilon$ which polarise to $\widetilde{B}$. We write $\mathcal{Q}(W) = \mathcal{Q}^+(W) \cup \mathcal{Q}^-(W)$, and use similar notation for quadratic forms on $V$.
The following theorem is a special case of \cite[Theorem C]{gill}, with the last assertion covered by the comments in \cite[Section 1]{gill}.
\begin{theorem}[\cite{gill}]
	\label{gilltype}
	If $\widetilde{B}$ is a symplectic form on $W$ then $B = \Tr \circ \widetilde{B}$ is a symplectic form on $V$. Additionally, if $\Phi$ is a quadratic form of type $\epsilon \in \{+,- \}$ on $W$, then $\phi = \Tr \circ \Phi$ is a quadratic form of type $\epsilon$ on $V$, and if $\widetilde{B}$ is the  polar form of $\Phi$, then $B$ is the  polar form of $\phi$. 
\end{theorem}

Recall from Notation \ref{notation} that if $\Phi$ is a quadratic form on a vector space $W$ and $c \in W$ then 
$\Phi_c$  is the quadratic form defined by $\Phi_c(w) = \Phi(w) + \widetilde{B}(w,c)^2$ for all $w \in W$.

\begin{lemma}
    \label{tracemaptranslation}
    Let $\Phi,\Psi \in \mathcal{Q}(W)$, and let $\phi = \tr\circ\Phi$ and $\psi = \tr\circ\Psi$. If $c \in W$ and $\Psi = \Phi_c$ then $\psi = \phi_c$.
\end{lemma}
\begin{proof}
We use Equation \eqref{parameq} in combination with properties of the trace map $\Tr$.  For all $v \in V$, $\psi(v)= \tr \circ \Psi(v)= \tr(\Phi(v) + \widetilde{B}(v,c)^2)$ is equal to
\[
\tr(\Phi(v)) + \tr(\widetilde{B}(v,c))^2 
        = \phi(v) + B(v,c)^2
        = \phi_c(v).
\]
Therefore $\psi = \phi_c$.
\end{proof}


\begin{lemma}
	\label{qfbijection}
	Define a mapping $T: \mathcal{Q}(W) \rightarrow \mathcal{Q}(V)$ by
		\begin{align}
	T(\Phi) &= \Tr \circ \Phi.
	\label{c3bijection}
\end{align}
    Then $T$ is a bijection. Moreover,  for each $\epsilon\in\{+,-\}$, $T$ induces a bijection $T:\mathcal{Q}^\epsilon(W) \to \mathcal{Q}^\epsilon(V)$, and for $\Phi \in \mathcal{Q}(W)$ and $c \in W$, $\Phi$ and $\Phi_c$ are of the same type if and only if $\Phi(c) \in \ker(\Tr)$.
\end{lemma}

\begin{proof}
Let $\Phi \in \mathcal{Q}(W)$, let $\widetilde{B}$ denote the polar form of $\Phi$, and let $\phi = T(\Phi)$, and $B=\Tr \circ \widetilde{B}$. Then by Theorem \ref{gilltype}, $\phi \in \mathcal{Q}(V)$, $\phi$ has the same type as $\Phi$, and $B$ is the polar form of $\phi$. In particular  $T(\mathcal{Q}(W))\subseteq \mathcal{Q}(V)$, and 
$T(\mathcal{Q}^\epsilon(W))\subseteq \mathcal{Q}^\epsilon(V)$ for each $\epsilon\in\{+,-\}$. 

To show that 
$T$ is onto, consider an arbitrary $\psi \in \mathcal{Q}(V)$.  By Lemma \ref{param}, there exists a unique $c \in V$ such that $\psi = \phi_c$. Moreover, again by Lemma~\ref{param}, the map $\Psi(x) = \Phi(x) + \widetilde{B}(x,c)^2$, for $x \in W$, is a quadratic form on $W$ which polarises to $\widetilde{B}$, and $\Psi = \Phi_{c}$ (by Notation~\ref{notation}). Then by Lemma \ref{tracemaptranslation} we have $T(\Psi) = T(\Phi_c) = \phi_c = \psi$. Therefore $T$ is onto.

 To show that $T$ is one-to-one, let $\Psi,\Psi' \in \mathcal{Q}(W)$ such that $T(\Psi)=T(\Psi')$.  By Lemma \ref{param} and Notation~\ref{notation}, there exist  unique $c, d \in W$ such that $\Psi=\Phi_c$ and $\Psi'=\Phi_d$, and then by Lemma \ref{tracemaptranslation}, $T(\Psi) = \phi_c$ and $T(\Psi')=\phi_d$. Thus $\phi_c=\phi_d$. It then follows from Lemma \ref{param}  that $c=d$, and hence $\Psi=\Psi'$. Therefore $T$ is one-to-one, and hence a bijection. Thus $|\mathcal{Q}(W)|=  |\mathcal{Q}(V)|$. Finally since, for each $\epsilon$, we have $T(\mathcal{Q}^\epsilon(W))\subseteq \mathcal{Q}^\epsilon(V)$, equality holds and $T$ induces a bijection $\mathcal{Q}^\epsilon(W) \to \mathcal{Q}^\epsilon(V)$.  
 
 Let  $c\in W$. By Theorem \ref{gilltype}, $\Phi$ and $\Phi_c$ are of the same type if and only if $\phi$ and $\phi_{c}$ are of the same type, and by Corollary~\ref{formvectorbijection}, $\phi$ and $\phi_{c}$ are of the same type if and only if $c\in\sing(\phi)$. Since $\phi(c) = \Tr(\Phi(c))$, this condition is equivalent to the condition $\Tr(\Phi(c)) = 0$, that is, $\Phi(c) \in \ker(\Tr)$.
\end{proof}


A \emph{semilinear transformation} of $W$ is a mapping $g: W \to W$ such that $(u + w)g = ug + wg$ for all $u,w \in W$, and there exists a field automorphism $\sigma \in \aut(\mathbb{F}_{q})$ such that $(\alpha w)g = \alpha^\sigma wg$ for all $\alpha \in \mathbb{F}_{q}$ and $w \in W$. The automorphism $\sigma$ depends on $g$ but not on $u,w$ or $\alpha$, and we say that $g$ is semilinear relative to $\sigma$. The group of all nonsingular semilinear transformations of $W$ is denoted $\gaml(W)$. We use these notions to describe the setup and what we will do in this and the next subsection. In particular the subgroup $G$ we introduce in Remark~\ref{remc3setup}(b)(c) is called a $\mathcal{C}_3$-subgroup, and is a maximal $\mathcal{C}_3$-subgroup if $b$ is prime -- hence our assumption on $b$ (see Table \ref{syp2n2max}).  Such subgroups are constructed in \cite[Section 2.2.3]{colva} by computing the intersection of $\syp(V)$ with the normaliser of a Singer subgroup in $\gl(V)$.

\begin{remark} [Our $\mathcal{C}_3$-action setup] \label{remc3setup}
We have $V=\mathbb{F}_2^{2n}$ and $W=\mathbb{F}_{2^b}^{2m}$, where $n=mb$ and $b$ is prime, and we identify the underlying sets. As in \cite[Definition 1.6.4]{colva} and \cite[Definition 1.6.17]{colva}:

\begin{enumerate}[(a)]

\item  A semilinear transformation $g$ on $W$, relative to an automorphism $\sigma$ of $\mathbb{F}_{q}$, is called a \emph{semi-isometry} of $\widetilde{B}$ if, for all $v, w\in W$, $\widetilde{B}(vg, wg)=\widetilde{B}(v,w)^\sigma$.
The set of all semi-isometries  of $\widetilde{B}$ forms a subgroup of $\Gamma{\rm L}(W)$ denoted $\gamsp(W, \widetilde{B}) \cong \syp_{2m}(2^b)\rtimes C_b$.

\item  Given $\Phi \in \mathcal{Q}^\epsilon(W)$, where $\epsilon\in\{+,-\}$, and  $g\in \gamsp(W, \widetilde{B})$ relative to $\sigma$, we define a function $\Phi^g: W \rightarrow \mathbb{F}_q$ by
\begin{equation}
	\label{generalisedaction}
	\Phi^g(v) = \Phi (vg^{-1})^\sigma.
\end{equation} 
Note that this generalises \eqref{action}. 
We check that 
\begin{align*}
\Phi^g(\lambda v)&= \Phi( (\lambda v)g^{-1})^\sigma=\Phi( \lambda^{\sigma^{-1}} (vg^{-1}))^\sigma=((\lambda^{\sigma^{-1}})^2\Phi(  vg^{-1}))^\sigma\\
&=\lambda^2\Phi(  vg^{-1})^\sigma=\lambda^2 \Phi^g(v),
\end{align*}
and 
\begin{align*}
\Phi^g(u+v)-\Phi^g(u)-\Phi^g(v)&= \Phi ((u+v)g^{-1})^\sigma- \Phi (ug^{-1})^\sigma- \Phi (vg^{-1})^\sigma\\
&=(\Phi (ug^{-1}+vg^{-1})- \Phi (ug^{-1})- \Phi (vg^{-1}))^\sigma\\&=
\widetilde{B}(ug^{-1},vg^{-1})^\sigma=\widetilde{B}(u,v).    
\end{align*}
Therefore $\Phi^g\in \mathcal{Q}(W)$, in fact (as we will show in Proposition\ref{c3permisos}(b)) $\Phi^g\in \mathcal{Q}^\epsilon(W)$, and  Equation \eqref{generalisedaction} defines a group action of $\gamsp(W, \widetilde{B})$ on $\mathcal{Q}(W)$ preserving $\mathcal{Q}^+(W)$ and $\mathcal{Q}^-(W)$. 

The element $g \in \gamsp(W, \widetilde{B}) $ is a \emph{semi-isometry of $\Phi$} if $\Phi^g=\Phi$, that is if, for all $v\in W$, $\Phi(vg)=\Phi(v)^\sigma$. 
The set of all semi-isometries  of  $\Phi$  forms a subgroup of $\gamsp(W, \widetilde{B})$,  denoted  $\gamoo^\epsilon(W,\Phi)\cong \oo^\epsilon(W,\Phi)\rtimes C_b$.

\item  We write $G=\gamsp(W, \widetilde{B})$
, and we will show in Proposition~\ref{c3permisos} that $G$ is a subgroup of $X=\syp(V, B)$; it is the maximal $\mathcal{C}_3$-subgroup listed in row (d) of Table~\ref{syp2n2max}. Also for $\Phi\in\mathcal{Q}^\epsilon(W)$, we write $G_\Phi=\gamoo^\epsilon(W,\Phi)$ for the semi-isometry group of $\Phi$; and for $\phi=T(\Phi)$ as in \eqref{c3bijection}, we will show in Proposition~\ref{c3permisos} that the stabiliser $G_\phi$ in $G$ of $\phi$ is equal to $G_\Phi$. 

\item For  a given $\epsilon\in\{+,-\}$, we choose $\Phi
\in\mathcal{Q}^\epsilon(W)$, set $\phi
:=T(\Phi
)\in\mathcal{Q}^\epsilon(V)$, and as in the proof of Lemma~\ref{qfbijection}, for each $c\in V$ we define $\Phi_c
\in\mathcal{Q}(W)$, 
and $\phi_c
\in\mathcal{Q}(V)$, so $\phi_c=T(\Phi_c)$, and we note that 
\[
\phi_c\in\mathcal{Q}^\epsilon(V) \quad \Longleftrightarrow\quad \Phi_c\in\mathcal{Q}^\epsilon(W)
\quad \Longleftrightarrow\quad \Phi
(c)\in \ker(\Tr).
\]

\item For $\lambda\in\mathbb{F}_{q}$, we write $[\lambda] :=\{ \lambda^\sigma \mid \sigma\in\Aut(\mathbb{F}_{q})\}$, and for $a\in V$, we define
    \begin{equation}\label{e:thetalam}
    \theta_{[\lambda]}(a) = \{ \phi_c \mid c\ne a,\ \mbox{and}\ \Phi_a(c+a)\in [\lambda]\},\  
    \end{equation}
and we also write $\theta_\lambda(a)     = \{ \phi_c  \mid c\ne a,\ \mbox{and}\ \Phi_a(c+a) = \lambda \}$. 
We also write simply $ \theta_{[\lambda]}$ and $ \theta_{\lambda}$ for $\theta_{[\lambda]}(0)$ and $\theta_{\lambda}(0)$ respectively.
In Proposition~\ref{suborthorbs}, we will show that, for $a\in V$, the $G_{\phi_a}$-orbits in $\mathcal{Q}\setminus\{\phi\}$ are the sets $\theta_{[\lambda]}(a)$. 
\end{enumerate}
\end{remark}

Note that, using \eqref{polar} and \eqref{parameq}, we have
 \begin{equation}\label{e:phi_a(c+a)}
    	\Phi_a(c+a)= \Phi(c+a) + \widetilde{B}(c+a, a)^2 = \Phi(c)+\Phi(a) + \widetilde{B}(c, a) + \widetilde{B}(c, a)^2.
    \end{equation}

\begin{proposition}
	\label{c3permisos}
	Let $G = \gamsp(W,\widetilde{B})$ and $X=\syp(V,B)$, where $B=\Tr\circ \widetilde{B}$, and let $T: \mathcal{Q}(W)\to \mathcal{Q}(V)$ as in \eqref{c3bijection}. Then $G\leq X$, and if $\phi=T(\Phi)\in\mathcal{Q}(V)$ and $g\in G$, then $\phi^g=T(\Phi^g)$. Moreover,
	\begin{enumerate}
	    	    \item[(a)] if $\epsilon\in\{+,-\}$ and $\phi = T(\Phi)= \tr\circ\Phi$, with  $\Phi\in\mathcal{Q}^\epsilon(W)$, then $\phi\in\mathcal{Q}^\epsilon(V)$ and $G_\phi = G_\Phi = \gamoo^\epsilon(W,\Phi)$; and
	    	    
	    	    \item[(b)] the $G$-actions on $\mathcal{Q}(W)$ and $\mathcal{Q}(V)$ are equivalent; there are two $G$-orbits in each action, namely $\mathcal{Q}^\epsilon(W)$ and $\mathcal{Q}^\epsilon(V)$, respectively, with $\epsilon\in\{+,-\}$.
	\end{enumerate}
\end{proposition}

\begin{proof} 
Let $g\in G$ relative to an automorphism $\sigma$ of $\mathbb{F}_{q}$, and let $v, w\in W$. Then
\[
B(vg, wg) =\Tr(\widetilde{B}(vg, wg)) = \Tr(\widetilde{B}(v, w)^\sigma) = \Tr(\widetilde{B}(v, w)) =B(v,w),
\]
and hence $g\in X$. Thus $G\leq X$. Also if  $\phi=T(\Phi)\in\mathcal{Q}(V)$, $v\in V$, then using \eqref{action} and \eqref{generalisedaction} we see that
\[
\phi^g(v)=\phi(vg^{-1})=\Tr(\Phi(vg^{-1}))=\Tr(\Phi^g(v)^{\sigma^{-1}})=\Tr(\Phi^g(v))=T(\Phi^g)(v),
\]
so $\phi^g=T(\Phi^g)$.

(a) Now  $G_\Phi = \gamoo^\epsilon(W,\Phi)$, and by Lemma~\ref{qfbijection}, $\phi=T(\Phi)\in\mathcal{Q}^\epsilon(V)$. 

To prove the rest of part (a), we first show that $G_\phi$ contains $G_\Phi$. To do this suppose that $g\in G_\Phi$, that is, $g\in G$ is a semi-isometry of $\Phi$ relative to some automorphism $\sigma$. 
Then, as we showed above, $\phi^g=T(\Phi^g)$, so 
\[
\phi^g=T(\Phi^g)=T(\Phi)=\phi,
\]
that is, $g\in G_\phi$, which shows that $G_\Phi\leq G_\phi$. 
Conversely, if $g\in G_\phi$, that is, $\phi^g=\phi$, then 
$T(\Phi)=\phi=\phi^g=T(\Phi^g)$. By Lemma \ref{qfbijection}, $T$ is a bijection, and so $\Phi=\Phi^g$, Thus $G_\phi\leq G_\Phi$. This completes the proof of part (a). 

(b) 
The $G$-orbits in $\mathcal{Q}(W)$ are $\mathcal{Q}^\epsilon(W)$, for $\epsilon\in\{+,-\}$, and it follows from the first assertion of Proposition~\ref{c3permisos}, proved above, that the $G$-actions on $\mathcal{Q}(W)$ and $\mathcal{Q}(V)$ are equivalent. Moreover, by Lemma~\ref{tracemaptranslation} the type of a form is preserved by the bijection $T$. It follows that the $G$-orbits in  $\mathcal{Q}(V)$  are $\mathcal{Q}^\epsilon(V)$ for $\epsilon\in\{+,-\}$.
%
\end{proof}

\subsubsection{Orbits of \texorpdfstring{$\mathcal{C}_3$}{C3}-subgroups in \texorpdfstring{$\mathcal{Q}^\varepsilon$}{Qepsilon} and an  application to codes}
\label{c3nocodes}
We continue with the notation and assumptions from Remark~\ref{remc3setup}, and note that all assertions in parts (a)--(d) of that remark have been proved..
In the main result Theorem \ref{c3conclusion} of this subsection we prove that there are no strongly incidence-transitive codes in $J(\mathcal{Q}^\epsilon(V),k)$ for which the stabiliser of a codeword is irreducible on $V$ and contained in $G$. 
Recall from Proposition~\ref{c3permisos}(b) that $G=\gamsp(W,\widetilde{B})$ acts transitively on $\mathcal{Q}^\epsilon=\mathcal{Q}^\epsilon(V)$, for  each $\epsilon \in \{ +,- \}$.


\begin{proposition}
	\label{suborthorbs}
	Using the notation in Remark~\ref{remc3setup}, we choose $\Phi
\in\mathcal{Q}^\epsilon(W)$, and set $\phi
:=T(\Phi
)\in\mathcal{Q}^\epsilon(V)$. Let $a\in V$ such that $\Phi(a)\in\ker(\Tr)$ so $\phi_a=T(\Phi_a)\in\mathcal{Q}^\epsilon(V)$, with stabiliser $G_{\phi_a}=\gamoo^\epsilon(W,\Phi_a)$. 
	\begin{enumerate}[(a)]
\item	Let  $g\in G$  such that $\Phi^g=\Phi_a$. Then, for each $\lambda\in\mathbb{F}_{2^b}$, $(\theta_{[\lambda]})^g=\theta_{[\lambda]}(a)$ (these sets are defined in Remark~\ref{remc3setup}(e)). 
	    \item The orbits of $G_{\phi_a}$ in $\mathcal{Q}(V)\setminus\{\phi_a\}$ are the sets $\theta_{[\lambda]}(a)$, for $\lambda\in V$.
	    
	    \item For $a=0$, 
	    the subset $\theta_{[\lambda]}=\theta_{[\lambda]}(0)\subseteq \mathcal{Q}^\epsilon(V)$ if and only if $\lambda\in\ker(\Tr)$, and for such $\lambda$,
	    \[
	    |\theta_{[\lambda]}| = \left\{\begin{array}{ll}
        (2^{n-b}+\varepsilon)(2^n-\varepsilon) & \text{if } \lambda=0,\\
        2^{n-b}(2^n-\varepsilon) & \text{if } \lambda=1 \text{ and } b=2, \\
        b\, 2^{n-b}(2^n-\varepsilon) & \text{otherwise. }
        \end{array}\right.
	    \]
	    Moreover, the largest $G_{\phi}$-orbit $\mathcal{O}$ in $\mathcal{Q}^\epsilon(V)$ satisfies $|\mathcal{O}|\geq |\mathcal{Q}^\epsilon(V)|/2$ if and only if $b=2$ or $3$ and one of the lines of Table~\ref{tabc3suborbs} holds.
	\end{enumerate}
\end{proposition}

\begin{table}[h!]
\centering
 \begin{tabular}{c c l l l} 
 \toprule
 $b$ & $\epsilon$ & $\mathcal{O}$ & $|\mathcal{O}|$ & $\mathcal{Q}^\epsilon\setminus(\mathcal{O}\cup\{\phi\})$ \\ 
 \midrule
 2 & $+$  & $\theta_{[0]}$ & $2^{2n-2} +3\cdot2^{n-2}-1$& $\theta_{[1]}$\\ [0.2ex] 
 2 & $-$  & $\theta_{[1]}$ & $2^{2n-2} +2^{n-2}$& $\theta_{[0]}$, empty if $n=2$\\[0.2ex] 
 3 & $\pm$& $\theta_{[\lambda]}$ ($\lambda\ne 0,1$) & $3\cdot2^{2n-3} -3\epsilon\, 2^{n-3}$ & $\theta_{[0]}$, empty if $(n,\epsilon)=(3,-)$  \\ 
 \bottomrule
 \end{tabular}
 \caption{Table for Proposition~\ref{suborthorbs}(b): largest $G_{\phi}$-orbit $\mathcal{O}$}
\label{tabc3suborbs}
\end{table}

\begin{proof}\leavevmode
\begin{enumerate}[(a)]
\item 
Let $\lambda\in\mathbb{F}_{2^b}$. Since $\Phi_a(c+a)=\Phi^g(c+a)=\Phi((c+a)g^{-1})^\sigma$ (for all $c$) we have
\begin{align*}
    \phi_c\in \theta_{[\lambda]}(a) 
    \iff& c\neq a \text{ and }  \Phi((c+a)g^{-1})\in [\lambda].
\end{align*}
Next, if we substitute $d$ for $(c+a)g^{-1}$, the above expression is equivalent to the following:
\begin{align*}
    \phi_c\in \theta_{[\lambda]}(a) \iff& d=(c+a)g^{-1}\neq 0, \text{ and }  \Phi(d)\in [\lambda].
\end{align*}
Alternatively, we could write $c=dg+a$ and $\phi_d\in\theta_{[\lambda]}$.
This then makes it apparent that $\phi_c\in \theta_{[\lambda]}(a)$ if and only if $\phi_c\in (\theta_{[\lambda]})^g$, because
\(\phi_c=\phi_{dg+a}=\phi_d^g\), by Corollary~\ref{affineactioncor}.
Therefore, $(\theta_{[\lambda]})^g=\theta_{[\lambda]}(a)$.

%


\item Recall that, for $g\in G=\gamsp(W,\widetilde{B})$ relative to an automorphism $\sigma$,
$\Phi^g(v) = \Phi (vg^{-1})^\sigma$ by \eqref{generalisedaction} and $\phi^g=T(\Phi^g)$ by Lemma \ref{c3permisos}.
If $g\in G$ is such that $\Phi^g=\Phi_a$ (which implies that $\phi^g=\phi_a$), recall that, for each $d\in V$, $\phi_d^g=\phi_{dg+a}$ by Corollary \ref{affineactioncor}.

\medskip
By part (a), it is sufficient to prove all the assertions for the case $a=0$ since $G_{\phi_a}$ is conjugate to $G_{\phi_0}=G_\phi$ in $G$. Thus from now on we assume that $a=0$, and we show that the orbits of $G_{\phi}$ in $\mathcal{Q}(V)\setminus\{\phi\}$ are the sets $\theta_{[\lambda]}$.

Let $g\in G_\phi\leq X_\phi$, that is, $g$ is a semi-isometry of $\Phi$ relative to some $\sigma\in\Aut(\mathbb{F}_{2^b})$. 
Consider $\phi_c$ with $c\neq 0$, and let  $\lambda=\Phi(c)$. Then  $\phi_c\in \theta_{\lambda}\subseteq \theta_{[\lambda]}$.  By Lemma \ref{formvectorpermiso}, $\phi_c^g=\phi_{cg}$, and by \eqref{generalisedaction}, $\Phi(cg)=(\Phi^{g^{-1}}(c))^\sigma = \Phi(c)^{\sigma}=\lambda^{\sigma}$, so $\phi_c^g\in \theta_{\lambda^{\sigma}}$. Hence $g$ maps
$\theta_{\lambda}$ to  $\theta_{\lambda^{\sigma}}$. 
Thus the $G_\phi$-orbit of $\phi_c$ is contained in $\theta_{[\lambda]}$.

Since $G_{\phi}$ induces $\Aut(\mathbb{F}_{2^b})$ on $[\lambda]$, it follows that $G_{\phi}$ permutes transitively the subsets $\theta_{\lambda'}$ for $\lambda'\in [\lambda]$. Hence to prove that $G_\phi$ is transitive on $\theta_{[\lambda]}$ it is sufficient to prove that $H_{\phi}=\oo(W,\Phi)$ (which consists of all semi-isometries of $\Phi$ relative to the trivial automorphism) is transitive on $\theta_{\lambda}$. Let $\phi_c, \phi_d\in \theta_{\lambda}$, so $\Phi(c)=\Phi(d)=\lambda$, and in particular,  $c, d$ are both $\Phi$-singular (if $\lambda=0$) or both $\Phi$-nonsingular (if $\lambda\ne 0$). Now $H_{\phi}$ is transitive on both the set of nonzero $\Phi$-singular vectors and  the set of vectors $v$ such that $\Phi(v)$ equals a given nonzero $\lambda$. To see the latter, note that $H_{\phi}$ is transitive on nonsingular 1-subspaces, so there exists $h\in H_{\phi}$ such that $ch$ is in the $1$-space spanned by $d$. Now $\Phi(ch)=\Phi(c)=\lambda$, so by Lemma \ref{uniqueon1space}, $ch=d$.
Hence $\phi_c^h=\phi_{ch}=\phi_d$,
so $H_{\phi}$ is transitive on $\theta_{\lambda}$, and part (a) is proved. 

\item We continue with $a=0$, and since $\phi_c\in\mathcal{Q}^\epsilon(V)$ if and only if $\Phi(c)\in\ker(\Tr)$, the $G_{\phi}$-orbit $\theta_{[\lambda]}\subseteq \mathcal{Q}^\epsilon(V)$ if and only if $\lambda\in\ker(\Tr)$. 
Suppose first that $\lambda=0$ so $[\lambda]=\{0\}$ and $\theta_{[\lambda]}=\theta_0=\{ \phi_c  \mid c\ne 0\ \mbox{and}\ \Phi(c) = 0\}$. Hence $|\theta_{[\lambda]}|$ is equal to the number of nonzero $\Phi$-singular vectors, and hence is equal to $(q^{m-1}+\varepsilon)(q^m-\varepsilon) = (2^{n-b}+\varepsilon)(2^n-\varepsilon)$ by \cite[Theorem 1.14]{ht}. Next, we note that $1\in\ker(\Tr)$ if and only if the prime $b=2$. By Lemma \ref{uniqueon1space}, for each nonzero $\lambda\in\mathbb{F}_{2^b}$, each nonsingular $1$-subspace of $W$ contains exactly one vector $v$ with $\Phi(v)=\lambda$. Hence for each $\lambda\ne 0$, the cardinality $|\theta_{\lambda}|$ is equal to the number of nonsingular $1$-subspaces of $W$, namely, 
	\[
		|\theta_\lambda| = \frac{q^{2m} - (q^{m-1}+\varepsilon)(q^m-\varepsilon) - 1}{q-1}
		                 = q^{m-1}(q^m-\varepsilon) = 2^{n-b}(2^n-\epsilon).
	\]
Thus $|\theta_{[1]}|=|\theta_{1}|= 2^{n-b}(2^n-\epsilon)$, and if $\lambda\ne 0, 1$, then $|[\lambda]|=b$ and hence
$|\theta_{[\lambda]}|=b\,2^{n-b}(2^n-\epsilon)$. 

Finally suppose that $\mathcal{O}= \theta_{[\lambda]}\subseteq \mathcal{Q}^\epsilon(V)$ and $|\mathcal{O}|\geq |\mathcal{Q}^\epsilon(V)|/2$. If $b=2$ then the only $G_{\phi}$-orbits in $\mathcal{Q}^\epsilon(V)$ are $\theta_{[0]}$ and $\theta_{[1]}$, and as we have just shown, their sizes are $2^{2n-2} +3\epsilon\,2^{n-2}-1$ and $2^{2n-2} -\epsilon\,2^{n-2}$, respectively. The larger of these depends on $\epsilon$ and is as given in the first two rows of Table~\ref{tabc3suborbs}; information in the final column is also valid, where we note that, if $(n,\epsilon)=(2,-)$, then $\theta_{[0]}$ is the empty set. Next suppose that $b=3$. Then there are exactly two $G_{\phi}$-orbits in $\mathcal{Q}^\epsilon(V)$, namely $\theta_{[0]}$ and $\theta_{[\lambda]}$, where $[\lambda]$ is the unique triple of elements of $\mathbb{F}_8$ with trace zero (except that $\theta_{[0]}$ is empty if $(n,\epsilon)=(3,-)$). We have shown that these orbits have lengths $2^{2n-3} +7\epsilon\,2^{n-3}-1$ and $3\cdot 2^{2n-3} -3\epsilon\,2^{n-3}$, respectively, and the latter is larger in all cases, so the third row of Table~\ref{tabc3suborbs} is valid. Thus we may assume that the prime $b\geq 5$. Then it is easy to check that the maximum orbit size $|\mathcal{O}|$ is $b\,2^{n-b}(2^n-\epsilon)$, and that the number of orbits of this length in $\mathcal{Q}^\epsilon(V)$ is $(|\ker(\Tr)|-1)/b=(2^{b-1}-1)/b \geq 3$. Thus 
$|\mathcal{Q}^\epsilon(V)| > 3|\mathcal{O}|$, which is a contradiction. This completes the proof.\qedhere
\end{enumerate}
\end{proof}

We complete the $\mathcal{C}_3$-analysis by showing that no $X$-strongly incidence-transitive codes have codeword stabilisers contained in a maximal $\mathcal{C}_3$-subgroup.

\begin{theorem}	\label{c3conclusion}
Let $X=\syp_{2n}(2)$ for $n\geq 2$, let $\Gamma$ be an $X$-strongly incidence-transitive code in $J(\mathcal{Q}^\varepsilon,k)$  for some $k$ satisfying $2 \leq k \leq |\mathcal{Q}^\varepsilon|-2$, such that the stabiliser $X_\Delta$ of a codeword $\Delta\in\Gamma$ is irreducible on $V$. Then $X_\Delta$ is not contained in a maximal $\mathcal{C}_3$-subgroup of $X$. 
\end{theorem}

\begin{proof}
Suppose to the contrary that $X_\Delta\leq G<X$ for some maximal $\mathcal{C}_3$-subgroup $G=\gamsp(W,\widetilde{B})\cong \syp_{2m}(2^b)\rtimes C_b$, for some prime $b$, as in Remark~\ref{remc3setup}. Since $G$ is transitive on $\mathcal{Q}^\epsilon(V)$, $X_\Delta$ is a proper subgroup of $G$. 
By Lemma \ref{smallinequality}, $n \geq 3$ and $2n+1\leq k \leq |\mathcal{Q}^\varepsilon|-(2n+1)$. Note also that $n=mb\geq b$.
Now, since $X_\Delta=X_{\overline{\Delta}}$, the complementary design $\overline{\Gamma}$ satisfies exactly the same conditions. So, interchanging the two codes if necessary, we may assume that $2n+1\leq k \leq |\mathcal{Q}^\varepsilon|/2$.

As in Remark~\ref{remc3setup}, we choose $\Phi$ such that $T(\Phi)=\phi\in\Delta$.  
By Lemma~\ref{choosecarefully}, $\overline{\Delta}$ is contained in a 
$G_{\phi}$-orbit $\mathcal{O}$, so $|\mathcal{O}| \geq |\overline{\Delta}|\geq |\mathcal{Q}^\epsilon|/2$. It then follows from Proposition~\ref{suborthorbs} that $b$ is 2 or 3 and $\mathcal{O}$ is as in one of the rows of Table~\ref{tabc3suborbs}, so $\mathcal{O}= \theta_{[\lambda]}$ for some $\lambda\in\ker(\Tr)$ as specified in the table. For each $\phi_a\in\mathcal{Q}^\epsilon(V)$, $\theta_{[\lambda]}(a)$ is the unique $G_{\phi_a}$-orbit in $\mathcal{Q}^\epsilon(V)$ of length at least $ |\mathcal{Q}^\epsilon(V)|/2$ by Proposition \ref{suborthorbs}(a), since $\theta_{[\lambda]}(a)=(\theta_{[\lambda]})^g$ for some $g$ such that $(G_\phi)^g=G_{\phi_a}$.  Thus, again by Lemma~\ref{choosecarefully}, 
$\overline{\Delta} \subseteq \cap_{\phi_a \in \Delta} \theta_{[\lambda]}(a)$. 
We treat each of the rows of  Table~\ref{tabc3suborbs} separately, but first we consider the exceptional case  $(n,b,\epsilon)=(3,3,-)$  where $\mathcal{Q}^\epsilon(V)=\mathcal{O}\cup \{\phi\}$.  Note that we do not need to consider the case $(n,b,\epsilon)=(2,2,-)$ since $n\geq 3$. 


\medskip\noindent{\it Case 1: $(n,b,\epsilon)=(3,3,-)$.}\quad 
Here $|\mathcal{Q}^-| = 28$ and $G=\syp_2(8)\rtimes C_3$, and we have $7=2n+1 \leq |\Delta| \leq 14$. Since $X_\Delta$ acts transitively on $\Delta\times \overline{\Delta}$ and  $X_\Delta\leq G$, the cardinality  $|G|=1512 = 2^3\cdot 3^3\cdot 7$ is divisible by $k(28-k)$, which is not the case for any $k$ such that $7\leq k\leq 14$.
Thus there are no examples in this case.

\medskip
To treat the other cases it is useful to choose an explicit quadratic form $\Phi=\Phi_0^\epsilon$. We take a symplectic basis $\mathscr{B}=(e_1,\dots, e_m, f_1,\dots, f_m)$ for $W$ relative to $\widetilde{B}$, and, writing a typical vector in $W$ as $x= \sum_{i=1}^m(x_i e_i + y_i f_i)$, with the $x_i, y_i\in\mathbb{F}_{2^b}$, we define $\Phi_0^\epsilon \in \mathcal{Q}^\epsilon(W)$ by
\begin{equation}\label{Phi0overq}
		\begin{array}{ll}
		\Phi_0^+(x) =	\sum_{i=1}^m x_i y_i   \\
		\Phi_0^-(x) =x_m^2 + y_m^2\zeta + \sum_{i=1}^m x_i y_i
		\end{array}
\end{equation}
where $\zeta\in\mathbb{F}_{2^{b}}$ and $t^2+t+\zeta$ is irreducible over $\mathbb{F}_{2^b}$ (see, for example, \cite[Lemma 2.5.2 and Proposition 2.5.3]{kl}).

\medskip\noindent{\it Case 2: Row $1$ of Table~$\ref{tabc3suborbs}$.}\quad 
Here $\epsilon=+$, $b=2$, $\overline{\Delta} \subseteq \cap_{\phi_a \in \Delta} \theta_{[0]}(a)$ and $\Delta$ contains $\theta_{[1]}\cup\{\phi_0\}$. We take $\Phi=\Phi_0^+$. Let $\alpha$ denote a primitive element of $\mathbb{F}_4$.
Then for each $a$ in the following set $A = \lbrace e_m + f_m, \alpha e_m + \alpha^2 f_m, \alpha^2 e_m + \alpha f_m \rbrace$, we have $\Phi_0^+(a) = 1$, and hence $\phi_a\in\Delta$. Suppose that $\phi_c\in\overline{\Delta}$. Then 
$\phi_c\in\cap_{a\in A\cup\{0\}}  \theta_{[0]}(a)$. By \eqref{e:thetalam},
$\Phi_0^+(c)=0$ and $\Phi_a(c+a)=0$ for each $a\in A$. Using equation \eqref{e:phi_a(c+a)} we have, for each $a\in A$,
	\begin{align*}
		\Phi_a(c+a) & = \Phi_0^+(c) + \Phi_0^+(a) + \widetilde{B}(c,a) + \widetilde{B}(c,a)^2 \\
		            & = 1 + \Tr(\widetilde{B}(c,a)) ).
	\end{align*}
	Therefore we require $\Tr(\widetilde{B}(c,a)) = 1$, that is, $\widetilde{B}(c,a) \in \mathbb{F}_4 \setminus \mathbb{F}_2$, for each $a\in A$. Writing  $c=\sum_{i=1}^m(c_ie_i+d_if_i)$, this implies

	\begin{enumerate}[(i)]
		\item $\widetilde{B}(c,e_m + f_m) = c_m + d_m \in \mathbb{F}_4 \setminus \mathbb{F}_2$,
		\item $\widetilde{B}(c,\alpha e_m + \alpha^2 f_m) = \alpha^2 c_m + \alpha d_m \in \mathbb{F}_4 \setminus \mathbb{F}_2$ , and
		\item $\widetilde{B}(c,\alpha^2 e_m + \alpha f_m) = \alpha c_m + \alpha^2 d_m \in \mathbb{F}_4 \setminus \mathbb{F}_2$.
	\end{enumerate}

	Now the sum of any two (possibly equal) elements of $\mathbb{F}_4\setminus \mathbb{F}_2=\{\alpha, \alpha+1\}$ lies in $\mathbb{F}_2$. Thus, adding requirements (i) and (ii) yields $c_m + d_m + \alpha^2 c_m + \alpha d_m = \alpha c_m + \alpha^2 d_m \in \mathbb{F}_2$, which clearly contradicts requirement (iii). Thus there are no possibilities for $\phi_c$, and we have a contradiction.

\medskip\noindent{\it Case 3: Row $2$ of Table~$\ref{tabc3suborbs}$.}\quad 
Here $\epsilon=-$, $b=2$, $n=2m\geq 4$, $\overline{\Delta} \subseteq \cap_{\phi_a \in \Delta} \theta_{[1]}(a)$ and $\Delta$ contains $\theta_{[0]}\cup\{\phi_0\}$, that is, $\Delta$ contains all forms $\phi_c$ such that $c$ is $\Phi$-singular. We take $\Phi=\Phi_0^-$ with $\zeta=\alpha$,  a primitive element of $\mathbb{F}_4$.  Suppose that $\phi_c\in\overline{\Delta}$. Then 
$\phi_c\in\cap_{a\in \sing(\Phi_0^-) }  \theta_{[1]}(a)$. Thus, by \eqref{e:thetalam}, $\Phi_0^-(c)=1$, and for each nonzero $a$ such that $\Phi_0^-(a)=0$, we have $\Phi_a(c+a)=1$.
 Using equation \eqref{e:phi_a(c+a)}, we have, for each nonzero $a$ with $\Phi_0^-(a)=0$,
	\begin{align*}
		\Phi_a(c+a)    & = \Phi_0^-(c) + \Phi_0^-(a) + \widetilde{B}(c,a) + \widetilde{B}(c,a))^2 \\
		            		            & = 1 + \Tr(\widetilde{B}(c,a)).
	\end{align*}
	Therefore we require $\Tr(\widetilde{B}(c,a)) = 0$, that is, $\widetilde{B}(c,a) \in \mathbb{F}_2$, for each nonzero $a$ with $\Phi_0^-(a)=0$. However, for any such $a$ we also have $\Phi_0^-(\alpha a) = \alpha^2 \Phi_0^-(a) = 0$, and hence we require $\widetilde{B}(c,\alpha a)=\alpha \widetilde{B}(c,a)  \in \mathbb{F}_2$. 
	This implies that $\widetilde{B}(c,a) = 0$ must hold for all nonzero $a$ with $\Phi_0^-(a)=0$. In particular this must hold for $a=e_i$  for $1\leq i\leq m-1$. Let  $c=\sum_{i=1}^m(c_ie_i+d_if_i)$. Then $d_i=\widetilde{B}(c,e_i)= 0$, for $1\leq i\leq m-1$, and it must also hold for $a=f_i$, and hence $c_i=\widetilde{B}(c,f_i)= 0$, for $1\leq i\leq m-1$. Hence $c=c_me_m+d_mf_m$, and we require $1=\Phi_0^-(c) = c_m^2+d_m^2\alpha +c_md_m$. The solutions to this equation are $(c_m,d_m)\in\{(0,\alpha), (1,0), (1,\alpha^2), (\alpha, \alpha), (\alpha, \alpha^2)\}$, and hence $|\overline{\Delta}|\leq 5$, which contradicts the fact that $|\overline{\Delta}| \geq \frac{1}{2}|\mathcal{Q}^-|=2^{n-2}(2^n-1)\geq 60$ for all $n \geq 4$. Thus there are no possibilities for $\Delta$ in this case.

\medskip
The final case involves computations over the field $\mathbb{F}_8$, and we make some preparatory comments. The polynomial $f(t)=t^3+t+1$ is irreducible over $\mathbb{F}_2$.  Let $\alpha$ denote a primitive element of $\mathbb{F}_8$, so $\alpha^3=\alpha+1$. The three solutions in $\mathbb{F}_8$ to the  equation $f(t)=0$ are the non-zero elements with trace zero, namely $\alpha, \alpha^2$, and $\alpha^4=\alpha +\alpha^2$; and we have $\ker(\Tr)=\{0\}\cup [\alpha]$ of size $4$. Further, the multiplicative group $\mathbb{F}_8^\ast = \langle\alpha\rangle\cong C_7$, and the map $\mathbb{F}_8\to \mathbb{F}_8$ given by $\alpha\to \alpha+\alpha^2$ has image $\ker(\Tr)$. We use the form $\Phi_0^-$ given by \eqref{Phi0overq}, and in the case $\epsilon=-$, the element $\zeta$ can be taken as any of $1,\alpha+1, \alpha^2+1, \alpha^4+1$ (all with trace $1$).
Recall that we have already dealt with the case $(n,b,\epsilon)= (3,3,-)$ in Case~1.

\medskip\noindent{\it Case 4: Row $3$ of Table~$\ref{tabc3suborbs}$.}\quad 
Here $\epsilon=\pm$, $b=3$, $n=3m$ and $(n,\epsilon)\ne (3,-)$.
Then $\mathcal{O}= \theta_{[\lambda]}$ for some $\lambda\in\ker(\Tr)\setminus\{0\}$ and by the preparatory comments, we may take $\lambda=\alpha$. 
Thus $\overline{\Delta} \subseteq \cap_{\phi_a \in \Delta} \theta_{[\alpha]}(a)$ and $\Delta$ contains $\theta_{[0]}\cup\{\phi_0\}$.  Suppose that $\phi_c\in\overline{\Delta}$. Then
for each  $a$ with $\Phi_0^\epsilon(a)=0$, we have $\phi_a\in\Delta$ and $\phi_c\in \theta_{[\alpha]}(a)$, so $\Phi_a(c+a)\in [\alpha]$. In particular $\Phi_0^\epsilon(c)\in [\alpha]$. Since $\theta_{[0]}$ is non-empty, there exists a non-zero $a$ with $\Phi_0^\epsilon(a)=0$. For each  $\mu \in \mathbb{F}_8$, we also have $\Phi_0^\epsilon(\mu a)=0$ and hence also $\Phi_{\mu a}(c+\mu a)\in [\alpha]$. Using  \eqref{e:phi_a(c+a)}, we have
	\begin{align*}
		\Phi_{\mu a}(c+\mu a) 
		  & = \Phi_0^\epsilon(c) + \Phi_0^\epsilon(\mu a ) + \widetilde{B}(c,\mu a ) + \widetilde{B}(c,\mu a)^2 \\
		  & = \Phi_0^\epsilon(c) + \mu \widetilde{B}(c,a) + \mu^2 \widetilde{B}(c,a)^2              
	\end{align*}
		and therefore we obtain the following condition:
	\begin{align}
		\label{annoyingcondition}
		\Phi_0^\epsilon(c) + \mu \widetilde{B}(c,a) + \mu^2 \widetilde{B}(c,a)^2  \in [\alpha]\quad  \text{ for all  } \mu \in \mathbb{F}_8.
	\end{align}
	We claim that $\widetilde{B}(c,a)=0$. Suppose to the contrary that $z=\widetilde{B}(c,a)\ne 0$. 
	Choose $\mu = \Phi_0^\epsilon(c)^2 z^{-1}$. Then  an easy computation gives $\Phi_0^\epsilon(c) + \mu \widetilde{B}(c,a) + \mu^2 \widetilde{B}(c,a)^2=\Tr(\Phi_0^\epsilon(c))=0$, contradicting \eqref{annoyingcondition}. Thus $\widetilde{B}(c,a)=0$ as claimed, and this must hold for all nonzero $a$ such that $\Phi_0^\epsilon(a)=0$. Write $c=\sum_{i=i}^m(c_ie_i+d_if_i)$, and use the form $\Phi_0^\epsilon$  given by \eqref{Phi0overq} (depending on $\epsilon$). In either case, taking $a=e_i$ or $f_i$ with $1\leq i\leq m-1$ implies that $c_i=d_i=0$ for these values of $i$, so $c=c_me_m+d_mf_m$. If $\epsilon=+$, we similarly conclude that $c_m=d_m=0$ so that $c=0$ and $\overline{\Delta}=\{\phi_0\}$, which is a contradiction. Thus $\epsilon=-$ and so $m\geq 2$ and $n\geq 6$. Then the number of possibilities for a nonzero $c=c_me_m+d_mf_m$ is less than $8^2=64$, and hence $64>|\overline{\Delta}|\geq |\mathcal{Q}^-|/2 = 2^{n-2}(2^n-1)\geq 16\cdot 63$, which is a contradiction. Thus there are no possibilities for $\Delta$ in this case. This completes the proof of Theorem~\ref{c3conclusion}.
\end{proof}

\subsection{\texorpdfstring{$\mathcal{C}_8$}{C8}-type codeword stabilisers}
\label{orthogonalsection}

In this subsection we complete the analysis of the geometric subgroups by showing, in Theorem \ref{c8conclusion}, that no $X$-strongly incidence-transitive codes have codeword stabilisers contained in a maximal $\mathcal{C}_8$-subgroup.
Let $V, B, X$ 
be as in Section~\ref{background}, and let $\mathcal{Q}=\mathcal{Q}(V)$, and $\mathcal{Q}^\epsilon$ ($\epsilon=\pm$) be as in Section~\ref{jordansteiner}. By Table~\ref{syp2n2max}, the maximal $\mathcal{C}_8$-subgroups of $X$ are the stabilisers $X_\psi$ of forms $\psi\in \mathcal{Q}$, and $X_\psi\cong \oo^{\epsilon'}_{2n}(2)$, where $\psi$ has type $\epsilon'=\pm$.
We will see in the proof of Theorem \ref{c8conclusion} that, if a codeword stabiliser in an $X$-strongly incidence-transitive code  in $J(\mathcal{Q}^\varepsilon(V),k)$ is contained in $X_\psi$, then $\psi \in\mathcal{Q}^{-\epsilon}$.  

So suppose that $\phi\in \mathcal{Q}^\epsilon$ and $\psi \in\mathcal{Q}^{-\epsilon}$. 
Recall that by Notation \ref{notation} each quadratic form is equal to $\phi_c$ for a unique $c\in V$. Moreover, by Corollary \ref{formvectorbijection}, 
$\psi=\phi_d$ for some $d\not\in\sing(\phi)$, and $\mathcal{Q}^\epsilon = \{\phi_c\,|\, c\in\sing(\phi)\}$. 
We use this notation throughout the subsection. By Lemma~\ref{formvectorpermiso}, the actions of $X_{\phi}$ (and hence of $X_{\psi\phi}$) on $\mathcal{Q}^\epsilon$ and $\sing(\phi)$ are equivalent, and also, by Corollary~\ref{affineactioncor}, if $g\in X$ and $\phi^g=\phi_a$, then for all $c\in V$, $(\phi_c)^g=\phi_{cg+a}$.  Our first crucial step is to determine the $X_{\psi\phi}$-orbits in $\sing(\phi)$ and $\mathcal{Q}^\epsilon$.

\begin{lemma}
	\label{c8vecorbs}
	Let $\phi,\psi, d$ be as above. Then the following hold.
	\begin{enumerate}
	    \item[(a)]  	The $X_{\phi\psi}$-orbits in $\sing(\varphi)$ are $\lbrace 0 \rbrace, S_0$ and $S_1$, where
	\begin{align*}
		S_0 & =  (\sing(\phi)\cap \sing(\psi)) \setminus \{ 0 \},\ \mbox{and}   \\
		S_1 & = \sing(\phi)\setminus \sing(\psi).
	\end{align*}
    \item[(b)] The $X_{\phi\psi}$-orbits in  $\mathcal{Q}^\varepsilon$ are $\lbrace \varphi_0 \rbrace$, $\omega_0$ and $\omega_1$, where
	\begin{align*}
		\omega_0 & = \lbrace \varphi_c \in \mathcal{Q} |\ \phi(c)=\phi_d(c)=0, \ c\ne 0 \rbrace,\ \mbox{and} \\
		\omega_1 & = \lbrace \varphi_c \in \mathcal{Q} |\  \phi(c)=0,\  \phi_d(c)=1 \rbrace.
	\end{align*}
	Moreover, $|\omega_0|=|S_0|=2^{2n-2}-1$, and $|\omega_1|=|S_1|= 2^{n-1}(2^{n-1}+\epsilon)$.
	\item[(c)] Let $\delta\in\{0,1\}$. Then for each $\phi_a\in\omega_{1-\delta}$, there exists $g\in X_\psi$ such that $\phi^g=\phi_a$ and
	\[
	(\omega_\delta)^g= \lbrace \varphi_{b} \in \mathcal{Q} |\ \phi(b)=0, \ \mbox{and}\ B(a+d,b)=1 \rbrace.
	\] 
	Given $\phi, \psi$ (and hence $d$), this set only depends on $a$, and we set $\omega_\delta(a):=(\omega_\delta)^g.$
	\end{enumerate}
\end{lemma}

\begin{proof}
By Lemma~\ref{formvectorpermiso}, the $X_{\phi}$-actions on  $\mathcal{Q}$  and $V$ are equivalent, and by Corollary~\ref{formvectorbijection}, under this equivalence, $\sing(\phi)$ corresponds to $\mathcal{Q}^\epsilon$. In particular, part (b) follows immediately from part (a), so it is sufficient to prove parts (a) and (c). Further, the equivalence of these actions implies that the subgroup $X_{\psi\phi}$ is equal to the stabiliser of the $\phi$-nonsingular vector $d$ in the  $X_{\phi}$-action on $V$, and hence $X_{\phi\psi}$ fixes setwise the sets $T_0, T_1$ of nonzero vectors $x\in\sing(\phi)$ which lie in $\langle d\rangle^\perp$, or do not lie in $\langle d\rangle^\perp$, respectively. By definition, $\sing(\phi)=\{0\}\cup T_0\cup T_1$.
Also, 
\[
x\in \langle d\rangle^\perp\ \Leftrightarrow\ B(x,d)=0 \ \Leftrightarrow\ \phi_d(x)=\phi(x)+B(x,d)^2=\phi(x),
\]
and hence, for $x\in\sing(\phi)\setminus\{0\}$, $x\in \langle d\rangle^\perp\ \Leftrightarrow\ x\in\sing(\psi)$. It follows that $S_0=T_0$ and  $S_1=T_1$. Next we show that, for each $i$, $S_i$ is an $X_{\phi\psi}$-orbit. 

Let $x\notin \langle d\rangle^\perp$ so $B(x,d)=1$. Then also $x+d \notin \langle d\rangle^\perp$, and 
$$\phi(x+d)=\phi(x)+\phi(d)+B(x,d)= \phi(x)+1+1=\phi(x),
$$
so either both or neither of the two elements of the coset 
$x+\langle d\rangle =\{ x,x+d\}$ lie in $S_1$. If $x\in S_1$ then $\phi(x) = \phi(x+d) = 0$ and $B(x,x+d)=1$, so $\{ x,x+d\}$ is what is called  a {\em hyperbolic pair} in the $\phi$-nonsingular $2$-subspace $ \langle x, d\rangle$ of plus type (while if $x\not\in S_1$ then the $\phi$-nonsingular $2$-subspace $ \langle x, d\rangle$ is of minus type). For any $\epsilon$, the group $X_\phi\cong\oo^{\epsilon}_{2n}(2)$ is transitive on the nonsingular $2$-subspaces of $V$ of plus type, and the stabiliser of such a subspace $W= \langle x, d\rangle$ fixes its unique $\phi$-nonsingular vector $d$, that is $X_{\phi,W}\leq X_{\phi, d}=X_{\phi\psi}$. Moreover the symplectic transvection $\tau_d$ (defined by $v\tau_d=v+B(v,d)d$) fixes $\phi$, $d$ and swaps $x$ and $x+d$ so that $\tau_d\in X_{\phi,W}$. 
It follows that $X_{\phi\psi}$ is transitive on the set $S_1$ of $\phi$-singular vectors in $V\setminus\langle d\rangle^\perp$.  

Now let $x\in \langle d\rangle^\perp$ so $B(x,d)=0$, and assume that $x\ne 0, d$. Then also $x+d \in \langle d\rangle^\perp$, and the computation above shows, this time, that $\phi(x+d)=\phi(x)+1$, so exactly one of the two elements of the 
coset $x+\langle d\rangle =\{ x,x+d\}$ lies in $S_0$. Thus $S_0$ consists of exactly one vector from each nontrivial coset of $\langle d\rangle$ in $\langle d\rangle^\perp$. The group $X_{\phi\psi}$ induces a symplectic group $\syp_{2n-2}(2)$ on the quotient 
$\langle d\rangle^\perp/\langle d\rangle$ (see for example \cite[Section 3.2.4(e)]{maximalfactorisations}), and hence acts transitively on the nontrivial cosets. Since $X_{\phi\psi}$ leaves $S_0$ invariant, it follows that $X_{\phi\psi}$ is transitive on $S_0$.  From this description of $S_0$ it is clear that $|S_0|=2^{2n-2}-1$, and hence that $|S_1|= |\mathcal{Q}^\epsilon|-1-|S_0| = 2^{n-1}(2^n+\epsilon)-1-(2^{2n-2}-1) = 2^{n-1}(2^{n-1}+\epsilon)$. This completes the proof of part (a), and hence also of part (b).

For part (c), let $\delta\in\{0,1\}$ and $\phi_a\in\omega_{1-\delta}\subset \mathcal{Q}^\epsilon$, so $\phi(a)=0$ and $\psi(a)=\phi_d(a)=1-\delta$. Since $X_\psi$ is transitive on $\mathcal{Q}^\epsilon$, there exists $g\in X_\psi$ such that $\phi^g=\phi_a$. 
We easily compute:
\begin{align*}
(\omega_\delta)^g 
&=\lbrace \varphi_c^g \in \mathcal{Q} |\ \phi(c)=0, \ c\ne 0 , \ \phi_d(c)=\delta \rbrace \\
  &=\lbrace \varphi_{cg+a} \in \mathcal{Q} |\ \phi(c)=0, \ c\ne 0 , \ \psi(c)=\delta \rbrace \\
    &=\lbrace \varphi_{b} \in \mathcal{Q} |\ \phi((a+b)g^{-1})=0, \ (a+b)g^{-1}\ne 0 , \ \psi((a+b)g^{-1})=\delta \rbrace \\
      &=\lbrace \varphi_{b} \in \mathcal{Q} |\ \phi^g(a+b)=0, \ b\ne a , \ \psi^g(a+b)=\delta \rbrace \\
       &=\lbrace \varphi_{b} \in \mathcal{Q} |\ \phi_a(a+b)=0, \ b\ne a , \ \psi(a+b)=\delta \rbrace \\
           &=\lbrace \varphi_{b} \in \mathcal{Q} |\ \phi_a(a)+\phi_a(b)+B(a,b)=0, \ b\ne a , \ \psi(a)+\psi(b)+B(a,b)=\delta \rbrace \\
            &=\lbrace \varphi_{b} \in \mathcal{Q} |\ \phi(a)+\phi(b)=0, \ b\ne a , \ \psi(b)+B(a,b)=1 \rbrace 
            \\
             &=\lbrace \varphi_{b} \in \mathcal{Q} |\ \phi(b)=0, \ b\ne a , \ \phi(b)+B(d,b)+B(a,b)=1 \rbrace 
             \\
              &=\lbrace \varphi_{b} \in \mathcal{Q} |\ \phi(b)=0, \ b\ne a , \ B(a+d,b)=1 \rbrace 
      \end{align*}
Notice in particular that this set is independent of the choice of $g$ and of $\delta$.
\end{proof}

\begin{theorem}	\label{c8conclusion}
Let $X=\syp_{2n}(2)$ for $n\geq 2$, let $\Gamma$ be an $X$-strongly incidence-transitive code in $J(\mathcal{Q}^\varepsilon(V),k)$, where $2\leq k\leq |\mathcal{Q}^\varepsilon|-2$, such that for $\Delta\in\Gamma$, the stabiliser $X_\Delta$  is irreducible on $V$. Then $X_\Delta$ is not contained  in a maximal $\mathcal{C}_8$-subgroup of $X$. 
\end{theorem}

\begin{proof}
Suppose that $X_\Delta$ is irreducible on $V$ and  $X_\Delta\leq G<X$ for some maximal $\mathcal{C}_8$-subgroup $G$ of $X$, that is,  
$G=X_\psi$, where $\psi\in\mathcal{Q}$.
Now, since $X_\Delta=X_{\overline{\Delta}}$, the complementary design $\overline{\Gamma}$ satisfies exactly the same conditions. So, by interchanging the two codes if necessary, we may assume that $2\leq k \leq |\mathcal{Q}^\varepsilon|/2$.

Since $X_\Delta$ has two orbits in $\mathcal{Q}^\varepsilon$, namely $\Delta$ and $\overline{\Delta}$, and each has size at least $2$, $X_\Delta$ does not fix any quadratic form in $\mathcal{Q}^\epsilon$. It follows that $\psi\in\mathcal{Q}^{-\epsilon}$. 
Choose $\phi\in\Delta$. Then $\psi=\phi_d$ for some  $\phi$-nonsingular $d\in V$ by Corollary \ref{formvectorbijection}. 
By \cite[Table 1, see also Sections 1.1 and 3.2.4(e)]{maximalfactorisations}, $X=X_\phi X_\psi$, and hence $G$ acts transitively on $\mathcal{Q}^\varepsilon$, and so $X_\Delta$ is a proper subgroup of $G$. It also follows that we may choose $d$ to be a fixed $\phi$-nonsingular vector. 

 Since $\Gamma$ is $X$-strongly incidence-transitive, $X_{\Delta, \phi}$ is transitive on $\overline{\Delta}$, and hence $\overline{\Delta}$ (of size at least $|\mathcal{Q}^\epsilon|/2$) is contained in an $X_{\psi\phi}$-orbit in $\mathcal{Q}^\epsilon$. 
Thus $\overline{\Delta}$ is contained in the larger of the two $X_{\phi\psi}$-orbits $\omega_\delta$, $\delta\in\{0,1\}$, in Lemma~\ref{c8vecorbs}(b). 


Considering the sizes of these two orbits given in Lemma~\ref{c8vecorbs}(b),  we observe that
	\[
		\frac{1}{2}|\mathcal{Q}^\varepsilon| - |\omega_0| =  2^{n-2}(2^n+\epsilon) - (2^{2(n-1)}-1) 
		= 1 + \epsilon 2^{n-2},
	\]
and hence the larger orbit is $\omega_1$ if $\epsilon=+$, and $\omega_0$ if $\epsilon=-$. We therefore treat the cases $\epsilon=+$ and $\epsilon=-$ separately.
As discussed above, we may choose our `favourite' form $\phi$ and our `favourite' $\phi$-nonsingular vector $d$ such that $\psi=\phi_d$. For our computations we choose the form $\phi=\phi_0^\epsilon \in \mathcal{Q}^\epsilon$ as in \eqref{eq-forms}, where we represent vectors $x = \sum_{i=1}^n(x_i e_i + y_i f_i) \in V$ in terms of the symplectic basis $B=\{e_1,\dots,e_n, f_1,\dots,f_n\}$.  	

Our strategy in both cases is as follows: we have $\overline{\Delta}\subseteq \omega_\delta$ (for some $\delta$), and hence $\Delta$ contains $\{\phi\}\cup\omega_{1-\delta}$. Thus, for all $\phi_a\in\omega_{1-\delta}$, $\overline{\Delta}\subseteq \omega_\delta(a)$, as given in Lemma~\ref{c8vecorbs}(c), by Lemma \ref{choosecarefully}. We select a `reasonably sized' subset $A$ of $\phi$-singular vectors such that $\{\phi_a\mid a\in A\}\subseteq \omega_{1-\delta}$, and hence such that
\begin{equation}\label{eq-theta}
\overline{\Delta}\subseteq \Omega:=\omega_\delta\cap \left(\cap_{a\in A} \omega_\delta(a)\right).    
\end{equation}
Apart from the exceptional case $(n, \epsilon)=(2,-)$, we are able to find $A$ such that $\Omega$ has size strictly less than $|\mathcal{Q}^\epsilon|/2$, which is a contradiction since $\overline{\Delta}\subseteq \Omega$ and $\overline{\Delta}$ has size at least $|\mathcal{Q}^\epsilon|/2$.

\medskip\noindent
{\it Case $\epsilon=+$}:\quad Here $\overline{\Delta}\subseteq \omega_1$ and $\Delta$ contains $\{\phi\}\cup\omega_0$. We take  $\phi(x)=\varphi_0^+ (x) = \sum_{i=1}^n x_iy_i$ for $x = \sum_{i=1}^n(x_i e_i + y_i f_i) \in V$, $\phi$-nonsingular vector $d = e_n+f_n$, and  subset $A = \lbrace e_i, f_i \mid 1 \leq i \leq n-1\rbrace$. Note that for each $a\in A$ we have 
$\phi_d(a)=\phi(a) + B(d,a)=0$, and hence $\phi_a\in\omega_0$. Suppose that $\phi_c\in \Omega$, with $\Omega$ as in \eqref{eq-theta}. Then $\phi(c)=0, \phi_d(c)=B(d,c)=1$, and $B(a+d,c)=1$ for all $a\in A$, whence $B(a,c)=0$ for all $a\in A$. Setting $c = \sum_{i=1}^n(c_i e_i + d_i f_i)$, these conditions imply that
\[
\sum_{i=1}^n c_id_i=0,\ c_n+d_n=1,\ \mbox{and}\ c_i=d_i=0\ \mbox{for all}\ i=1,\dots,n-1.
\]
Thus $c\in\{e_n, f_n\}$, and we conclude that $\overline{\Delta}$ is contained in a set of size $2$, which is a contradiction since $\overline{\Delta}$ has size at least $|\mathcal{Q}^\epsilon|/2\geq 5$ for $n\geq 2$.

\medskip\noindent
{\it Case $\epsilon=-$}:\quad  Here $\overline{\Delta}\subseteq \omega_0$ and  $\Delta$ contains $\{\phi\}\cup\omega_1$. We take $\phi(x)=\varphi_0^- (x) = x_n+y_n+\sum_{i=1}^n x_iy_i$ for $x = \sum_{i=1}^n(x_i e_i + y_i f_i) \in V$, $\phi$-nonsingular vector $d = e_n+f_n$, and  subset $A = \lbrace e_i + f_i + f_n \mid 1 \leq i \leq n-1\rbrace$. Note that for each $a\in A$ we have $\phi(a)=0$ and $\phi_d(a)=\phi(a) + B(d,a)= 0+1= 1$, and hence $\phi_a\in\omega_1$. Suppose that $\phi_c\in \Omega$, with $\Omega$ as in \eqref{eq-theta}.
Then $\phi(c)=0, \phi_d(c)=B(d,c)=0$, and $B(a+d,c)=1$ for all $a\in A$, whence $B(a,c)=1$ for all $a\in A$. Setting $c = \sum_{i=1}^n(c_i e_i + d_i f_i)$, these conditions imply that
\[
c_n+d_n+\sum_{i=1}^n c_id_i=0,\ c_n+d_n=0,\ \mbox{and}\ c_i+d_i+c_n=1\ \mbox{for all}\ i=1,\dots,n-1.
\]
Thus, for any given $c_1,\dots, c_n$ the remaining coefficients $d_1,\dots, d_n$ are uniquely determined, and hence there are at most $2^n$ choices for $c$. 
Thus 
\[
  2^{n-2}(2^n-1)=\frac{1}{2}|\mathcal{Q}^-|\leq  |\overline{\Delta}|\leq |\Omega|\leq 2^n,
\]
which implies that $n= 2$, which contradicts Lemma \ref{smallinequality}.
\end{proof}

\section{Proof of Theorem~\ref{mainthm}}\label{proofmain}

	

Finally we prove the main theorem. Let $X=\syp_{2n}(2)$ acting on $\mathcal{Q}^\epsilon=\mathcal{Q}^\epsilon(V)$ of degree $|\mathcal{Q}^\epsilon|=2^{n-1}(2^n+\epsilon)$, where $V=\mathbb{F}_2^{2n}$. Suppose that  $\Gamma$ is an $X$-strongly incidence-transitive code in $J(\calq^\epsilon,k)$, where $1< k < |\mathcal{Q}^\varepsilon|-1$. 
Then $|\mathcal{Q}^\varepsilon|\geq 4$, and so $n\geq 2$.
Further, suppose  that $\Delta\in\Gamma$ and $X_\Delta$ is contained in a geometric subgroup of $X$ as in one of the lines of Table~\ref{syp2n2max}. 

If $X_\Delta$ is reducible on $V$, then by Theorems~\ref{nondegencase}, \ \ref{ticodeissit},  and~\ref{almostsit},  $\Gamma$ is one of the codes in Construction~\ref{ndcode} or Construction~\ref{ticode}.
On the other hand, suppose that  $X_\Delta$ is irreducible on $V$, and lies in a maximal geometric subgroup of $X$ of type $\mathcal{C}_2, \mathcal{C}_3$, or $\mathcal{C}_8$ as in Table~\ref{syp2n2max}. In the case of
$\mathcal{C}_2$-subgroups, it follows from Theorem~\ref{c2smallcases}, that again 
$\Gamma$ arises from Construction~\ref{ndcode}. For the other types, there are no examples, by Theorems~\ref{c3conclusion} and~\ref{c8conclusion}.
Thus, Theorem~\ref{mainthm} is proved.


\end{document}